\documentclass[preprint]{imsart}

\RequirePackage{amsthm,amsmath,amsfonts,amssymb}
\RequirePackage[numbers]{natbib}
\RequirePackage[colorlinks,citecolor=blue,urlcolor=blue]{hyperref}
\RequirePackage{graphicx}
\usepackage{cleveref}
\usepackage{subfigure}
\usepackage{romannum}
\usepackage{romanbar}
\usepackage{todonotes}

\usepackage{geometry}
\startlocaldefs
\theoremstyle{plain}
\newtheorem{theorem}{Theorem}
\newtheorem{korollar}[theorem]{Corollary}
\newtheorem{lemma}[theorem]{Lemma}
\theoremstyle{remark}

\newtheorem{assumption}{Assumption}
\newtheorem{remark}[theorem]{Remark}

\newcommand{\norm}[1]{\left\lVert#1\right\rVert}
\def\rmd{\mathrm{d}}

\newcommand{\lb}[1]{{\lfloor #1 \rfloor_h}}
\newcommand{\ub}[1]{{\lceil #1 \rceil_h}}

\endlocaldefs

\begin{document}

\begin{frontmatter}
\title{Convergence of unadjusted Hamiltonian Monte Carlo \\ for mean-field models}
\runauthor{N. Bou-Rabee, K. Schuh}
\runtitle{Unadjusted HMC for mean-field models}

\begin{aug}
\author[A]{\fnms{Nawaf} \snm{Bou-Rabee}\ead[label=e1]{ nawaf.bourabee@rutgers.edu}},
\author[B]{\fnms{Katharina} \snm{Schuh}\ead[label=e2]{Katharina.Schuh@uni-bonn.de}}
\address[A]{Department of Mathematical Sciences \\ Rutgers University Camden \\ 311 N 5th Street \\ Camden, NJ 08102, USA,\\ \printead{e1}}
\address[B]{Institut f\"ur Angewandte Mathematik \\  Universit\"at Bonn \\
 Endenicher Allee 60 \\ 53115 Bonn, Germany\\  \printead{e2}}
\end{aug}

\begin{abstract}
We present dimension-free convergence and discretization error bounds for the unadjusted Hamiltonian Monte Carlo algorithm applied to high-dimensional probability distributions of mean-field type.  These bounds require the discretization step to be sufficiently small, but do not require strong convexity of either the unary or pairwise potential terms present in the mean-field model.  To handle high dimensionality, our proof uses a particlewise coupling that is contractive in a complementary particlewise metric.
\end{abstract}



\end{frontmatter}

\pagenumbering{arabic}
\section{Introduction}
Markov Chain Monte Carlo (MCMC) methods 
are used to sample from a target probability distribution of the form $\mu(\rmd x)\propto \exp(-U(x))\rmd x$.
The simplest methods (e.g., Gibbs and random walk Metropolis) display random walk behavior which slow their convergence to equilibrium. 
This slow convergence motivates  
the Hamiltonian Monte Carlo (HMC) method, first established in \cite{duanekennedypendletonroweth87}, which offers the potential to converge faster, particularly in high dimension
\cite{ neal11,GiCa2011,BePiRoSaSt2013,cheng2018sharp, dunson2020hastings}.

The convergence properties of HMC have received increasing interest.
Ergodicity was proven in \cite{Schuette99,CaLeSt2007,Stoltz2007}. By drift/minorization conditions, geometric ergodicity was demonstrated in \cite{bourabeesanzserna17,livingstonebetancourtbyrnegirolami19, durmusmoulinessaksmann17}. 
 In \cite{BoEbZi2020, mangoubismith17, chenvempala19}, the convergence behavior is analyzed for a strongly convex potential $U$ and explicit bounds on convergence rates  are obtained using a synchronous coupling approach.
In \cite{BoEbZi2020}, contraction bounds were obtained for more general potentials $U$ by developing a coupling tailored to HMC.  
However, these convergence bounds deteriorate in high dimension for mean-field models (see, in particular, \eqref{eq:contractionrate_standardcase} for the precise form of these contraction bounds for high-dimensional mean-field models).
Therefore, a new approach is needed to obtain convergence bounds for non-strongly convex potentials of mean-field type that are \textit{dimension-free}, i.e., independent of the number of particles in the mean-field model. 

Mean-field models play an important role in understanding statistical properties of high-di\-men\-sio\-nal systems. This connection was introduced by Kac in \cite{kac60} as propagation of chaos and has been investigated amongst others in \cite{mckean66,sznitman91,meleard96}, for very recent related work on second-order mean-field Langevin dynamics see
 \cite{GuLiWuZh2019kinetic,GuMo2020}.
A key component in Kac’s program was to establish bounds on relaxation times of many-body dynamical systems that are dimension-free, see Section 1.4 of \cite{MiMo2013} for a fuller discussion.

The behavior of HMC in high-dimensional mean-field models is also relevant, at least conceptually, to molecular dynamics (MD), see \cite{AlTi1987} and \cite{FrSm2002}, or \cite{LeRoSt2010} for a mathematical perspective. 
MD involves the time integration of high-dimensional Hamiltonian dynamics often coupled to a heat or pressure bath \cite{AlTi1987,FrSm2002}. The corresponding process typically admits a stationary distribution. 
Time discretization introduces an error in the numerically sampled stationary distribution. 
In general, one might hope that this discretization error is dimension-free for ergodic averages of measurable functions (“observables”) that are intensive (e.g., energy per particle) as opposed to extensive (e.g., total energy).  A key contribution of this paper is to demonstrate that this is indeed the case for particles with weak mean-field interactions (see \Cref{thm:ergodicaverages} and \Cref{rem:intensive_observables}). 

In this paper, we consider high-dimensional mean-field models, where the potential $U:\mathbb{R}^{d n}\to \mathbb{R}$ is a function of the form \begin{align*}
U(x)=\sum_{i=1}^n \Big( V(x^i)+\frac{\epsilon}{n}\sum_{\substack{j=1\\ j\neq i}}^n W(x^i-x^j)\Big).
\end{align*}
Here, $V:\mathbb{R}^d\to\mathbb{R}$ and $W:\mathbb{R}^{d}\to\mathbb{R}$ are twice differentiable functions, $\epsilon$ is a real constant and $x=(x^1,...,x^n)$ where $x^i\in\mathbb{R}^d$ represents the position of the $i$-th particle. 
Usually, $d$ is a small fixed number that represents the dimension per particle, whereas the number $n$ of particles is large. We call the unary potential $V$ the  \textit{confinement potential per particle} and the pairwise potential $W$ the \textit{interaction potential}. While we focus on mean-field $U$ with pairwise interactions in this paper, our results can be readily extended to potentials $U$ with more general mean-field interactions (see \Cref{rem:mean_field}).

In its simplest form, every step of HMC uses the Hamiltonian dynamics $(q_t(x,v),p_t(x,v))$ of the mean-field particle system with unit masses defined as the solution to the ordinary differential equations
\begin{equation} \label{hamiltoniandynamics} \begin{aligned} & \frac{\rmd}{\rmd t} q^i_t=p^i_t 
\\ & \frac{\rmd}{\rmd t} p^i_t=-\nabla_i U(q_t)=-\nabla V(q^i_t)-\frac{\epsilon}{n}\sum_{\substack{j=1 \\j\neq i}}^n \Big(\nabla W(q^i_t-q^j_t)-\nabla W(q^j_t-q^i_t)\Big),    \end{aligned} \end{equation}
for $i=1,...,n$ with initial value $(q_0,p_0)=(x,v)$. 
The transition step of the Markov chain in $\mathbb{R}^{d n}$ corresponding to HMC is given by
\begin{align*}
\mathbf{X}(x)=q_T(x,\xi),
\end{align*}
where the initial velocity $\xi\sim\mathcal{N}(0,I_{d n})$ is sampled independently per HMC step, and the integration time $T>0$ is a fixed constant, determining 
the duration of the Hamiltonian dynamics per HMC step.
The corresponding Markov chain is known as {\em exact HMC} because it uses the exact Hamiltonian dynamics and therefore, leaves invariant the target measure $\mu$, cf. ~\cite{bourabeesanzserna18}. 

Generally, the choice of the duration $T$ has a large impact on the performance per HMC step. If $T$ is too small, we obtain a highly correlated chain indicative of random walk behavior. Whereas, if $T$ is chosen too large, due to periodicities and near-periodicities, $q_T(x,v)$ can realize U-turns even as the computational cost of the algorithm increases. This issue was observed by Mackenzie in \cite{mackenzie89}, and motivated duration randomization \cite{ neal11, CaLeSt2007, bourabeesanzserna17} and the No-U-Turn sampler \cite{hoffmangelman14}. 
In contraction bounds for HMC, this issue leads to conditions that limit the duration $T$ of the Hamiltonian dynamics, e.g., for $U$ stronlgy convex $LT^2 \le \text{constant}$ where $L$ is the Lipschitz constant of $\nabla U$ \cite{chenvempala19}.  As we discuss more below, non-convexity of $U$ leads to additional restrictions on the duration $T$.

Since the Hamiltonian dynamics cannot be simulated exactly in general, a numerical version of these dynamics comes into play to approximate the exact dynamics, and normally, the velocity Verlet algorithm is used, cf.~\cite{liu01, bourabeesanzserna18}. The numerical version contains an additional parameter, the discretization step $h>0$ satisfying $T\in h\mathbb{Z}$. 
Note that in the numerical version of HMC without adjusting the algorithm by an additional acceptance-rejection step (see e.g. \cite{neal11,bourabeesanzserna18}), the corresponding Markov chain does not exactly preserve the target measure.  This chain is called  {\em unadjusted HMC}. 
In this article we focus on unadjusted HMC because both from the viewpoint of theory and practice the acceptance-rejection step in adjusted HMC may lead to difficulties in high dimension. Indeed, in the product case (when $\epsilon=0$), a dimension-dependent time step size ($h \propto n^{-1/4}$) is needed to ensure that the acceptance rate in adjusted HMC is bounded away from zero as $n \uparrow \infty$, cf. \cite{BePiRoSaSt2013,gupta1988tuning}. Further, as far as we know only a local contraction result for adjusted HMC is known (see \Cref{rem:adjustedHMC_contr}).
We stress that both adjusted and unadjusted HMC are implementable on a computer, whereas exact HMC is not.

The main result of this paper gives dimension-free convergence bounds for unadjusted HMC applied to mean-field models, i.e., bounds that are independent of the number of particles in the mean-field model.  
Our proof is motivated by the coupling approach in \cite{BoEbZi2020}, but with a new `particlewise' coupling and a complementary particlewise metric.  
We now state a simplified version of our main result, which holds in the special case of exact HMC where $h=0$.

We assume that $\nabla V$ and $\nabla W$ are Lipschitz continuous with Lipschitz constants $L$ and $\tilde{L}$, respectively. Further, we assume that $V$ is $K$-strongly convex outside a Euclidean ball of radius $R$, but possibly non-convex inside this ball. Let $\pi(x,dy)$ be the transition kernel of exact HMC, and let $\mathcal{W}_{\ell^1}$ denote the Kantorovich/$L^1$-Wasserstein distance on $\mathbb{R}^{dn}$ based on an $\ell^1$-metric $\ell^1(x,y)=\sum_{i=1}^n|x^i-y^i|$. 
Then for any two probability measures $\eta$ and $\nu$ on $\mathbb{R}^{dn}$, we show that
\begin{align} \label{eq:Wasserstein_contraction}
\mathcal{W}_{\ell^1}(\eta\pi^m,\nu\pi^m)\leq M e^{-cm} \mathcal{W}_{\ell^1}(\eta,\nu).
\end{align}
Here, $M=\exp\Big(\frac{5}{2}\Big(1+\frac{4R}{T}\sqrt{\frac{L+K}{K}}\Big)\Big)$ and the contraction rate $c$ is of the form
\begin{align*}
c= \frac{1}{156}KT^2\exp\Big(-10\frac{R}{T}\sqrt{\frac{L+K}{K}}\Big).
\end{align*}
This bound holds provided the duration $T$ and the interaction parameter $\epsilon$ are sufficiently small, i.e.,  
\begin{align*}
\frac{5}{3} LT^2&\leq \min\Big(\frac{1}{4},\frac{3K}{10 L},\frac{3K}{256\cdot 5\cdot 2^6 L R^2(L+K) }\Big), \text{ and} 
\\ |\epsilon|\tilde{L} &< \min\Big(\frac{K}{6}, \frac{1}{2}\Big(\frac{K}{36\cdot149}\Big)^2\Big(T+8R\sqrt{\frac{L+K}{K}}\Big)^2\exp\Big(-40\frac{R}{T}\sqrt{\frac{L+K}{K}}\Big)\Big).
\end{align*}
Note that both the contraction rate $c$ and the conditions above are dimension-free, i.e., independent of the number $n$ of particles.  A restriction on the strength of interactions $\epsilon$ cannot be avoided because for large values of $\epsilon$ multiple invariant measures and phase transition phenomena can occur, which typically leads to an exponential deterioration in the rate of convergence as the number of particles tends to infinity \cite{oelschlager1984martingale,sznitman91,tugaut2013convergence}. Roughly speaking, the factor $L R^2$ appearing in the condition on $T$ measures the degree of non-convexity of $U$ and excludes the possibility of high energy barriers.
To obtain this result, we first show contraction for a modified Wasserstein distance that is based on a specially designed particlewise metric $\rho$ on $\mathbb{R}^{dn}$, i.e. , $\mathcal{W}_\rho(\eta\pi^m,\nu\pi^m)\leq e^{-cm} \mathcal{W}_\rho(\eta,\nu)$, and by using that $\rho$ is equivalent to $\ell^1$, we obtain \eqref{eq:Wasserstein_contraction}. 
From this result we deduce a quantitative bound 
for the number $m$ of steps required to approximate the target measure $\mu$ up to a given error $\tilde{\epsilon}$, i.e., $\mathcal{W}_{\ell^1}(\eta\pi^m,\mu)\leq\tilde{\epsilon}$. 
This bound may depend logarithmically on the number $n$ of particles through the distance between the initial distribution and the target measure. 
Finally, we show quantitative dimension-free bounds on the bias for ergodic averages of intensive observables of the form $f(x)=\frac{1}{n}\sum_i\hat{f}(x^i)$.

For unadjusted HMC, we show the same contraction result provided the discretization step $h$ is chosen small enough and deduce that there exists a unique invariant measure $\mu_h$ of unadjusted HMC.
Since unadjusted HMC does not exactly preserve the target measure $\mu$, 
we prove that $\mathcal{W}_{\ell^1}(\mu,\mu_h)=\mathcal{O}(h^2n)$ provided enough regularity for $U$ is assumed, i.e., $V$ and $W$ are three times differentiable and have bounded third derivatives.
If less regularity is assumed, i.e., $V$ and $W$ are only twice differentiable, an $\mathcal{O}(hn)$ bound is obtained. 
Invariant measure accuracy of numerical approximations for related second-order measure preserving dynamics has been extensively investigated in the literature
\cite{RoTw1996B,MaStHi2002,Ta2002,mattingly2010convergence,
BlCaSa2014,bou2010long,leimkuhler2016computation,
abdulle2014high,abdulle2015long}, but according to our knowledge, 
it is new to obtain bounds on $\mathcal{W}_{l^1}$ with a precise dimension dependence 
(see \Cref{cor:bound_mu}).
Durmus and Eberle \cite{DuEb21}, using partially the same approach, generalize these results on invariant measure accuracy to a broader class of both models and inexact (or unadjusted) MCMC methods.

\subsection*{Other work on HMC in high dimension} The study of the behavior of HMC as dimensionality increases is carried out in other settings, too.
For example, in Bayesian inference problems with a large number of observations where the posterior itself is not necessarily high-dimensional. 
In this setting, sampling the posterior directly using HMC is computationally intractable, which motivates stochastic gradient HMC \cite{chen2014stochastic}, the zig-zag process \cite{bierkens2019zig} and the bouncy particle sampler \cite{deligiannidis2019exponential}.
 In \cite{vono2019efficient},
an ADMM-type splitting of the posterior in conjunction with a split Gibbs sampler are proposed, and a dimension-free convergence rate for the split Gibbs sampler is obtained.

Considering the truncation of infinite dimensional probability distributions having a density with respect to a Gaussian reference measure leads to another class of high-dimensional target measures, which arises for instance in path integral MD, cf. \cite{KoRoBoMi2020,bolhuis2002transition,
pinski2010transition}, and statistical inverse problems, cf. \cite{dashti2017bayesian}. 
Dimension-free convergence bounds are obtained for the Metropolis adjusted Langevin Algorithm \cite{eberle2014} and for preconditioned Crank Nicholson (pCN) \cite{hairer2014spectral}. 
Moreover, preconditioned HMC was introduced in \cite{BePiSaSt2011}. The convergence of pHMC was analyzed under strong convexity using a synchronous coupling \cite{pidstrigach20}, and by using a two-scale coupling, dimension-free convergence bounds are obtained for semi-discrete pHMC applied to potential energies that are not necessarily globally strongly convex \cite{BoEb2019}.

Another standard approach to analyze convergence properties in high dimension is optimal scaling of MCMC, see \cite{GeGiRo1997,RoRo1998,BeRoSt2009,
dunson2020hastings}.
This theory of optimal scaling provides a general way to tune the time step size in HMC \cite{gupta1988tuning,BePiRoSaSt2013}.

While our object of study is the simplest version of HMC applied to mean-field models, there are other variants of HMC available including one that uses a general reversible approximation of the Hamiltonian dynamics \cite{fang2014compressible}, HMC with partial randomization of momentum \cite{Ho1991,AkRe2008}, preconditioned HMC using a position dependent mass matrix \cite{GiCa2011}, and adjusted HMC with delayed rejection \cite{CaSa2015}. 

\subsection*{Outline} The rest of the paper is organized as follows. In \Cref{section_setting}, we state the considered framework before presenting our main results in \Cref{section_mainresults}. In \Cref{section_estimates}, estimates used to prove the main results are stated. Finally, \Cref{section_proof_oflemmata}, \Cref{section_proof_oftheorems} and \Cref{section:proofs_errorbounds} contain the proofs.

\section{Preliminaries} \label{section_setting}
We first give the definition of unadjusted HMC applied to mean-field models and state assumptions for the mean-field model before constructing the particlewise coupling used to obtain the contraction result in the next section.

\subsection{Hamiltonian Monte Carlo Method} \label{subsec:exactHMC}
Consider a function $U\in\mathcal{C}^2(\mathbb{R}^{d n})$ of the form
\begin{align} \label{eq:definition_U}
U(x)=\sum_{i=1}^n \Big( V(x^i)+\frac{\epsilon}{n}\sum_{\substack{j=1\\ j\neq i}}^n W(x^i-x^j)\Big)
\end{align}
such that $\int \exp(-U(x))\rmd x<\infty$ holds. Assuming all particles have unit masses, the corresponding Hamiltonian is defined by 
$H(x,v)=U(x)+\frac{1}{2}|v|^2 
$ for $x,v \in\mathbb{R}^{dn}$.
The HMC method is an MCMC method for sampling from a `target' probability distribution
\begin{equation} \label{eq:probabilitymeasure_sampling}
\mu(\rmd x)= Z^{-1} \exp(-U(x))\rmd x, 
\end{equation}
on $\mathbb{R}^{d n}$ with normalizing constant $Z=\int\exp(-U(x))\rmd x$. 
In particular, the HMC method generates a Markov chain on $\mathbb{R}^{d n}$.

Since \eqref{hamiltoniandynamics} is not exactly solvable, a discretized version is considered. Here, we consider the velocity Verlet integrator with discretization step $h>0$, cf. \cite{bourabeesanzserna18}.
The numerical solution produced by the velocity Verlet integrator is interpolated by the flow $({q}_t(x,v),{p}_t(x,v))$ of the ODE
\begin{align}\label{eq:hamdyn_num}
\frac{\rmd }{\rmd t}{q}_t^i={p}^i_{\lfloor t\rfloor_h}-\frac{h}{2}\nabla_i U({q}_{\lfloor t \rfloor_h}), \hspace{1cm} \frac{\rmd }{\rmd t}{p}^i_t=-\frac{1}{2}(\nabla_i U({q}_{\lfloor t\rfloor_h})+\nabla_i U({q}_{\lceil t\rceil _h}))
\end{align}
with initial condition $({q}_0,{p}_0)=(x,v)$ where
\begin{align*}
\lfloor t\rfloor_h=\max\{s\in h\mathbb{Z}:s\leq t\}, \hspace{5mm} \lceil t\rceil_h=\min\{s\in h\mathbb{Z}:s\geq t\},
\end{align*}
and where $\nabla_i U:\mathbb{R}^{dn}\to\mathbb{R}^d$ is the gradient in the $x^i$-th direction, i.e., $\frac{\partial U}{\partial x^i}$.
The transition step of {\em unadjusted HMC} is
given by $x \mapsto \mathbf{X}_h(x)$ where ${\mathbf{X}_h}(x)={q}_T(x,\xi)$, $T/h\in\mathbb{Z}$ for $h>0$ and $\xi\sim \mathcal{N}(0,I_{dn})$ is a random variable, where $\mathcal{N}(0,I_{dn})$ denotes the centered normal distribution on $\mathbb{R}^{dn}$ with covariance given by the $dn\times dn$ identity matrix. 
The transition kernel of the Markov chain on $\mathbb{R}^{d n}$ induced by the unadjusted HMC algorithm is denoted by ${\pi_h}(x,B)=P[{\mathbf{X}_h}(x)\in B]$.
 
If $h>0$ is fixed, we write the abbreviation $\lfloor t\rfloor$ and $\lceil t\rceil$ instead of $\lfloor t\rfloor_h$ and $\lceil t\rceil_h$ and omit the $h$ dependence in $\mathbf{X}_h(x)$. 
For $h=0$ we consider the solution $(q_t(x,\xi),p_t(x,\xi))$ of \eqref{hamiltoniandynamics} and obtain {\em exact HMC} with transition step $\mathbf{X}(x) := \mathbf{X}_0(x)=q_T(x,\xi)$ and transition kernel $\pi(x,B) := \pi_0(x,B)$.
As the Hamiltonian is not preserved by the numerical flow with $h>0$, unadjusted HMC does not preserve the target measure $\mu$. 
Therefore, after we study convergence of unadjusted HMC, we then bound the error between exact and unadjusted HMC in \Cref{section_mainresults}.

\subsection{Mean-field particle model} \label{subsec:mean-field}
Let $U:\mathbb{R}^{d n}\to \mathbb{R}$ be a potential function of the form  \eqref{eq:definition_U} 
where $V:\mathbb{R}^d\to\mathbb{R}$ and $W:\mathbb{R}^{d}\to\mathbb{R}$ are twice continuously differentiable functions such that $\int \exp(-U(x))\mu(\rmd x)<\infty$. Without loss of generality we assume that $\epsilon$ is a non-negative constant. Otherwise we change the sign of the interaction potential $W$.
The following conditions are imposed on the functions $V$ and $W$ for proving the contraction results for exact HMC.
\begin{assumption} \label{ass_V_locmin}
	$V$ has a global minimum at $0$, $V(0)=0$ and $V(\mathsf{x})\geq 0$ for all $\mathsf{x}\in\mathbb{R}^d$.
\end{assumption}

\begin{assumption} \label{ass_V_lipschitz}
	$V$ has bounded second derivatives, i.e., $L:=\sup \|\nabla^2V\|<\infty$.
\end{assumption}

\begin{assumption} \label{ass_V_strongconv}
	$V$ is strongly convex outside a Euclidean ball: 
 There exists $K\in(0,\infty)$ and $R\in[0,\infty)$ such that for all $\mathsf{x},\mathsf{y}\in\mathbb{R}^d$ with $|\mathsf{x}-\mathsf{y}|\geq R$,
\begin{align*}
(\mathsf{x}-\mathsf{y})\cdot(\nabla V(\mathsf{x})-\nabla V(\mathsf{y}))\geq K|\mathsf{x}-\mathsf{y}|^2.
\end{align*}

\end{assumption}
\begin{assumption} \label{ass_W_lipschitz}
		$W$ has bounded second derivatives, i.e., $\tilde{L}:=\sup \|\nabla^2 W\|<\infty$.
\end{assumption}

We note that \Cref{ass_V_locmin} is stated for simplicity, since \Cref{ass_V_strongconv} implies that $V$ has a local minimum and so \Cref{ass_V_locmin} can always be obtained by adjusting the coordinate system appropriately and adding a constant to $V$. 
Since $V$ is a unary confinement potential per particle and $W$ is a pairwise interaction potential, note that the strong convexity constant $K$, the Lip\-schitz constants $L$, $\tilde{L}$ and the radius $R$ are dimension-free, i.e., independent of the number of particles.
By \Cref{ass_V_locmin}, \Cref{ass_V_lipschitz} and \Cref{ass_W_lipschitz},
\begin{align}
&|\nabla V(\mathsf{x})|=|\nabla V(\mathsf{x})-\nabla V(0)|\leq L|\mathsf{x}|, \hspace{20mm} \text{ and}  \label{eq:conseq_ass1}
\\ &|\nabla W(\mathsf{x}-\mathsf{y})-\nabla W(\mathsf{y}-\mathsf{x})|\leq 2\tilde{L}|\mathsf{x}-\mathsf{y}|\leq 2\tilde{L}(|\mathsf{x}|+|\mathsf{y}|) \label{eq:conseq_ass2}
\end{align}
for all $\mathsf{x},\mathsf{y}\in\mathbb{R}^d$. From \eqref{eq:conseq_ass1} and \Cref{ass_V_strongconv}, it follows that $K$ is smaller than $L$,
\begin{align} \label{eq:KandL}
K/L\leq 1.
\end{align}
Further, we deduce from \Cref{ass_V_lipschitz} and \Cref{ass_V_strongconv} that
 for all $\mathsf{x},\mathsf{y}\in\mathbb{R}^{d}$,
\begin{align} \label{eq:cond_strongconv}
(\mathsf{x}-\mathsf{y})\cdot(\nabla V(\mathsf{x})-\nabla V(\mathsf{y}))\geq K|\mathsf{x}-\mathsf{y}|^2-\hat{C}
\end{align} 
with $\hat{C}:=R^2(L+K)$ and so $V$ is asymptotically strongly convex.

\begin{remark} \label{rem:mean_field}
In this work, we focus on a pairwise mean-field interaction energy $W$.  However, the results can be readily extended to the situation where the Hessian of the mean-field potential $U$ satisfies: \begin{align}
    L = \sup_{\substack{1 \le i \le d n \\ x \in \mathbb{R}^{d n}}} \left| \frac{\partial^2 U}{\partial x_i^2}(x) \right| \;, \quad  \tilde{L} = \sup_{\substack{1 \le i < j \le d n \\ x \in \mathbb{R}^{d n}}} \left| \frac{\partial^2 U}{\partial x_i \partial x_j}(x) \right| \;
\end{align} 
and the parameter $\tilde{L}$ scales like $1/n$ as $n \to \infty$ which corresponds to the standard mean-field limit \cite{oelschlager1984martingale,sznitman91,tugaut2013convergence,DuEbGuZi20}.
\end{remark}

For proving discretization error bounds, 
we suppose additionally for the confinement potential $V$ and for the interaction potential $W$:
\begin{assumption} \label{ass_V_thirdfourthder}
	$V$ is three times differentiable and has bounded third derivatives, i.e., $L_H:=\sup\|\nabla^3 V\|<\infty$. 
\end{assumption}
\begin{assumption} \label{ass_W_thirdfourthder}
	$W$ is three times differentiable and has bounded third derivatives, i.e., $\tilde{L}_H=\sup\|\nabla^3 W\|<\infty$. 
\end{assumption} 
This additional regularity gives a better order in the error bounds between exact HMC and unadjusted HMC, see \Cref{thm:flowdifference_exact_num}. 

Possible interaction potentials meeting \Cref{ass_W_lipschitz} and \Cref{ass_W_thirdfourthder} are the Morse potential \cite{Wales2010jp} and the harmonic (or linear) bonding potential \cite[Section 7.4.1.1]{Amarante2017}, which are both used to model interactions between particles in molecular dynamics.

\begin{remark} \label{rem:effectiveLip}
Let us note that by \eqref{eq:conseq_ass2} and \eqref{eq:cond_strongconv} it holds for the potential $U$ that
\begin{align*}
( x&-y)\cdot (\nabla U(x)-\nabla U(y))
=\sum_{i=1}^n\Big( (x^i- y^i)\cdot(\nabla V(x^i)-\nabla V(y^i))
\\ & +\frac{\epsilon}{n}\sum_{j\neq i}(x^i-y^i)\cdot(\nabla W(x^i-x^j)-\nabla W(y^i-y^j)-\nabla W(x^j-x^i)+\nabla W(y^j-y^i))\Big)
\\ & \ge K|x-y|^2-n(K+L)R^2-\frac{2\epsilon\tilde{L}}{n}\sum_{i}\sum_{j\neq i} |x^i-y^i-(x^j-y^j)||x^i-y^i|
\\ & \ge  (K-4\epsilon\tilde{L})|x-y|^2-n(K+L)R^2.
\end{align*}
Hence, the potential $U$ is strongly convex if $R=0$ and $K-4\epsilon\tilde{L}>0$ holds. Moreover, a similar calculation shows that $\nabla U$ is globally Lipschitz continuous with an effective Lipschitz constant of $L+4 \epsilon \tilde L$.  In this case, \cite[Theorem 2.1]{BoEbZi2020} and \cite[Theorem 1]{mangoubismith17} have already shown contraction for exact HMC with the  dimension-free rate $c=(1/2)(K-4\epsilon\tilde{L})T^2$
if $(L+4\epsilon\tilde{L})T^2\leq (K-4\epsilon\tilde{L})/(L+4\epsilon\tilde{L})$ holds. Recently, the latter condition on $T$ has been improved to $(L+4\epsilon\tilde{L})T^2\leq (1/4)$,  cf. \cite[Theorem 3]{chenvempala19}. Whereas, if $R>0$, then the potential $U$ is only asymptotically strongly convex provided $K-4\epsilon\tilde{L}>0$, and in this case,
\begin{align*}
( x-y)\cdot(\nabla U(x)-\nabla U(y)) \ge ((K-4\epsilon\tilde{L})/2)|x-y|^2  
\end{align*}
for all $|x-y|\ge R_n=R\sqrt{2n(L+K)/(K-4\epsilon\tilde{L})}$.
Thus, by \cite[Theorem 2.3]{BoEbZi2020} we obtain the following contraction rate for exact HMC 
\begin{align} \label{eq:contractionrate_standardcase}
c_n=(1/10)\min(1,(1/4)(K-4\epsilon\tilde{L})T^2(1+(R_n/T))e^{-R_n/(2T)})e^{-2R_n/T}
\end{align}
provided $(L+4\epsilon\tilde{L})T^2\leq \min(1/4,(K-4\epsilon\tilde{L})/(L+4\epsilon\tilde{L}),1/(2^6(L+4\epsilon\tilde{L})R_n^2))$ holds.
The condition on $T$ is dependent on the number $n$ of particles and the rate $c_n$ decreases exponentially fast in the number of particles. This dimension dependence motivates the particlewise coupling stated next. 
\end{remark}

\subsection{Construction of coupling} \label{subsec:coupling}
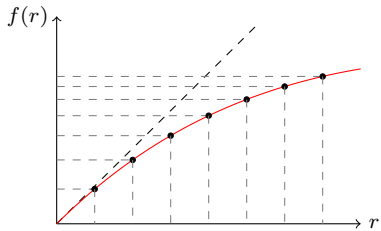
\begin{figure} 
%
%

\begin{subfigure}{
\begin{tikzpicture}
\draw[->] (0,0)--(0,2.75);
\draw[->] (0,0)--(4,0);
\draw[dashed, -] (0,0)--(2.65,2.65);
\draw (0,2.75) node[left]{\footnotesize $f(r)$};
\draw (4,0) node[right]{\footnotesize $r$};
\draw[scale=1, domain=0:4, smooth, variable=\x, red]  plot ({\x},{\x*e^(-(1/6)*\x)});

\draw[dashed,-, gray] (0,{0.5*e^(-(1/6)*0.5)})--(0.5,{0.5*e^(-(1/6)*0.5)});
\draw[dashed,-, gray] (0,{1*e^(-(1/6)*1)})--(1,{1*e^(-(1/6)*1)});
\draw[dashed,-, gray] (0,{1.5*e^(-(1/6)*1.5)})--(1.5,{1.5*e^(-(1/6)*1.5)});
\draw[dashed,-, gray] (0,{2*e^(-(1/6)*2)})--(2,{2*e^(-(1/6)*2)});
\draw[dashed,-, gray] (0,{2.5*e^(-(1/6)*2.5)})--(2.5,{2.5*e^(-(1/6)*2.5)});
\draw[dashed,-, gray] (0,{3*e^(-(1/6)*3)})--(3,{3*e^(-(1/6)*3)});
\draw[dashed,-, gray] (0,{3.5*e^(-(1/6)*3.5)})--(3.5,{3.5*e^(-(1/6)*3.5)});
\filldraw[black] (0.5,{0.5*e^(-(1/6)*0.5)}) circle(1pt);
\filldraw[black] (1,{1*e^(-(1/6)*1)}) circle(1pt);
\filldraw[black] (1.5,{1.5*e^(-(1/6)*1.5)}) circle(1pt);
\filldraw[black] (2,{2*e^(-(1/6)*2)}) circle(1pt);
\filldraw[black] (2.5,{2.5*e^(-(1/6)*2.5)}) circle(1pt);
\filldraw[black] (3,{3*e^(-(1/6)*3)}) circle(1pt);
\filldraw[black] (3.5,{3.5*e^(-(1/6)*3.5)}) circle(1pt);
\draw[dashed,-, gray] (0.5,{0.5*e^(-(1/6)*0.5)})--(0.5,0);
\draw[dashed,-, gray] (1,{1*e^(-(1/6)*1)})--(1,0);
\draw[dashed,-, gray] (1.5,{1.5*e^(-(1/6)*1.5)})--(1.5,0);
\draw[dashed,-, gray] (2,{2*e^(-(1/6)*2)})--(2,0);
\draw[dashed,-, gray] (2.5,{2.5*e^(-(1/6)*2.5)})--(2.5,0);
\draw[dashed,-, gray] (3,{3*e^(-(1/6)*3)})--(3,0);
\draw[dashed,-, gray] (3.5,{3.5*e^(-(1/6)*3.5)})--(3.5,0);

\end{tikzpicture} }
\end{subfigure}
\caption{Under an increasing concave distance function $f$, a decrease in $r$ has a larger impact on $f(r)$ than an increase in $r$, 
i.e., $f(r)-f(r-\Delta)\geq f(r+\Delta)-f(r)$ for $r, \Delta>0$.}\label{fig:concavefunction}
\end{figure}

We establish a coupling between the transition probabilities $\pi_h(x,\cdot)$ and $\pi_h(y,\cdot)$ of unadjusted HMC with discretization step $h$ for two states $x,y\in\mathbb{R}^{d n}$.
The key idea for the coupling is to locally couple the velocity randomizations, i.e., for the $i$-th particles in each component of the coupling separately and independently of the other particles. 
A particlewise coupling approach was used before in \cite{Eberle16A,DuEbGuZi20} and enables us here to show a dimension-free contraction rate, i.e. a rate that does not depend on the number $n$ of particles. 
The idea for the construction for the $i$-th particles in each component of the coupling is adapted from \cite{BoEbZi2020}, see also \cite{eberleguillinzimmer19}.
The coupling transition step for unadjusted HMC is given by
\begin{align} \label{eq:coupling_transitionstep}
\mathbf{X}(x,y)=q_T(x,\xi)  \hspace{5mm}\text{and}\hspace{5mm} \mathbf{Y}(x,y)=q_T(y,\eta) 
\end{align}
with $q_T$ defined in \eqref{eq:hamdyn_num} and where $\xi$ and $\eta$ are the corresponding velocity refreshments for the position $x$ and $y$ given in the following way:
Let $\xi\in\mathbb{R}^{d n}$ be a normally distributed random variable.
Let $\mathcal{U}_i\sim \text{Unif}[0,1]$ be independent uniformly distributed random variables that are independent of $\xi$. Let $\gamma$ be a constant that is specified later. 
If $|x^i-y^i|\geq\tilde{R}$, where $\tilde{R}$ is a positive constant specified later, we apply a synchronous coupling for the $i$-th particle by setting $\eta^i=\xi^i$. If $|x^i-y^i|<\tilde{R}$, the $i$-th velocity refreshment of $y$ is given by 
\begin{equation}
\begin{aligned} \label{eq:def_contractivecoupling}
\eta^i:=\begin{cases} \xi^i+\gamma z^i & \text{if } \mathcal{U}_i\leq \frac{\varphi_{0,1}(e^i\cdot\xi^i+\gamma|z^i|)}{\varphi_{0,1}(e^i\cdot\xi^i)}, \\ \xi^i-2(e^i\cdot \xi^i)e^i &\text{otherwise,} \end{cases}
\end{aligned}
\end{equation}
where $\varphi_{0,1}$ denotes the density of the standard normal distribution, $z^i=x^i-y^i$, and $e^i=z^i/|z^i|$ if $|z^i|\neq 0$. If $|z^i|=0$, $e^i$ is some arbitrary unit vector.
If we consider the free dynamics, i.e.,  $U\equiv 0$, then the first case in \eqref{eq:def_contractivecoupling} leads to a decrease in the difference of the positions in the $i$-th component provided the duration $T$ is sufficiently small, i.e., $|\mathbf{X}^{ i}(x,y)-\mathbf{Y}^{ i}(x,y)|=|x^i-y^i||1-T\gamma|$. 
When $U$ does not vanish, we obtain contractivity of this coupling in a metric equivalent to the standard $\ell_1$ metric that involves a concave distance function, see \Cref{fig:concavefunction}.

We note that each of the components $\eta^i$ are normally distributed random variables by \cite[Section 2.3]{BoEbZi2020} and that the components $\eta^i$ are independent by the  independent particlewise construction. This implies $\eta\sim\mathcal{N}(0,I_{d n})$, which is sufficient to verify that the constructed transition step given by \eqref{eq:coupling_transitionstep} is a coupling of the transition probabilities $\pi_h(x,\cdot)$ and $\pi_h(y,\cdot)$.


\subsection{Numerical simulations}\label{sec:numerical_simulations}

We next present a numerical illustration of some properties of the particlewise coupling which supports the main results for unadjusted HMC stated in the next section. 

\begin{figure}
\begin{subfigure}{
	\includegraphics[scale=0.33]{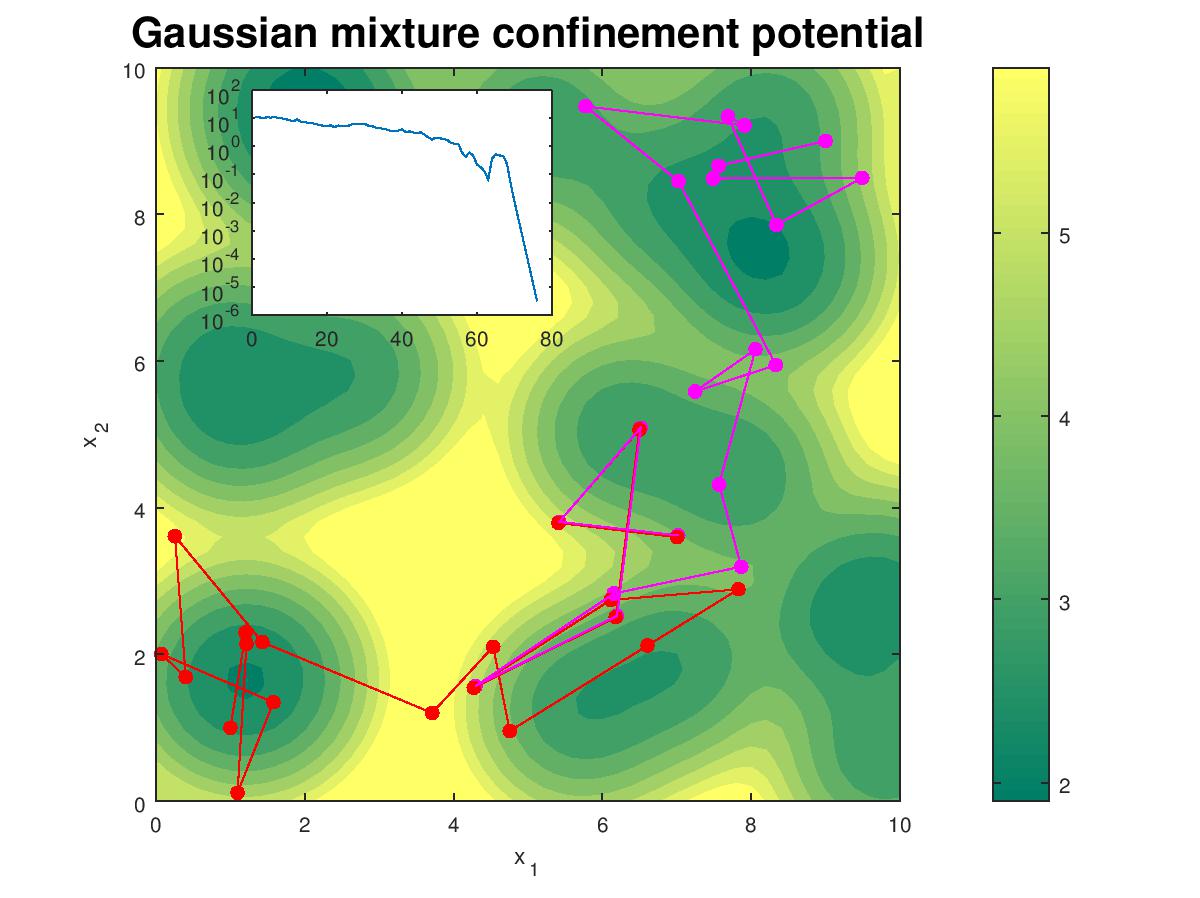}}
\end{subfigure}
\begin{subfigure}{
	\includegraphics[scale=0.33]{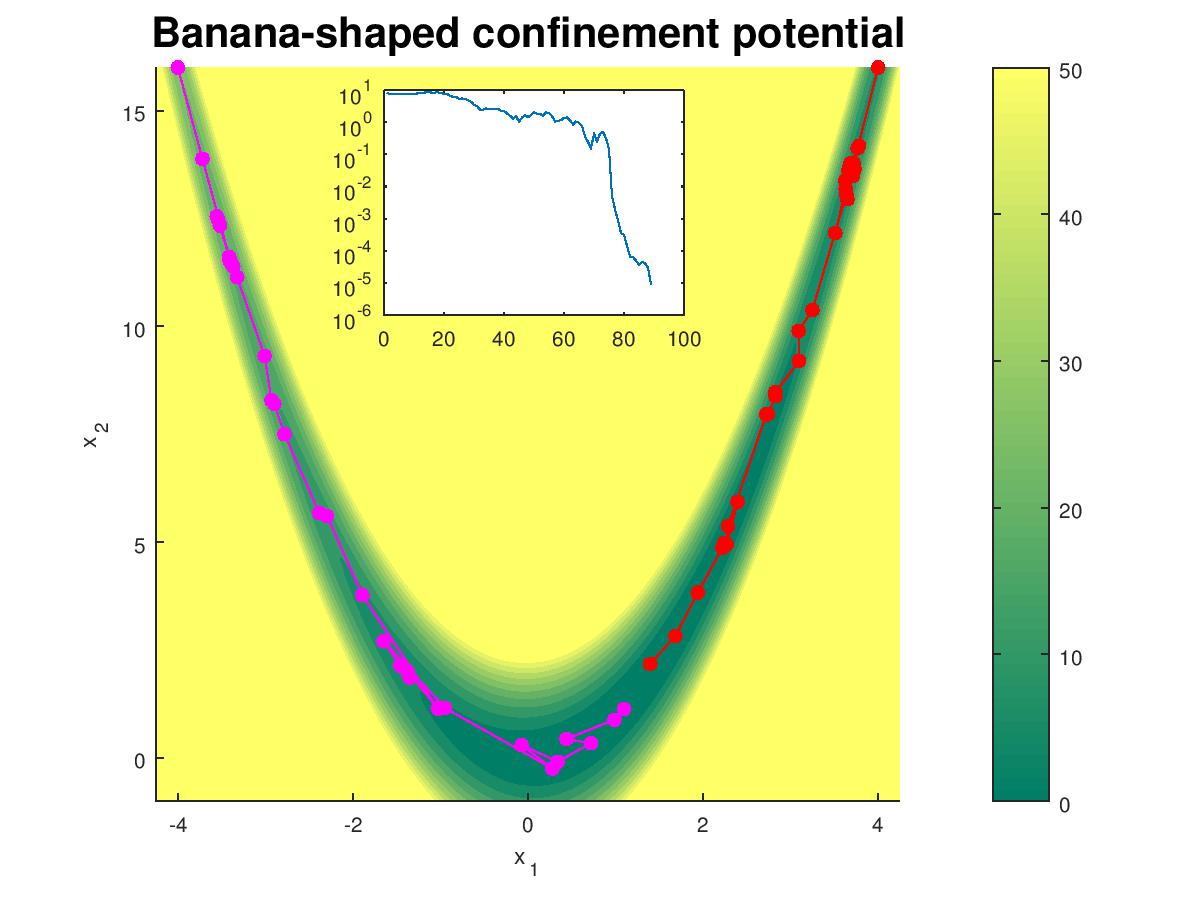}}\end{subfigure}
\caption{\small{\textit{Coupling of HMC applied to mean-field models with $n=10$ particles. The confinement potential is the potential of a Gaussian mixture distribution in the left plot and of a banana-shaped distribution in the right plot. The projection to one particle of the Markov chain is plotted on the contour graph of the potentials and connected by a linear interpolation; the inset shows the mean distance between the two components of the coupling on a log-scale. 
} }}\label{figure}
\end{figure}

\begin{figure}
\begin{subfigure}{
	\includegraphics[scale=0.335]{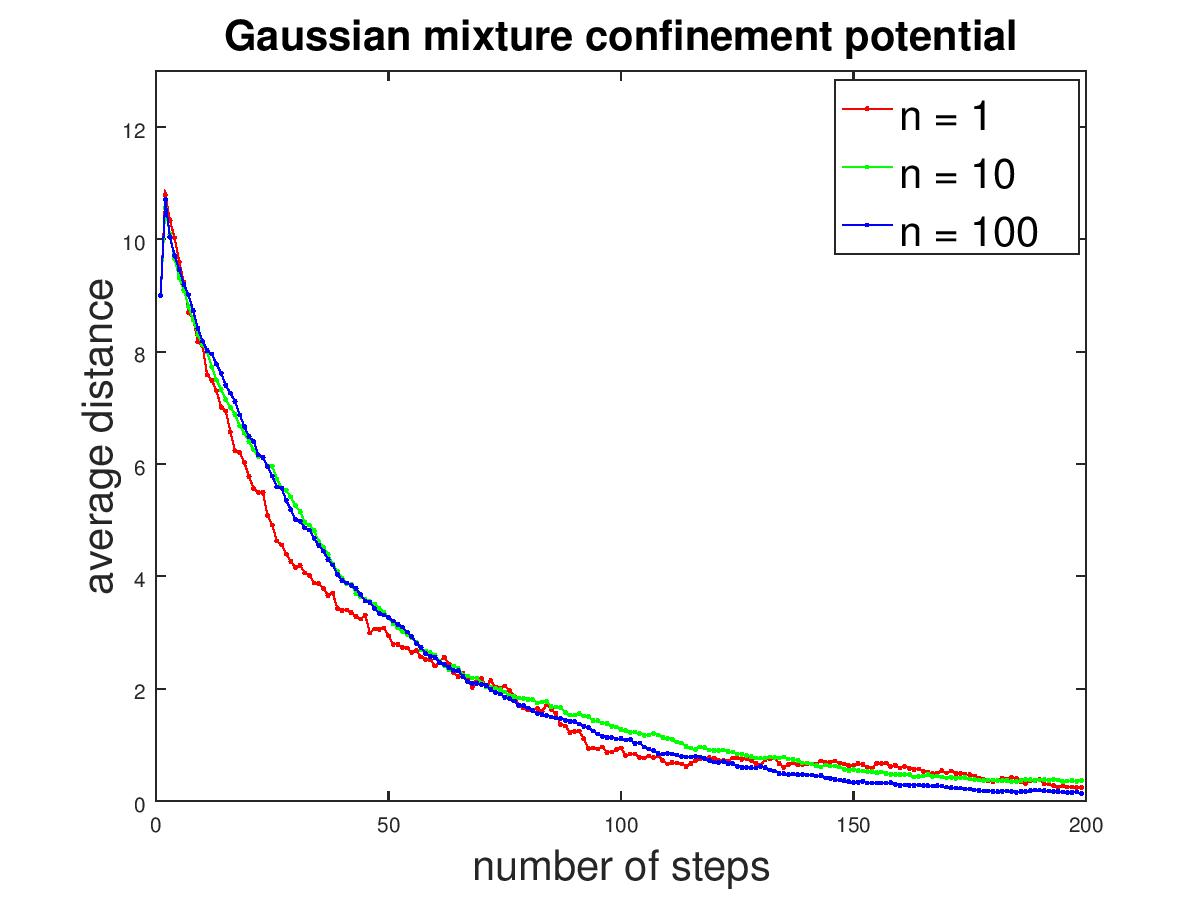} }
\end{subfigure}
\begin{subfigure}{
	\includegraphics[scale=0.335]{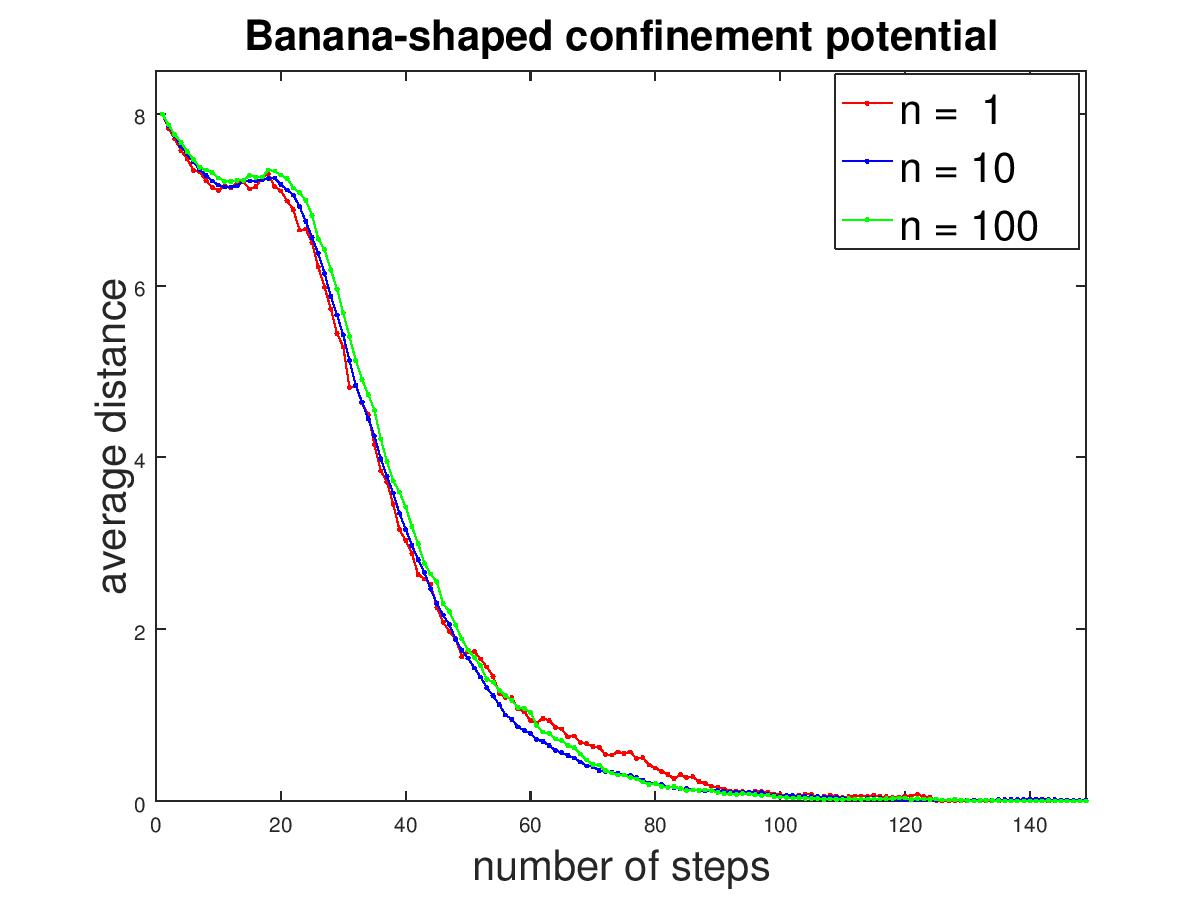} }\end{subfigure}
\caption{\small{\textit{Evolution of the mean distance $\frac{1}{n}\sum_{i=1}^n |\mathbf{X}^i_k-\mathbf{Y}^i_k|$ between the two components of the coupling for HMC after $k$ steps with $n\in\{1,10,100\}$ particles.}}}\label{figure2}
\end{figure}

\begin{figure}
\begin{subfigure}{
	\includegraphics[scale=0.335]{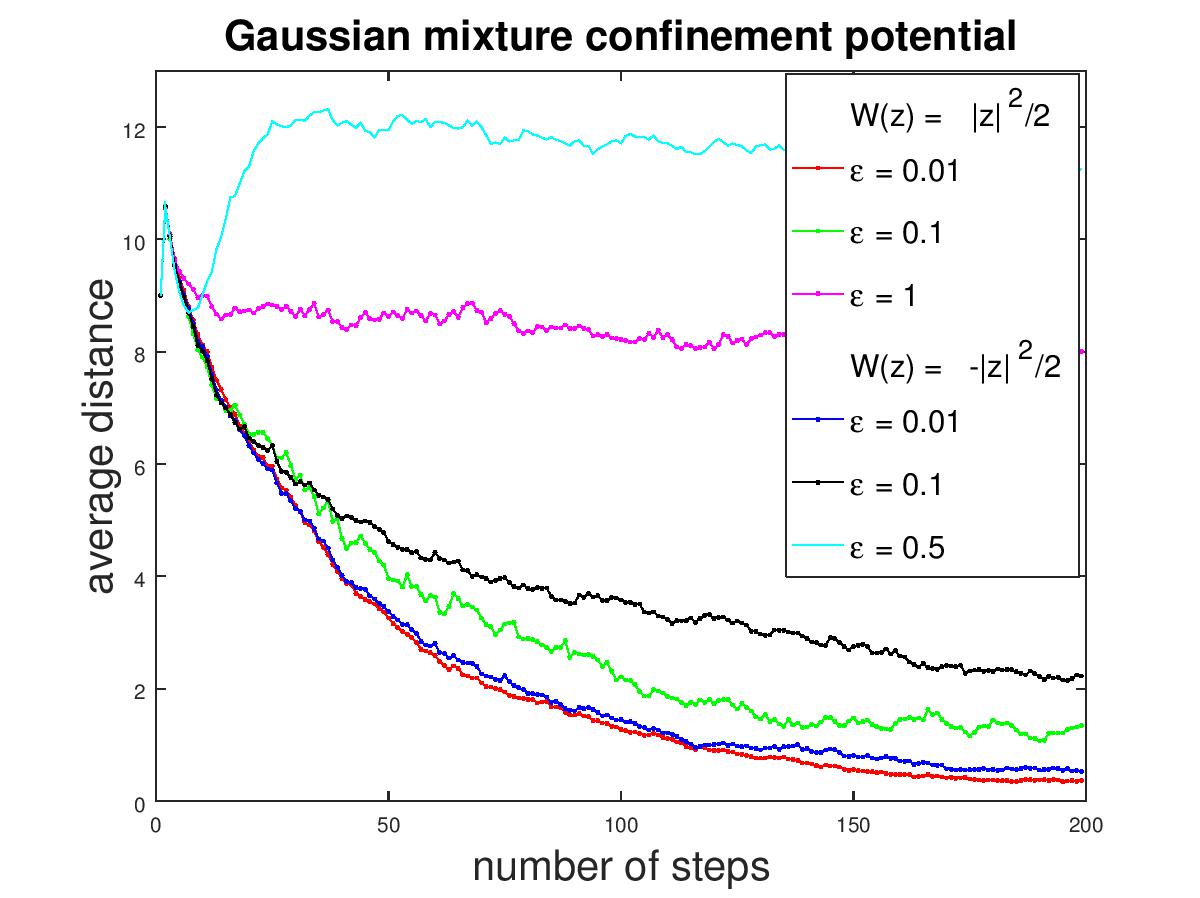} }\end{subfigure}
\begin{subfigure}{
	\includegraphics[scale=0.335]{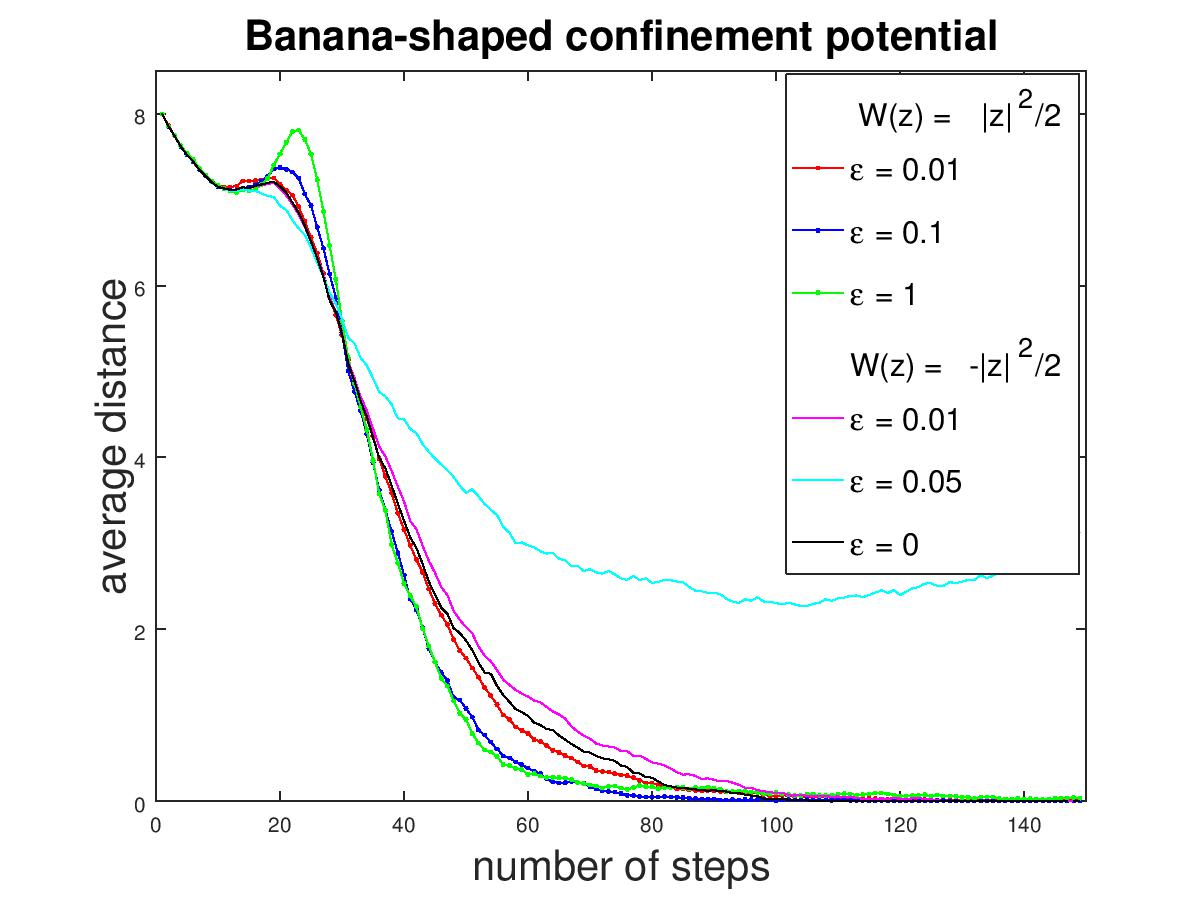}}\end{subfigure}
\caption{\small{\textit{Evolution of the mean distance $\frac{1}{n}\sum_{i=1}^n |\mathbf{X}^i_k-\mathbf{Y}^i_k|$ between the two components of the coupling for HMC after $k$ steps with $n=10$ particles for various interaction parameters $\epsilon$. This figure suggests that the particlewise coupling does not converge if the interaction is too large.}}}\label{figure3}
\end{figure}

 We simulate the coupling for mean-field potentials with non-strongly convex confinement potential to illustrate the coupling and to support our theoretical results stated in the next subsection. 

We consider two mean-field models with two different confinement potentials. The first potential is the negative logarithm of a Gaussian mixture distribution. Here, we take a mixture of 20 two-dimensional Gaussian distributions whose means are independent uniformly distributed random variables on the rectangle $[0,10]\times[0,10]$ and whose covariance matrices are the identity matrix, cf. \cite{liangwong01,kouzhouwong06,BoEbZi2020}.
The second confinement potential is the negative logarithm of a banana-shaped distribution. In particular, $V:\mathbb{R}^2\to\mathbb{R}$ is given by the Rosenbrock function $V(\mathsf{x})=(1-\mathsf{x}_1)^2+10(\mathsf{x}_2-(\mathsf{x}_1)^2)^2$, cf. \cite{BoEbZi2020}.

For the interaction between particle $i$ and $j$, we take the function $W(x^i-x^j)=(1/2)|x^i-x^j|^2$ and $\epsilon=0.01$ in \Cref{figure} and \Cref{figure2}. In \Cref{figure3}, we vary $\epsilon$ and $W$, as indicated in the legend.

The plots in \Cref{figure} show realizations of the coupling with $T=1$, $\gamma=1$ and $n=10$. The evolution of a selected particle of the coupling is drawn on a contour plot of the confinement potential. To visualize the order of the projected points they are connected by linear interpolation.
The evolution of the distance function $\frac{1}{n}\sum_{i=1}^n |\mathbf{X}_k^i-\mathbf{Y}_k^i|$ is given in the inset.
Here, $\mathbf{X}_k^i$ and $\mathbf{Y}_k^i$ are the positions of the $i$-th particles of the two realizations of the coupling after $k$ HMC steps of duration $T=1$.
The simulation terminates when the distance is smaller than $\tilde{\epsilon}=10^{-5}$.
\Cref{figure2} shows the sample average of the mean distance $\frac{1}{n}\sum_{i=1}^n |\mathbf{X}_k^i-\mathbf{Y}_k^i|$ for different numbers $n\in\{1,10,100\}$ of particles.
For $n\in\{1,10\}$ we sampled the mean distance a hundred times and for $n=100$ thirty times, since the statistical error is smaller for $n$ large.
We observe that the mean distance decreases exponentially fast after a short time, which reflects a factor $M$ appearing in the bounds in \Cref{cor_Wassersteindist} given below, and that the rate is dimension-free, i.e., independent of the number of particles. 
In \Cref{figure3}, the impact of the size of the interaction parameter $\epsilon$ is illustrated. 
We observe that for small attractive and repulsive interaction the mean coupling distance appears to converge to zero, whereas for larger interaction, particularly for large repulsive interaction (corresponding to $W(x^i-y^i)=-(1/2)|x^i-y^i|^2$) this convergence is not observed.

\section{Main results} \label{section_mainresults}

\subsection{Dimension-free contraction rate for unadjusted HMC}

To prove contraction for unadjusted HMC, we introduce a modified distance function. Define
\begin{align}
\tilde{R}&:= 8R\sqrt{(L+K)/K}, \label{eq:tildeR}
\\ \gamma &:=\min(T^{-1},\tilde{R}^{-1}/4), \label{def_gamma}
\\ R_1&:=(5/4)(\tilde{R}+2T). \label{def_R_1}
\end{align}
Note that the constants are dimension-free, i.e. independent of the number of particles.
Let $f:\mathbb{R}_+\to \mathbb{R}_+$ be given by 
\begin{align} \label{eq:definition_f}
f(r):=\int_0^r \exp(-\min(R_1,s)/T) \rmd s.
\end{align}
This function is concave and strictly increasing with $f(0)=0$ and $f'(0)=1$.
We define a metric $\rho: \mathbb{R}^{d n}\times \mathbb{R}^{d n}\to [0,\infty) $ by
\begin{align} \label{eq:distancefunction}
\rho(x,y):=\sum_{i=1}^n f(|x^i-y^i|).
\end{align}
This definition is motivated by \cite{Eberle16A} where it was introduced to obtain optimal contraction rates for weakly interacting diffusions.  This metric is equivalent to the $\ell^1$-metric, 
\begin{align} \label{eq:l^1}
 \ell^1(x,y):=\sum_i|x^i-y^i|.
\end{align}
More precisely, since $rf'(r)\leq f(r)\leq r$,
\begin{align} \label{eq:equiv_metric}
&\rho(x,y)\leq \ell^1(x,y)\leq M \rho(x,y), \qquad \qquad \text{with}
\\ & M=f'(R_1)^{-1}=\exp((5/4)(\tilde{R}/T+2)).  \label{eq:M}
\end{align}

The following theorem gives a contraction result for unadjusted HMC with respect to the metric $\rho$.

\begin{theorem}[Global contractivity for unadjusted HMC] \label{thm:generalcase}
Suppose that \Cref{ass_V_locmin}, \Cref{ass_V_lipschitz}, \Cref{ass_V_strongconv} and \Cref{ass_W_lipschitz} hold. Let $\tilde{R}$, $\gamma$, $R_1$ and $f$ be given as in \eqref{eq:tildeR}, \eqref{def_gamma}, \eqref{def_R_1} and \eqref{eq:definition_f}.
Let $T\in(0,\infty)$ and $h_1\in[0,\infty)$ satisfy
\begin{align} 
L(T+h_1)^2&\leq \frac{3}{5}\min\Big(\frac{1}{4},\frac{3K}{10 L},\frac{3}{256\cdot 5L\tilde{R}^2}\Big), \label{cond_T}
\\  h_1&\leq \frac{KT}{525L+235K}. \label{eq:cond_h_T}
\end{align} 
Let $\epsilon\in[0,\infty)$ satisfy
\begin{align} 
\epsilon\tilde{L} &< \min\Big( \frac{K}{6} , \frac{1}{2}\Big(\frac{K(\tilde{R}+T)}{36\cdot 149}\Big)^2\exp\Big(-5\frac{\tilde{R}}{T}\Big)\Big). \label{cond_epsilon}
\end{align}
Then for all $x,y\in\mathbb{R}^{d n}$ and for any $h\in[0, h_1]$ such that $h=0$ or  $T/h\in\mathbb{N}$,
\begin{align*}
\mathbb{E}\Big[\rho(\mathbf{X}(x,y),\mathbf{Y}(x,y))\Big]\leq (1-c)\rho(x,y)
\end{align*}
with contraction rate
\begin{align} \label{eq_c}
c= \frac{1}{156}K T^2\exp\Big(-\frac{5\tilde{R}}{4T}\Big).
\end{align}

\end{theorem} 
A proof is given in \Cref{section_proof_oftheorems}.
\begin{remark} \label{rem:dimensionfree}
The parameter c is dimension-free, i.e., independent of the number of particles, which is an improvement compared to the contraction rate given in \eqref{eq:contractionrate_standardcase} obtained by applying \cite[Theorem 2.3]{BoEbZi2020}. However,
it might depend implicitly on the number of degrees of freedom per particle $d$ through the parameter $\tilde{R}$.

Further, note that the contraction result holds only if the interaction parameter $\epsilon$ is sufficiently small.
For larger $\epsilon$, contraction with a dimension-free contraction rate is not guaranteed, as illustrated in \Cref{figure3}. 

\end{remark}

\begin{remark} \label{rem:adjustedHMC_contr}
For adjusted HMC one can show local contraction  by precisely bounding the effect of the accept-reject step. The case is considered  for a general potential in \cite{BoEbZi2020}. In the mean-field model for a large number $n$ of particles, an analogous local contraction result for adjusted HMC is only obtained for a restrictive choice of $h$. In particular, using the estimate for the rejection probability of \cite[Theorem 3.8]{BoEbZi2020} the discretization step  $h$ has to be chosen of order $\mathcal{O}(n^{-2})$.
\end{remark}

\begin{remark} 
\Cref{thm:generalcase} holds in particular for the product case with $\epsilon=0$. As the interaction terms vanish and some calculations simplify in that case, the condition in $T$ becomes $L(T+h_1)^2\le \min(1/4,K/L,1/(256L\tilde{R}^2))$ as in \cite{BoEbZi2020}, the condition in $h_1$ relaxes to $h_1\le 4KT/(165L)$ and the contraction rate improves to $c^{prod}=(1/39)KT^2\exp(-5\tilde{R}/(4T))$. 
If $V$ is a quadratic function, the mean-field model can be treated as a perturbation of the product model and the difference $|\mathbf{X}^{prod}(x,y)-\mathbf{Y}^{prod}(x,y)-(\mathbf{X}(x,y)-\mathbf{Y}(x,y))|$ of a coupling between to copies of the product model and two copies of the mean-field model can be bounded in terms of $\epsilon \tilde{L}\sum_{i=1}^n|x^i-y^i|$. This term can be controlled for sufficiently small $\epsilon$ by the obtained contraction for the product case. See \Cref{appendix} for the complete argument. 
\end{remark}

\subsection{Quantitative bounds for distance to the target measure} \label{subsection:dimesionfree_distance_bounds}
We deduce from \Cref{thm:generalcase} global contractivity of the transition kernel $\pi_h(x,dy)$ with respect to the Kantorovich distance based on $\rho$
\begin{align*}
\mathcal{W}_\rho (\nu,\eta) =\inf_{\omega\in\Gamma(\nu,\eta)} \int \rho(x,y)\omega(\rmd x \rmd y)
\end{align*}
on probability measures $\nu, \eta$ on $\mathbb{R}^{d n}$, where $\Gamma(\nu,\eta)$ denotes the set of all couplings of $\nu$ and $\eta$. Since the metric $\rho$ is equivalent to the $\ell^1$-distance $\ell^1$ on $(\mathbb{R}^{d})^n$ given in \eqref{eq:l^1}, contractivity with respect to $\mathcal{W}_\rho$ yields a quantitative bound on the Kantorovich distance based on $\ell^1$ on $(\mathbb{R}^{d})^n$,
\begin{align*}
\mathcal{W}_{\ell^1}(\nu{\pi_h}^m,\mu_h):=\inf_{\omega\in\Gamma(\nu{\pi_h}^m,\mu_h)} \int\sum_{i=1}^n |x^i-y^i|\omega(\rmd x\rmd y)
\end{align*}
between the law after $m$ HMC steps with initial distribution $\nu$ and invariant measure $\mu_h$.

\begin{korollar} \label{cor_Wassersteindist}
Suppose that \Cref{ass_V_locmin}, \Cref{ass_V_lipschitz}, \Cref{ass_V_strongconv} and \Cref{ass_W_lipschitz} hold. Let $T\in(0,\infty)$ and $h_1\in[0,\infty)$ satisfy \eqref{cond_T} and \eqref{eq:cond_h_T}. Let $\epsilon\in[0,\infty)$ satisfy \eqref{cond_epsilon}.
Then, for any $m\in\mathbb{N}$, for any probability measures $\nu, \eta$ on $\mathbb{R}^{d n}$, and for any $h\in[0, h_1]$ such that $h=0$ or $T/h\in\mathbb{N}$,
\begin{align}
\mathcal{W}_\rho (\nu{\pi_h}^m,\eta{\pi_h}^m)&\leq e^{-cm}\mathcal{W}_\rho(\nu,\eta), \label{eq:contraction_Wasserstein_1}
\\ \mathcal{W}_{\ell^1}(\nu{\pi_h}^m,\eta{\pi_h}^m)&\leq M e^{-cm}\mathcal{W}_{\ell^1}(\nu,\eta) \label{eq:contraction_Wasserstein_2}
\end{align}
with $c$ given by \eqref{eq_c} and $M$ given by \eqref{eq:M}. 
Further,
there exists a unique invariant probability measure ${\mu}_h$ on $\mathbb{R}^{dn}$  for the transition kernel $\pi_h$ of unadjusted HMC and 
\begin{align} \label{eq:contraction_Wasserstein_3}
\mathcal{W}_{\ell^1}(\nu{\pi_h}^m,{\mu_h})&\leq M e^{-cm}\mathcal{W}_{\ell^1}(\nu,{\mu_h}).
\end{align} 
Thus, for any constant $\tilde{\epsilon}\in(0,\infty)$ and for any initial probability distribution $\nu$ the Kantorovich distance $\Delta(m)=\mathcal{W}_{\ell^1}(\nu{\pi_h}^m,\mu_h)$ 
satisfies $\Delta(m)\leq \tilde{\epsilon}$ provided
\begin{align} \label{eq_stepnumber}
m\geq \frac{1}{c}\Big(\frac{5}{2}+\frac{5\tilde{R}}{4 T}+\log\Big(\frac{\Delta(0)}{\tilde{\epsilon}}\Big)\Big).
\end{align}
\end{korollar}
A proof is given in \Cref{section:proofs_errorbounds}. 
We note that we obtain the same bound as in \eqref{eq:contraction_Wasserstein_2} and \eqref{eq:contraction_Wasserstein_3} for the Kantorovich distance with respect to the $\ell^1$-distance averaged over all particles, $\tilde{\ell}^1(x,y)=\frac{1}{n}\sum_i|x^i-y^i|$. Then, the term $\Delta(0)/\tilde{\epsilon}$ in \eqref{eq_stepnumber} differs by a factor $1/n$. In this case, if we consider for example a product measure as initial distribution, the bound in terms of this metric does not depend logarthmically on the number of particles.


To give quantitative results of the accuracy of unadjusted HMC with respect to the target measure $\mu$, we bound the strong accuracy of velocity Verlet. The exact dynamics started in $(x,\xi)$ with $h=0$ is denoted by $(q_s(x,\xi),p_s(x,\xi))$ and the position of the dynamics started in $(x,\xi)$ with $h>0$ is denoted by $(\tilde{q}_s(x,\xi),\tilde{p}_s(x,\xi))$.

\begin{theorem}[Strong accuracy of velocity Verlet] \label{thm:flowdifference_exact_num}
Suppose that \Cref{ass_V_locmin}, \Cref{ass_V_lipschitz} and \Cref{ass_W_lipschitz} hold. Let $T\in(0,\infty)$ satisfy $(L+4\epsilon\tilde{L})T^2\leq (1/4)$. For $x\in\mathbb{R}^{d n}$, for any $h\in(0,\infty)$ with $T/h\in\mathbb{N}$ and $k\in\mathbb{N}$ with $kh\leq T$, it holds
\begin{align} \label{eq:first_estimate}
\mathbb{E}_{\xi\sim \mathcal{N}(0,I_{d n})}\Big[\sum_i|q^i_{kh}(x,\xi)-\tilde{q}^i_{kh}(x,\xi)|\Big]\leq h C_2\Big(d^{1/2}n+\sum_i |x^i|\Big)
\end{align}
with $C_2$ depending on $L$, $\tilde{L}$, $\epsilon$ and $T$.
If additionally \Cref{ass_V_thirdfourthder} and \Cref{ass_W_thirdfourthder} are supposed, then 
for $x\in\mathbb{R}^{d n}$, for any $h>0$ with $T/h\in\mathbb{N}$ and $k\in\mathbb{N}$ with $kh\leq T$,
\begin{align} \label{eq:second_estimate}
\mathbb{E}_{\xi\sim \mathcal{N}(0,I_{d n})}\Big[\sum_i|q^i_{kh}(x,\xi)-\tilde{q}^i_{kh}(x,\xi)|\Big]\leq h^2 \tilde{C}_2\Big(dn+\sum_i |x^i|+\sum_i |x^i|^2\Big)
\end{align}
with $\tilde{C}_2$ depending on $L$, $\tilde{L}$, $\epsilon$, $L_H$, $\tilde{L}_H$ and $T$.
\end{theorem}
A proof is given in \Cref{section:proofs_errorbounds}. 

We obtain a bound on the difference between the invariant measure $\mu_h$ and the target measure $\mu$, by using the contraction result of \Cref{thm:generalcase} and by applying a triangle inequality trick, which is mentioned in \cite[Remark 6.3]{mattingly2010convergence} and has been used in many other works.
In particular, it holds
\begin{align*}
\mathcal{W}_\rho(\mu,\mu_h)&=\mathcal{W}_\rho(\mu\pi,\mu_h{\pi_h})\leq \mathcal{W}_\rho(\mu\pi,\mu\pi_h)+\mathcal{W}_\rho(\mu\pi_h,\mu_h\pi_h)
\\ & \leq \mathcal{W}_\rho(\mu\pi,\mu\pi_h)+(1-c)\mathcal{W}_\rho(\mu,\mu_h).
\end{align*}
Hence, by \eqref{eq:equiv_metric} 
\begin{align*}
\mathcal{W}_{\ell^1}(\mu,\mu_h)\leq Mc^{-1}\mathcal{W}_{\ell^1}(\mu\pi,\mu\pi_h)\leq Mc^{-1}\mathbb{E}_{x\sim\mu, \ \xi\sim \mathcal{N}(0,I_{d n})}\Big[\sum_i|q^i_{kh}(x,\xi)-\tilde{q}^i_{kh}(x,\xi)|\Big]
\end{align*}
with $M$ given in \eqref{eq:M}.
Inserting \eqref{eq:first_estimate}, respectively \eqref{eq:second_estimate}, yields the following result.

\begin{korollar}[Asymptotic Bias]\label{cor:bound_mu} Suppose that \Cref{ass_V_locmin}, \Cref{ass_V_lipschitz}, \Cref{ass_V_strongconv} and \Cref{ass_W_lipschitz} hold. Let $T$ and $h_1$ satisfy \eqref{cond_T}. Let $\epsilon$ satisfy \eqref{cond_epsilon}. Let  $C_2$ and $\tilde{C}_2$ be as in \Cref{thm:flowdifference_exact_num}. Then for $h\in(0, h_1]$ with $T/h\in\mathbb{N}$, 
\begin{align*}
\mathcal{W}_{\ell^1} (\mu,{\mu_h})\leq hc^{-1}MC_2\Big(d^{1/2}n+\int_{\mathbb{R}^{nd}}\sum_i |x^i|\mu(\rmd x)\Big) 
\end{align*}
with $c$ given by \eqref{eq_c} and $M$ given by \eqref{eq:M}.
If additionally \Cref{ass_V_thirdfourthder} and \Cref{ass_W_thirdfourthder} are assumed, then for $h\in(0, h_1]$ with $T/h\in\mathbb{N}$,
\begin{align*}
\mathcal{W}_{\ell^1} (\mu,{\mu_h})\leq h^2c^{-1}M \tilde{C}_2\Big(dn+\int_{\mathbb{R}^{nd}}\sum_i |x^i|\mu(\rmd x)+\int_{\mathbb{R}^{nd}}\sum_i |x^i|^2\mu(\rmd x)\Big). 
\end{align*}
\end{korollar}

Note that the bound in \Cref{cor:bound_mu} is linear in the number $n$ of particles. 

For unadjusted HMC, \Cref{cor_Wassersteindist} gives exponential convergence to the invariant measure ${\mu_h}$. In the next theorem, we give a bound on the number of steps to reach the target measure $\mu$ up to a given error.

\begin{theorem}[Complexity Guarantee] \label{thm:error_invmeas_stepnumber}
Suppose that \Cref{ass_V_locmin}, \Cref{ass_V_lipschitz}, \Cref{ass_V_strongconv} and \Cref{ass_W_lipschitz} hold. Let $T\in(0,\infty)$ and $h_1\in (0,\infty)$ satisfy \eqref{cond_T} and \eqref{eq:cond_h_T}. Let $\epsilon\in[0,\infty)$ satisfy \eqref{cond_epsilon}. Let $\nu$ be a probability measure on $\mathbb{R}^{d n}$, and let $\Delta (m) =\mathcal{W}_{\ell^1}(\nu{\pi_h}^m,\mu)$ denote the Kantorovich distance with respect to $\ell^1$ to the target probability measure $\mu$ after $m$ steps with initial distribution $\nu$. For some $\tilde{\epsilon}\in(0,\infty)$, let $m\in\mathbb{N}$ be such that 
\begin{align} \label{eq:estimateofstepnumbers_num}
m\geq \frac{1}{c}\Big(\frac{5}{2}+\frac{5R}{4T}+\log\Big(\frac{2\mathcal{W}_{\ell^1}(\mu_h,\nu)}{\tilde{\epsilon}}\Big)^+ \Big)
\end{align}
with $c$ given by \eqref{eq_c}. Then, there exists $h_2$ such that
for $h\in(0,\min(h_1,h_2)]$ with $T/h\in\mathbb{Z}$, 
\begin{align} \label{eq:deltam}
\Delta(m)\leq \tilde{\epsilon}
\end{align}
where for fixed $K$, $L$, $\tilde{L}$, $\epsilon$, $R$ and $T$, $h_2^{-1}$ is of order $\mathcal{O}(\tilde{\epsilon}^{-1}(d^{1/2}n+\int\sum_i|x^i|\mu(\rmd x)))$.
If additionally \Cref{ass_V_thirdfourthder} and \Cref{ass_W_thirdfourthder} are assumed, then there exists $\tilde{h}_2$ such that for $h\in(0,\min(h_1,\tilde{h}_2)]$ with $T/h\in\mathbb{Z}$, \eqref{eq:deltam} holds, where for fixed $K$, $L$, $\tilde{L}$, $L_H$, $\tilde{L}_H$, $\epsilon$, $R$ and $T$, $\tilde{h}_2^{-1}$ is of order $\mathcal{O}(\tilde{\epsilon}^{-1/2}((nd)^{1/2}+\sqrt{\int\sum_i|x^i|\mu(\rmd x)}+\sqrt{\int\sum_i|x^i|^2\mu(\rmd x)}))$.
\end{theorem}
A proof is given in \Cref{section:proofs_errorbounds}.
If we consider the averaged distance $\tilde{\ell}^1$ instead of $\ell^1$, the argument in the logarithmic term in \eqref{eq:estimateofstepnumbers_num} changes by a factor $1/n$ and the logarithmic dependence on $n$ in $h_2$ and $\tilde{h}_2$ vanishes.

\begin{remark} 
We note that $h^{-1}$ is $\mathcal{O}(n^{1/2})$ in \Cref{thm:error_invmeas_stepnumber} and hence it grows sublinear in $n$. 
%
Further, the constant $\tilde{C}_2$ obtained in the proof of \Cref{thm:flowdifference_exact_num} is $\mathcal{O}(T^{-1})$. For the numerical method uLA, which forms a special case of unadjusted HMC with $h=T$ (see  \cite[Section 5.2]{neal11}), we obtain that $h^{-1}=T^{-1}$ has to be chosen of order $\mathcal{O}(n)$, which corresponds to the results in \cite[Example 18]{DuEb21}.
Therefore, an $\tilde{\epsilon}$-accurate approximation of the target measure in the $\mathcal{W}_{\ell_1}$ distance can be achieved by  uHMC applied to the $n$-particle mean-field  system  with $O(n^{1/2} \tilde{\epsilon}^{-1/2} \log(n/\tilde{\epsilon}))$ gradient evaluations; whereas the corresponding complexity of uLA is $O(n \tilde{\epsilon}^{-1} \log(n/\tilde{\epsilon}))$.

\end{remark}

\begin{remark} From \Cref{thm:error_invmeas_stepnumber}, note that the number of evaluations of  the gradient $\nabla U(x)$ in each step of duration $T$ is $\mathcal{O}(n^{1/2})$ for fixed $K$, $L$, $\tilde{L}$, $\epsilon$, $T$, $R$, $d$ and $h$.
If we assume that the computation of the gradient in one step is $\mathcal{O}(n)$, then the overall complexity of unadjusted HMC is $\mathcal{O}(n^{3/2})$.
\end{remark}

\subsection{Dimension-free bounds for ergodic averages of intensive observables}\label{sec:dimension_free_bounds}

Next, we define the ergodic averages $A_{m,b}g$, which approximate $\mu(g)=\int g(x)\mu(\rmd x)$, by
\begin{align} \label{eq:estimator}
A_{m,b}g:=\frac{1}{m}\sum_{i=b}^{b+m-1} g(\mathbf{X}_i),
\end{align}
for some function $g:\mathbb{R}^{dn}\to\mathbb{R}$ 
and for $b,m\in\mathbb{N}$, where $(\mathbf{X}_n)$ is the Markov chain given by unadjusted HMC. The parameter $b$ corresponds to the burn-in time. Here, we consider bounded and continuously differentiable observables, i.e.,  $g\in\mathcal{C}_b^1(\mathbb{R}^{dn})$. 
Quantitative bounds on the bias of the ergodic averages follow by the exponential convergence in the Kantorovich distance with respect to the $\ell^1$ metric given in \eqref{eq:l^1} and the bound on the accuracy of unadjusted HMC.  
\begin{theorem}[Bias of Ergodic Averages] \label{thm:ergodicaverages}
Let $g\in\mathcal{C}_b^1(\mathbb{R}^{dn})$ with $\max_i\|\nabla_i g\|_{\infty}<\infty$. Suppose that \Cref{ass_V_locmin}, \Cref{ass_V_lipschitz}, \Cref{ass_V_strongconv} and \Cref{ass_W_lipschitz} hold. Let $T\in(0,\infty)$ and $h_1\in[0,\infty)$ satisfy \eqref{cond_T} and \eqref{eq:cond_h_T}. Let $\epsilon\in[0,\infty)$ satisfy \eqref{cond_epsilon}. Let $\nu$ be a probability measure on $\mathbb{R}^{dn}$. Let $C_2$ and $\tilde{C}_2$ be given as in \Cref{thm:flowdifference_exact_num}, and let $c$ be given as in \eqref{eq:contractionrate_standardcase}.
Then for $h\in[0,h_1]$ such that $h=0$ or $T/h\in\mathbb{N}$,
\begin{align*}
|\mathbb{E}_\nu[A_{m,b}g]-\mu(g)|\leq \frac{1}{m} \max_i\|\nabla_i g\|_{\infty}\frac{e^{-cb}}{1-e^{-c}}\mathcal{W}_{\ell^1}(\nu,\mu_h)+h \max_i\|\nabla_i g\|_{\infty}C_3,
\end{align*}
where $C_3=\exp(\frac{5}{4}(2+\tilde{R}/T)) c^{-1}C_2\Big(d^{1/2}n+\int \sum_i |x^i|\mu(\rmd x)\Big)$.
If additionally \Cref{ass_V_thirdfourthder} and \Cref{ass_W_thirdfourthder} are supposed, then
\begin{align*}
|\mathbb{E}_\nu[A_{m,b}g]-\mu(g)|&\leq \frac{1}{m}\max_i\|\nabla_i g\|_{\infty}\frac{e^{-cb}}{1-e^{-c}}\mathcal{W}_{\ell^1}(\nu,\mu_h)  +h^2 \max_i\|\nabla_i g\|_{\infty}\tilde{C}_3,
\end{align*}
where $\tilde{C}_3=\exp(\frac{5}{4}(2+\tilde{R}/T))c^{-1}\tilde{C}_2\Big(dn+\int \sum_i |x^i|\mu(\rmd x)+\int \sum_i |x^i|^2\mu(\rmd x)\Big)$.
\end{theorem}
A proof is given in \Cref{section:proofs_dim_free_bounds}. 

\begin{remark} \label{rem:intensive_observables}
We note that provided $\max_i \|\nabla_i g\|_\infty$ is $\mathcal{O}(1/n)$ the bound of the bias of the ergodic averages is independent of the number $n$ of particles. Hence for intensive observables of the form $g(x)=\frac{1}{n}\sum_i\hat{g}(x^i)$ where $\hat{g}\in\mathcal{C}_b^1(\mathbb{R}^d)$ with $\|\nabla \hat{g}\|_\infty<\infty$, \Cref{thm:ergodicaverages} gives quantitative bounds on the bias of their ergodic averages which are dimension-free, i.e., independent of the number $n$ of particles. 
Whereas, for extensive observables, where $\max_i \|\nabla_i g\|_\infty$ is $\mathcal{O}(1)$, the bound depends on the number $n$ of particles. 
\end{remark}

\section{Estimates for the dynamics \texorpdfstring{\eqref{eq:hamdyn_num}}{}} \label{section_estimates}

\subsection{Deviation from free dynamics}
Here we apply the Lipschitz conditions in \Cref{ass_V_lipschitz} and \Cref{ass_W_lipschitz} to obtain bounds on how far the dynamics in \eqref{eq:hamdyn_num} deviates from the free dynamics, $U \equiv 0$.
To obtain these bounds, we assume in the following that $t, h\in[0,\infty)$ are such that $t/h\in\mathbb{N}$ for $h>0$ and such that 
\begin{align} \label{eq:cond_t}
(L+4\epsilon\tilde{L})(t^2+th)\leq 1.
\end{align} 
This condition essentially states that the duration of the Hamiltonian dynamics in \eqref{eq:hamdyn_num} is small with respect to the fastest characteristic time-scale of the mean-field particle system represented by $\sqrt{ \sup \| \mathrm{Hess} U \| } \le \sqrt{L+4 \epsilon \tilde L}$ (see \Cref{rem:effectiveLip}). This bound follows from \Cref{ass_V_lipschitz} and \Cref{ass_W_lipschitz}. 
The $i$-th component of the solution to \eqref{eq:hamdyn_num} is denoted by $(x_s^i,v_s^i)$.
 
\begin{lemma} \label{lemma:positionestimate}
Let $x,v\in\mathbb{R}^{d n}$. 
Then for $i\in\{1,...,n\}$,
\begin{align}
\max_{s\leq t}|x_s^i|&\leq (1+(L+2\epsilon\tilde{L})(t^2+th))\max(|x^i|,|x^i+t v^i|) \label{eq:position_estimate2}
\\ &  +\frac{2\epsilon\tilde{L}(t^2+th)}{n}\max_{s\leq t}\sum_{j\neq i} |x^j_s|,  \nonumber
\\ \max_{s\leq t} |v^i_s|  &\leq  |v^i|+(L+2\epsilon\tilde{L})t(1+(L+2\epsilon\tilde{L})(t^2+th))\max(|x^i|,|x^i+t v^i|)) \label{eq:velocity_estimate2}
\\ & +\frac{2\epsilon\tilde{L}t}{n}(1+(L+2\epsilon\tilde{L})(t^2+th))\max_{s\leq t} \sum_{j\neq i} |x^j_s|. \nonumber 
\end{align}
Moreover,
\begin{align}
 \max_{s\leq t}\sum_i|x_s^i|&\leq (1+(L+4\epsilon\tilde{L})(t^2+th))\sum_i\max(|x^i|,|x^i+t v^i|), \label{eq:position_estimate3}
\\ \max_{s\leq t} \sum_i|v^i_s|&\leq  (L+4\epsilon\tilde{L})t(1+(L+4\epsilon\tilde{L})(t^2+th))\sum_i\max(|x^i|,|x^i+t v^i|)
\label{eq:velocity_estimate3}
\\ &  +\sum_i|v^i|.  \nonumber
\end{align}

\end{lemma}
A proof of \Cref{lemma:positionestimate} is provided in \Cref{section_proof_oflemmata}.

Let two processes $(x_s,v_s)$, $(y_s,u_s)$ with initial values $(x,v)$ and $(y,u)$ be driven by the Hamiltonian dynamics in \eqref{eq:hamdyn_num}. 
We set $(z_s,w_s):=(x_s-y_s,v_s-u_s)$.
Since $(x_s,v_s)$ and $(y_s,u_s)$ depend on $(x,v)$ and $(y,u)$, respectively, $(z_s,w_s)$ depends on $(x,v,y,u)$. By \eqref{eq:hamdyn_num}, the dynamics of the $i$-th component of $(z_s,w_s)$ is given by 
\begin{equation}
\begin{aligned}
\frac{\rmd }{\rmd t}z_t^i&=w_{\lfloor t\rfloor }^i-(h/2)(\nabla_i U(x_{\lfloor t \rfloor})-\nabla_i U(y_{\lfloor t \rfloor})) 
\\ \frac{\rmd }{\rmd t}w_t^i&=(1/2)(-\nabla_i U(x_{\lfloor t\rfloor})-\nabla_i U(x_{\lceil t\rceil})+\nabla_i U(y_{\lfloor t \rfloor})+\nabla_i U(y_{\lceil t \rceil})). \label{eq:hamdyn_ofdifferenceprocess}
\end{aligned} 
\end{equation}
Next, we bound the distance between the process $(z_s^i,w_s^i)$ and the process given by the free dynamics, where $U \equiv 0$.
As the particlewise coupling in \Cref{subsec:coupling} is designed with respect to the free dynamics, 
this bound plays an important role in proving the contraction results of \Cref{section_mainresults}.
It explains why the particlewise coupling works when the distance between $i$-th particles is small, i.e., when $|x^i-y^i|<\tilde{R}$, and when the duration $T$ and the time step $h$ are small, i.e., when \eqref{eq:cond_t} is assumed.

\begin{lemma} \label{lemma:differencepositionestimate}
Let $x,y,v,u\in\mathbb{R}^{d n}$. 
Then for all $i\in\{1,...,n\}$,
\begin{align}
&\max_{s\leq t } |z_s^i-z^i-sw^i|
\leq (L+2\epsilon \tilde{L})(t^2+th)\max (|z^i+t w^i|,|z^i|)\label{eq:position_estimate4}
\\ &\hspace{2cm}+\frac{2\epsilon \tilde{L}(t^2+th)}{n} \max_{s\leq t } \sum_{j\neq i} |z_s^j|, \nonumber
\\ 
&\max_{s\leq t } |z_s^i|
\leq (1+(L+2\epsilon \tilde{L})(t^2+th))\max (|z^i+t w^i|,|z^i|)\label{eq:position_estimate5}
\\ &\hspace{2cm}+\frac{2\epsilon \tilde{L}(t^2+th)}{n} \max_{s\leq t } \sum_{j\neq i} |z_s^j|, \nonumber
\\
&\max_{s\leq t}  |w_s^i-w^i|
\leq (L+2\epsilon \tilde{L})t(1+(L+2\epsilon\tilde{L})(t^2+th)) \max(|z^i+t w^i|,|z^i|)) \label{eq:velocity_estimate4}
\\ &\hspace{2cm}+ \frac{2\epsilon\tilde{L}t}{n}(1+(L+2\epsilon\tilde{L})(t^2+th))\max_{s\leq t } \sum_{j\neq i} |z_s^j|, \nonumber
\\
& \max_{s\leq t }  |w_s^i|
\leq  |w^i|+(L+2\epsilon \tilde{L})t(1+(L+2\epsilon\tilde{L})(t^2+th) )\max(|z^i+t w^i|,|z^i|)) \label{eq:velocity_estimate5}
\\ &\hspace{2cm}+ \frac{2\epsilon\tilde{L}t}{n}(1+(L+2\epsilon\tilde{L})(t^2+th))\max_{s\leq t } \sum_{j\neq i} |z_s^j|. \nonumber
\end{align}
Moreover,
\begin{align}
\max_{s\leq t }\sum_i|z_s^i|
&\leq  (1+(L+4\epsilon\tilde{L})(t^2+th))\sum_i\max(|z^i+t w^i|,|z^i|), \label{eq:position_estimate6}
\\
\max_{s\leq t }\sum_i|w_s^i|
&\leq  (L+4\epsilon\tilde{L})t(1+(L+4\epsilon\tilde{L})(t^2+th))\sum_i\max(|z^i+t w^i|,|z^i|)
\label{eq:velocity_estimate6}
\\ &  +\sum_i|w^i|. \nonumber
\end{align}

\end{lemma}
A proof of \Cref{lemma:differencepositionestimate} is provided in \Cref{section_proof_oflemmata}.

\subsection{Bounds in region of strong convexity}
Next, we obtain a bound for the difference
between the positions of the $i$-th particles provided that $|x^i-y^i|>\tilde{R}$ and $v^i=u^i$.
We assume that
\begin{align} \label{eq:cond_t_num}
(L+4\epsilon\tilde{L})(t^2+th)\leq\min\Big(\frac{\kappa}{L+4\epsilon\tilde{L}},\frac{1}{4}\Big),
\end{align}
where $\kappa$ is given by
\begin{align}\label{eq:def_kappa}
\kappa:=K-3\epsilon\tilde{L}.
\end{align}
Further, we assume that
\begin{align} \label{eq:cond_h_t}
h\leq  \frac{Kt}{525L+235K}.
\end{align}

\begin{lemma}\label{lemma:contraction_num} Suppose that \Cref{ass_V_locmin}, \Cref{ass_V_lipschitz}, \Cref{ass_V_strongconv} and \Cref{ass_W_lipschitz} hold. Let $\epsilon\in[0,\infty)$ be such that $\epsilon\tilde{L}< K/6$ holds. Let $\tilde{R}$ be given in \eqref{eq:tildeR}. 
Let $t,h \in[0,\infty)$ be such that $h=0$ or $t/h\in\mathbb{N}$, and such that \eqref{eq:cond_t_num} and \eqref{eq:cond_h_t} holds.
Then, for all $x,y,v,u\in\mathbb{R}^{dn}$ and $i\in\{1,...,n\}$ such that $|x^i-y^i|\geq\tilde{R}$ and $v^i=u^i$,
\begin{align}
|x_t^i-y_t^i|^2&\leq \Big(1-\frac{1}{4}\kappa t^2\Big)|x^i-y^i|^2 + 2
\frac{\epsilon\tilde{L}t^2}{n^2}\Big(\max_{s\leq t}\sum_{j\neq i}|x_s^j-y_s^j|\Big)^2. \label{eq:contraction_num_bound}
\end{align}
\end{lemma}
A proof of \Cref{lemma:contraction_num} is given in \Cref{section_proof_oflemmata}.

In the strongly convex case with only one particle (i.e., $R=0$, $\epsilon=0$ and $n=1$), an improved version of \Cref{lemma:contraction_num} with less restrictive assumptions on $T$ and $h$ is given in \Cref{appendixA} in \Cref{lem:contr_verlet_flow}. This bound provides directly contraction in $L^p$ Wasserstein distance provided $T>0$ and $h\ge 0$ satisfy $LT^2\le 20^{-1}$ and $T/h\in \mathbb{Z}$ if $h>0$, see \Cref{appendixA}.

\section{Proof of results from Section \ref{section_estimates}} \label{section_proof_oflemmata}

Before stating the proofs of \Cref{section_estimates}, note that by 
\eqref{eq:conseq_ass1} and \eqref{eq:conseq_ass2} 
for all $x, y\in\mathbb{R}^{dn}$,
\begin{align} 
&|\nabla_i U(x)| \leq L|x^i|+\frac{2\epsilon\tilde{L}}{n}\sum_{j\neq i} |x^i - x^j| \leq (L+2\epsilon\tilde{L})|x^i|+\frac{2\epsilon\tilde{L}}{n}\sum_{j\neq i} |x^j| ,
 \label{eq:bound_deltaiU}
 \end{align}
 and by \Cref{ass_V_lipschitz} and \Cref{ass_W_lipschitz}
 \begin{align}
&|\nabla_i U(x)-\nabla_i U(y)| \leq (L+2\epsilon\tilde{L})|x^i-y^i|+\frac{2\epsilon\tilde{L}}{n}\sum_{j\neq i} |x^j-y^j|. \label{eq:bound_deltaiU2}
\end{align}
Further by \eqref{eq:cond_strongconv} and \eqref{eq:conseq_ass2}, it holds for all $x,y\in\mathbb{R}^{d n}$,
\begin{align}
-(x^i-y^i)&\cdot(\nabla_i U(x)-\nabla_i U(y)) \nonumber
\\  & \leq -(K-2\epsilon\tilde{L})|x^i-y^i|^2+2\epsilon\tilde{L}|x^i-y^i|\frac{1}{n}\sum_{j\neq i}|x^j-y^j| +\hat{C} \nonumber
\\ & \leq  -\kappa|x^i-y^i|^2+\epsilon\tilde{L}\Big(\frac{1}{n}\sum_{j\neq i}|x^j-y^j|\Big)^2 +\hat{C}.  \label{eq:conseq_strongconv}
\end{align}
It follows from the definition \eqref{eq:tildeR} of $\tilde{R}$ and the condition $\epsilon\tilde{L}< K/6$, which is assumed in \Cref{lemma:contraction_num}, that
for all $\mathsf{x},\mathsf{y}\in\mathbb{R}^d$ with $|\mathsf{x}-\mathsf{y}|\geq\tilde{R}$,
\begin{align} \label{eq:strong_conv_hatC}
\hat{C}=R^2(L+K)<\frac{1}{64}K|\mathsf{x}-\mathsf{y}|^2\leq \frac{1}{32}\kappa|\mathsf{x}-\mathsf{y}|^2.
\end{align} 

\begin{proof}[Proof of \Cref{lemma:positionestimate}]
Fix $x,v\in\mathbb{R}^{d n}$. Let $s\leq t$. 
We have from \eqref{eq:hamdyn_num} 
\begin{align*}
x^i_s-x^i-sv^i&= \int_0^s\int_0^{\lfloor r\rfloor}\Big( -\frac{1}{2}\nabla_i U(x_{\lfloor u \rfloor})-\frac{1}{2}\nabla_i U(x_{\lceil u \rceil})\Big) \rmd u \ \rmd r- \int_0^s \frac{h}{2}\nabla_i U(x_{\lfloor r\rfloor})\rmd r.
\end{align*}
We apply \eqref{eq:bound_deltaiU} to obtain
\begin{align*}
|x^i_s-x^i-sv^i| &\leq \frac{(L+2\epsilon\tilde{L})(t^2+th)}{2}\max_{r\leq t } (|x_r^i-x^i-rv^i|+|x^i+rv^i|)
\\ & +\frac{2\epsilon\tilde{L}(t^2+th)}{2n}\max_{r\leq t }\sum_{j\neq i}|x_r^j|. 
\end{align*}
Invoking condition \eqref{eq:cond_t}, we get
\begin{align*}
\max_{s\leq t }|x^i_s-x^i&-sv^i|\leq (L+2\epsilon\tilde{L})(t^2+th) \max_{s\leq t }|x^i+sv^i|+\frac{2\epsilon\tilde{L}(t^2+th) }{n}\max_{s\leq t }\sum_{j\neq i}|x_s^j|
\\ & = (L+2\epsilon\tilde{L})(t^2+th) \max(|x^i|,|x^i+t v^i|)+\frac{2\epsilon\tilde{L}(t^2+th) }{n}\max_{s\leq t }\sum_{j\neq i}|x_s^j|.
\end{align*}
By applying the triangle inequality, \eqref{eq:position_estimate2} is obtained. 
From \eqref{eq:hamdyn_num} and \eqref{eq:bound_deltaiU}, we have
\begin{align}
|v_s^i-v^i|&\leq \int_0^s \max_{u\leq t }|\nabla_i U(x_u)|\rmd r \leq (L+2\epsilon\tilde{L})t\max_{u\leq t } |x_u^i|+\frac{2\epsilon\tilde{L}t}{n}\max_{u\leq t }\sum_{j\neq i}|x_u^j|.
\label{eq:velocity_intermediate_estimate}
\end{align} 
We insert \eqref{eq:position_estimate2} in \eqref{eq:velocity_intermediate_estimate} to obtain
\begin{align*}
|v_s^i-v^i|&\leq  (L+2\epsilon\tilde{L})t(1+(L+2\epsilon\tilde{L})(t^2+th))\max(|x^i|,|x^i+t v^i|)
\\ & +\frac{2\epsilon\tilde{L}t}{n}(1+(L+2\epsilon\tilde{L})(t^2+th))\max_{u\leq t }\sum_{j\neq i}|x_u^j|.
\end{align*}
By applying the triangle inequality, \eqref{eq:velocity_estimate2} is obtained. 
Equation \eqref{eq:position_estimate3} and \eqref{eq:velocity_estimate3} follow by considering the sum over all particles, i.e., by \eqref{eq:hamdyn_num} we have
\begin{align*}
&\sum_i|x_s^i-x^i-sv^i|
\\ &\leq \int_0^s\int_0^r \frac{1}{2} \sum_i|\nabla_i U(x_{\lfloor u \rfloor})+\nabla_i U(x_{\lceil u\rceil})| \rmd u \rmd r +\frac{h}{2}\int_0^s \sum_i|\nabla_i U(x_{\lfloor r \rfloor})|\rmd r
\\ & \leq \frac{(L+4\epsilon\tilde{L})(t^2+th)}{2}\max_{r\leq t}\Big( \sum_i|x_r^i|\Big)
\end{align*}
and hence analogous to the estimate obtained for the $i$-th particle,
\begin{align*}
\max_{s\leq t}\sum_i |x_s^i-x^i-sv^i|&\leq (L+4\epsilon\tilde{L})(t^2+th)\max_{r\leq t}\sum_i|x^i+rv^i|
\\ & \leq (L+4\epsilon\tilde{L})(t^2+th)\sum_i\max(|x^i|,|x^i+tv^i|).
\end{align*}
By applying the triangle inequality, \eqref{eq:position_estimate3} is obtained.
By \eqref{eq:hamdyn_num} and \eqref{eq:position_estimate3},
\begin{align*}
\sum_i|v_s^i-v^i|&\leq (L+4\epsilon\tilde{L})t\max_{r\leq t } \Big(\sum_i|x_r^i|\Big)
\\ & \leq (L+4\epsilon\tilde{L})t(1+(L+4\epsilon\tilde{L})(t^2+th))\sum_i\max(|x^i|,|x^i+tv^i|),
\end{align*}
and \eqref{eq:velocity_estimate3} is obtained by the triangle inequality.
\end{proof}

\begin{proof} [Proof of \Cref{lemma:differencepositionestimate}]
By \eqref{eq:bound_deltaiU2} and \eqref{eq:hamdyn_ofdifferenceprocess},
\begin{align*}
&|z_s^i-z^i-sw^i|
\\ & \leq \int_0^s\int_0^r \max_{v\leq t }|-\nabla_i U(x_v)+\nabla_i U(y_v)|\rmd u\ \rmd r
+\frac{h}{2}\int_0^s \max_{v\leq t}|-\nabla_i U(x_v)+\nabla_i U(y_v)| \rmd r
\\ &\leq \frac{(L+2\epsilon\tilde{L})(t^2+th)}{2}\max_{r\leq t }|z_r^i|+\frac{2\epsilon\tilde{L}(t^2+th)}{2n}\max_{r\leq t } \sum_{j\neq i}|z_r^j|.
\end{align*}
Hence, we obtain similar to the previous proof
\begin{align*}
\max_{s\leq t }|z_s^i-z^i-sw^i|
&\leq (L+2\epsilon\tilde{L})(t^2+th)\max(|z^i|,|z^i+t w^i|)
\\&  +\frac{2\epsilon\tilde{L}(t^2+th)}{n}\max_{s\leq t } \sum_{j\neq i}|z_s^j|,
\end{align*}
which gives \eqref{eq:position_estimate4}. Then \eqref{eq:position_estimate5} is obtained by applying triangle inequality.
Next, we consider
\begin{align*}
|w^i_s-w^i|&\leq \frac{1}{2}\int_0^s(|-\nabla_i U(x_{\lfloor r\rfloor})+\nabla_i U(y_{\lfloor r\rfloor })|+|-\nabla_i U(x_{\lceil r\rceil})+\nabla_i U(y_{\lceil r\rceil })|) \rmd r
\\ & \leq (L+2\epsilon\tilde{L})t\max_ {r\leq t } |z_r^i|+\frac{2\epsilon\tilde{L}t}{n}\max_{r\leq t } \sum_{j\neq i} |z_r^j|,
\end{align*}
where we again used \eqref{eq:bound_deltaiU2} and \eqref{eq:hamdyn_ofdifferenceprocess}.
Hence, we obtain by \eqref{eq:position_estimate5},
\begin{align*}
\max_{s\leq t} |w^i_s-w^i|
&\leq (L+2\epsilon\tilde{L})t(1+(L+2\epsilon\tilde{L})(t^2+th))\max(|z^i|,|z^i+t w^i|)
\\ & + \frac{2\epsilon\tilde{L}t}{n}(1+(L+2\epsilon\tilde{L})(t^2+th)) \max_{s\leq t } \sum_{j\neq i}|z_s^j|,
\end{align*} 
which gives \eqref{eq:velocity_estimate4} and \eqref{eq:velocity_estimate5}. Estimates \eqref{eq:position_estimate6} and \eqref{eq:velocity_estimate6} hold similarly by considering the sum over all particles instead of considering only the $i$-th particle.
\end{proof}

\begin{proof}[Proof of \Cref{lemma:contraction_num}]
As before, write $(z_s,w_s)=(x_s-y_s,v_s-u_s)$ whose dynamics is given by \eqref{eq:hamdyn_ofdifferenceprocess}. Then, $z_0=x-y$ and $w^i_0=0$ since the velocities of the $i$-th component are synchronized.

Define $a^i(t)=|z_t^i|^2$ and $b^i(t)=2z_t^i\cdot w_t^i$. We set up an initial value problem of the two deterministic processes $a^i(t)$ and $b^i(t)$ and solve it to obtain the required bound for $a^i(t)$.
By \eqref{eq:hamdyn_ofdifferenceprocess}, we have

\begin{align*}
\frac{\rmd}{\rmd t}{a}^i(t)&=b^i(t)+2z_t^i\cdot(w_{\lfloor t\rfloor}^i-w_t^i)-hz_t^i\cdot(\nabla_iU(x_{\lfloor t \rfloor})-\nabla_i U(y_{\lfloor t \rfloor}))=b^i(t)+\delta^i(t)
\\ \frac{\rmd}{\rmd t}{b}^i(t)&=-z_t^i\cdot(\nabla_i U(x_{\lfloor t \rfloor})-\nabla_i U(y_{\lfloor t \rfloor}) +\nabla_i U(x_{\lceil t \rceil}) -\nabla_i U(y_{\lceil t \rceil}))
\\ & \indent +2w_{t}^i\cdot w_{\lfloor t \rfloor}^i-hw_t^i\cdot(\nabla_i U(x_{\lfloor t\rfloor})-\nabla_i U(y_{\lfloor t\rfloor}))  
\\ & =-z_{\lfloor t\rfloor}^i\cdot(\nabla_i U(x_{\lfloor t\rfloor})-\nabla_i U(y_{\lfloor t\rfloor}))-z_{\lceil t\rceil}^i\cdot(\nabla_i U(x_{\lceil t\rceil})-\nabla_i U(y_{\lceil t\rceil}))
\\ & \indent +2|w_t^i|^2-2\kappa|z_t^i|^2+\kappa(|z_{\lfloor t \rfloor}^i|^2+|z_{\lceil t \rceil}^i|^2)+\varepsilon^i(t)
\end{align*}
where $\varepsilon^i(t)=\varepsilon_1^i(t)+\varepsilon_2^i(t)+\varepsilon_3^i(t)+\varepsilon_4^i(t) $ and
\begin{align*}
\delta^i(t)&= z_t^i\cdot(2(w_{\lfloor t \rfloor}^i-w_t^i)-h (\nabla_i U(x_{\lfloor t \rfloor })-\nabla_i U(y_{\lfloor t \rfloor })))
\\ \varepsilon_1^i(t)&=-(z_t^i-z_{\lfloor t\rfloor}^i)\cdot(\nabla_i U(x_{\lfloor t\rfloor})-\nabla_i U(y_{\lfloor t\rfloor}))
\\ \varepsilon_2^i(t)&=-(z_t^i-z_{\lceil t\rceil}^i)\cdot(\nabla_i U(x_{\lceil t\rceil})-\nabla_i U(y_{\lceil t\rceil}))
\\ \varepsilon_3^i(t)&=w_t^i\cdot(2(w^i_{\lfloor t\rfloor}-w_t^i)-h(\nabla_i U(x_{\lfloor t\rfloor})-\nabla_i U(y_ {\lfloor t \rfloor})))
\\ \varepsilon_4^i(t)&=\kappa(2|z_t^i|^2-|z_{\lfloor t \rfloor}^i|^2-|z_{\lceil t \rceil}^i|^2). 
\end{align*}
By \eqref{eq:conseq_strongconv} the derivative of $b^i(t)$ is bounded by
\begin{align*}
\frac{\rmd}{\rmd t}{b}^i(t)
 &\leq -2\kappa|z_t^i|^2+\frac{\epsilon\tilde{L}}{n^2}\Big(\sum_{j\neq i} |z_{\lfloor t \rfloor}^j|\Big)^2+\frac{\epsilon\tilde{L}}{n^2}\Big(\sum_{j\neq i} |z_{\lceil t \rceil}^j|\Big)^2
+2|w_t^i|^2
 +\varepsilon^i(t) +2\hat{C}.
\end{align*} 
The previous estimate leads to an initial value problem of the form
\begin{align*}
&\frac{\rmd }{\rmd t}a^i(t)=b^i(t)+\delta^i(t), && a^i(0)=|z_0^i|^2,
\\ & \frac{\rmd }{\rmd t} b^i(t)=-2\kappa a^i(t)+\beta^i(t)+\varepsilon^i(t),&& b^i(0)=0,
\end{align*}
where 
\begin{align}
\beta^i(t)&\leq \frac{\epsilon\tilde{L}}{n^2}\Big(\sum_{j\neq i} |z_{\lfloor t \rfloor}^j|\Big)^2+\frac{\epsilon\tilde{L}}{n^2}\Big(\sum_{j\neq i} |z_{\lceil t \rceil}^j|\Big)^2+2|w_t^i|^2+2\hat{C}. \label{eq:def_beta}
\end{align}
Note that when $h=0$, $\varepsilon^i(t)=\delta^i(t)=0$. 
By variation of parameters, $a^i(t)$ can be written as
\begin{equation}
\begin{aligned}
a^i(t)&=\cos(\sqrt{2\kappa} \ t)|z_0^i|^2+\int_0^t\cos(\sqrt{2\kappa}(t-r))\delta^i(r)\rmd r 
\\ & +\int_0^t\frac{1}{\sqrt{2\kappa}}\sin(\sqrt{2\kappa}(t-r))(\beta^i(r)+\varepsilon^i(r))\rmd r. \label{eq:uniquesolution}
\end{aligned}
\end{equation}
Taylor’s integral formula, i.e., 
$\cos( \sqrt{2 \kappa} \ t) = 1-\kappa t^{2}+(1/6) \int_0^t (t-s)^3 \cos( \sqrt{2 \kappa} \ s) (2 \kappa)^2 ds \le 1 - \kappa t^2 + \kappa^2 t^4 / 6$, and the fact that by \eqref{eq:cond_t_num} and \eqref{eq:KandL} $\kappa^2t^4\leq (L+2\epsilon\tilde{L})^2t^4\leq \kappa t^2$ yield 
\begin{align} \label{eq:taylor_cosinus}
\cos (\sqrt{2\kappa}t)\leq 1-(5/6)\kappa t^2.
\end{align}
Further, we get  by \eqref{eq:cond_t_num} and \eqref{eq:KandL} 
\begin{align} \label{condition_t2}
\kappa  t^2\leq (L+2\epsilon\tilde{L})t^2\leq 1\leq \pi^2/2, \text{ and so } t\leq (\pi/\sqrt{2\kappa}).
\end{align} 
Therefore, $\sin(\sqrt{2\kappa}(t-r))\geq 0$ for all $r\in[0,t]$.
Further, 
\begin{align} \label{eq:sinus_estimate}
\frac{1}{\sqrt{2\kappa}}\sin(\sqrt{2\kappa}(t-r))\leq (t-r).
\end{align}
Inserting \eqref{eq:taylor_cosinus} and \eqref{eq:sinus_estimate} in \eqref{eq:uniquesolution} yields
\begin{equation}
\begin{aligned}
a^i(t)&\leq(1-(5/6)\kappa t^2)|z_0^i|^2+\int_0^t|\delta^i(r)|\rmd r  +\int_0^t(t-r)(|\beta^i(r)|+|\varepsilon^i(r)|)\rmd r. \label{eq:uniquesolution2}
\end{aligned}
\end{equation}

For $\beta^i(t)$, we note that by \eqref{eq:def_beta}, \eqref{eq:velocity_estimate5} with $w^i=0$ and \eqref{eq:cond_t_num},
\begin{align}
|\beta^i(t)|& 
\leq 2\Big((L+2\epsilon\tilde{L})t\frac{5}{4}|z^i_0|+\frac{2\epsilon\tilde{L}t}{n}\frac{5}{4}\max_{s\leq t } \sum_{j\neq i}|z_s^j|\Big)^2+\frac{2\epsilon\tilde{L}}{n^2}\Big(\max_{s\leq t } \sum_{j\neq i} |z_s^j|\Big)^2+2\hat{C} \nonumber
\\ & \leq\frac{25}{4}(L+2\epsilon\tilde{L})^2t^2|z^i_0|^2+\Big(25\frac{\epsilon^2\tilde{L}^2t^2}{n^2}+\frac{2\epsilon\tilde{L}}{n^2}\Big)\Big(\max_{s\leq t } \sum_{j\neq i}|z_s^j|\Big)^2+2\hat{C}. \label{eq:betai}
\end{align}
Note that by \eqref{eq:KandL}, \eqref{eq:cond_t_num} and  since by assumption $\epsilon\tilde{L}< K/6$,
\begin{align} \label{eq:bound_epstildeL}
\epsilon\tilde{L}(t^2+th)\leq (1/10)(K+4\epsilon\tilde{L}) (t^2+th)\leq (1/10)(L+4\epsilon\tilde{L})(t^2+th)\leq 40^{-1}.
\end{align}
Hence, by \eqref{eq:betai} we obtain for the integral containing $\beta^i(t)$ in \eqref{eq:uniquesolution2} 
\begin{align}
\int_0^t &(t-r)|\beta^i(r)|\rmd r \nonumber
\\ &\leq \int_0^t (t-r) \Big(\frac{25}{4}r^2(L+2\epsilon\tilde{L})^2|z^i_0|^2
 +\Big(25\frac{r^2\epsilon^2\tilde{L}^2}{n^2}+\frac{2\epsilon\tilde{L}}{n^2}\Big)\Big(\max_{s\leq t } \sum_{j\neq i}|z_s^j|\Big)^2+2\hat{C}\Big) \rmd r \nonumber
\\ & \leq \frac{25}{48}t^4(L+2\epsilon\tilde{L})^2|z^i_0|^2+ \Big(\frac{25}{12\cdot40}+1\Big)\frac{\epsilon\tilde{L}t^2}{n^2}\Big(\max_{s\leq t }\sum_{j\neq i}|z_s^j|\Big)^2+\hat{C}t^2, \label{eq:integral_beta}
\end{align}
where the last step follows by \eqref{eq:bound_epstildeL}.

Next, we bound $\delta^i(t)$ and $\varepsilon^i(t)$. To bound $\delta^i(t)$ and $\varepsilon_3^i(t)$, we note that by \eqref{eq:hamdyn_ofdifferenceprocess} and \eqref{eq:bound_deltaiU2},
\begin{align*} 
|w^i_{\lfloor t\rfloor}-w_t^i| &\leq \Big|\int_{\lfloor t \rfloor}^t \frac{\rmd}{\rmd s}w_s^i \rmd s\Big|\leq \frac{h}{2} |\nabla_i U(x_{\lfloor t\rfloor})+\nabla_i U(x_{\lceil t\rceil})-(\nabla_i U(y_{\lfloor t \rfloor})+\nabla_i U(y_{\lceil t \rceil}))|
\\ & \leq h\Big((L+2\epsilon\tilde{L}) z_t^{i,*}+\frac{2\epsilon\tilde{L}}{n}\max_{s\leq t }\sum_{j\neq i}|z_s^j|\Big)
\end{align*}
where $z_t^{i,*}=\max_{s\leq t}|z_s^i|$.
Hence, by \eqref{eq:bound_deltaiU2}, \eqref{eq:position_estimate5} with $w^i=0$ and \eqref{eq:cond_t_num},
\begin{align}
|2(w^i_{\lfloor t\rfloor}-w_t^i)-h(\nabla_i U(x_{\lfloor t\rfloor})-\nabla_i U(y_{\lfloor t \rfloor}))|&\leq 3h\Big((L+2\epsilon\tilde{L}) z_t^{i,*}+\frac{2\epsilon\tilde{L}}{n}\max_{s\leq t}\sum_{j\neq i}|z_s^j|\Big) \nonumber
\\ & \leq 3h\Big(\frac{5}{4}(L+2\epsilon\tilde{L})|z_0^i|+\frac{5}{2}\frac{\epsilon\tilde{L}}{n}\max_{s\leq t}\sum_{j\neq i}|z_s^j|\Big). \label{eq:precond_delta^i}
\end{align}
Hence by \eqref{eq:precond_delta^i} and \eqref{eq:position_estimate5} with $w^i=0$, and then by  \eqref{eq:bound_epstildeL} and \eqref{eq:cond_t_num},
\begin{align}
\max_{s\leq t}|\delta^i(s)|&\leq 3h\Big(\frac{5}{4}|z_0^i|+\frac{2\epsilon\tilde{L}(t^2+th)}{n}\max_{s\leq t}\sum_{j\neq i}|z_s^j|\Big)  \nonumber
\\ & \cdot\Big(\frac{5}{4}(L+2\epsilon\tilde{L})|z_0^i|+\frac{5}{2}\frac{\epsilon\tilde{L}}{n}\max_{s\leq t}\sum_{j\neq i}|z_s^j|\Big) \nonumber
\\ & \leq 3h\Big(\frac{25}{16}(L+2\epsilon\tilde{L})|z_0^i|^2+\frac{15}{4}\epsilon\tilde{L}|z_0^i|\frac{1}{n}\max_{s\leq t}\sum_{j\neq i}|z_s^j|+\frac{\epsilon\tilde{L}}{8n}\max_{s\leq t}\Big(\sum_{j\neq i}|z_s^j|\Big)^2 \Big) \label{eq:delta^i_2}
\\ & \leq h\Big(\frac{75L}{16}+15\epsilon\tilde{L}\Big)|z_0^i|^2+h\frac{6\epsilon\tilde{L}}{n^2}\max_{s\leq t}\Big(\sum_{j\neq i}|z_s^j|\Big)^2, \label{eq:delta^i}
\end{align}
Note that Young's product inequality is used in \eqref{eq:delta^i_2} to bound the cross term.
Similarly, by \eqref{eq:precond_delta^i}, \eqref{eq:position_estimate5} with $w^i=0$, \eqref{eq:velocity_estimate5} with $w^i=0$ and \eqref{eq:bound_epstildeL},
\begin{align}
\max_{s\leq t}|\varepsilon^i_3(s)|\frac{t}{2}&\leq 3h\frac{t}{2} \Big(\frac{5}{4}(L+2\epsilon\tilde{L})t|z_0^i|+\frac{5}{4}\frac{2\epsilon\tilde{L}t}{n}\max_{s\leq t} \sum_{j\neq i}|z_s^j|\Big) \nonumber
\\ & \cdot\Big(\frac{5}{4}(L+2\epsilon\tilde{L})|z_0^i|+\frac{5}{2}\frac{\epsilon\tilde{L}}{n}\max_{s\leq t}\sum_{j\neq i}|z_s^j|\Big) \nonumber
\\ & \leq 3h \frac{t^2}{2}\Big(\frac{25}{16}(L+2\epsilon\tilde{L})^2|z_0^i|^2+\frac{25}{4}(L+2\epsilon\tilde{L})|z_0^i|\frac{\epsilon\tilde{L}}{n}\max_{s\leq t}\sum_{j\neq i}|z_s^j| \label{eq:eps^i_3_2}
\\ & \indent +\frac{25}{4}\Big(\frac{\epsilon\tilde{L}}{n}\max_{s\leq t}\sum_{j\neq i}|z_s^j|\Big)^2\Big) \nonumber
\\ & \leq h\frac{75(L+2\epsilon\tilde{L})}{64}|z_0^i|^2+h\frac{15}{32}\frac{\epsilon\tilde{L}}{n^2}\max_{s\leq t}\Big(\sum_{j\neq i}|z_s^j|\Big)^2. \label{eq:eps^i_3}
\end{align}
Note that Young's product inequality is used to bound the cross term in \eqref{eq:eps^i_3_2}.

To bound $\varepsilon_1^i(t)$, $\varepsilon_2^i(t)$ and $\varepsilon_4^i(t)$, we note that by \eqref{eq:hamdyn_ofdifferenceprocess} and \eqref{eq:bound_deltaiU2},
\begin{align} \label{eq:estimateforeps_gamma1}
|z_{\lfloor t \rfloor}^i-z_t^i|&=\Big|\int_{\lfloor t \rfloor }^t \frac{\rmd}{\rmd s} z_s^i \rmd s\Big|  \leq h |w_{\lfloor t \rfloor}^i-\frac{h}{2}(\nabla_i U(x_{\lfloor t \rfloor})-\nabla_i U(y_{\lfloor t \rfloor}))| \nonumber
\\ & \leq hw_t^{i,*}+\frac{h^2}{2}(L+2\epsilon\tilde{L})z_t^{i,*}+\frac{h^2\epsilon\tilde{L}}{n}\max_{s\leq t}\sum_{j\neq i}|z_s^j|
\end{align}
where $w_t^{i,*}=\max_{s\leq t}|w_s^i|$. Similarly, 
\begin{align} \label{eq:estimateforeps_gamma2}
|z_{\lceil t \rceil}^i-z_t^i|& \leq hw_t^{i,*}+\frac{h^2}{2}(L+2\epsilon\tilde{L})z_t^{i,*}+\frac{h^2\epsilon\tilde{L}}{n}\max_{s\leq t}\sum_{j\neq i}|z_s^j|.
\end{align}
Hence, by applying \eqref{eq:estimateforeps_gamma1}, \eqref{eq:estimateforeps_gamma2} and \eqref{eq:bound_deltaiU2} in the first step, and \eqref{eq:position_estimate5} and \eqref{eq:velocity_estimate5} with $w^i=0$ in the second step,
\begin{align} 
\max_{s\leq t} (|\varepsilon^i_1(s)&+\varepsilon^i_2(s)|)\frac{t}{2}\leq th\Big((L+2\epsilon\tilde{L})z_t^{i,*}+\frac{2\epsilon \tilde{L}}{n}\max_{s\leq t} \sum_{j\neq i}|z_s^j|
 \Big) \nonumber
 \\ & \indent \cdot\Big(w_t^{i,*}+\frac{h}{2}(L+2\epsilon\tilde{L})z_t^{i,*}+\frac{h\epsilon\tilde{L}}{n}\max_{s\leq t}\sum_{j\neq i}|z_s^j|\Big) \nonumber
\\ & \leq th \Big(\frac{5}{4}(L+2\epsilon\tilde{L})|z_0^i|+\frac{5\epsilon\tilde{L}}{2n}\max_{s\leq t}\sum_{j\neq i}|z_s^j|\Big) \nonumber
\\ &  \cdot\Big((L+2\epsilon\tilde{L})\Big(t+\frac{h}{2}\Big)\frac{5}{4}|z_0|
 +\frac{2\epsilon\tilde{L}}{n}\frac{5}{4}\Big(t+\frac{h}{2}\Big)\max_{s\leq t}\sum_{j\neq i}|z_s^j|\Big) \nonumber
\\ & \leq ht\Big(t+\frac{h}{2}\Big)\Big((L+2\epsilon\tilde{L})^2\frac{25}{16}|z_0^i|^2+\frac{25}{4}(L+2\epsilon\tilde{L})\epsilon\tilde{L}|z_0^i|\max_{s\leq t}\frac{1}{n}\sum_{j\neq i}|z_s^j| \nonumber
\\&  +\frac{25}{4}\frac{(\epsilon\tilde{L})^2}{n^2}\max_{s\leq t}\Big(\sum_{j\neq i}|z_s^j|\Big)^2\Big) \nonumber
\\ & \leq ht\Big(t+\frac{h}{2}\Big)\Big(\frac{25}{8}(L+2\epsilon\tilde{L})^2|z_0^i|^2+\frac{25(\epsilon\tilde{L})^2}{2n^2}\max_{s\leq t}\Big(\sum_{j\neq i}|z_s^j|\Big)^2 \Big) \nonumber
\\ & \leq h\frac{25}{32}(L+2\epsilon\tilde{L})|z_0^i|^2+h\frac{5\epsilon\tilde{L}}{16n^2}\max_{s\leq t}\Big(\sum_{j\neq i}|z_s^j|\Big)^2 . \label{eq:epsilon12}
\end{align}
Note that Young's product inequality is used to bound the cross term in the third step and \eqref{eq:cond_t_num} and \eqref{eq:bound_epstildeL} are used in the last step.
For $\varepsilon_4^i(t)$, we obtain by \eqref{eq:estimateforeps_gamma1} and \eqref{eq:estimateforeps_gamma2},
\begin{align}
\max_{s\leq t}|\varepsilon_4^i(s)|\frac{t}{2}&  \leq \frac{t}{2}\kappa\max_{s\leq t}| (z_s^i+z_{\lfloor s \rfloor}^i)\cdot(z_s^i-z_{\lfloor t \rfloor}^i)+(z_s^i+z_{\lceil s \rceil}^i)\cdot(z_s^i-z_{\lceil s \rceil}^i)| \nonumber
\\ & \leq 2th\kappa z_t^{i,*}\Big(w_t^{i,*}+\frac{h}{2}(L+2\epsilon\tilde{L})z_t^{i,*}+\frac{h\epsilon\tilde{L}}{n}\max_{s\leq t}\sum_{j\neq i}|z_s^j|\Big) \nonumber
\\ & \leq 2th  \Big(\frac{5}{4}(L+2\epsilon\tilde{L})|z_0^i|+\frac{5\epsilon\tilde{L}}{2n}\max_{s\leq t}\sum_{j\neq i}|z_s^j|\Big) \label{eq:gamma_estimate}
\\ &\indent \cdot\Big((L+2\epsilon\tilde{L})\Big(t+\frac{h}{2}\Big)\frac{5}{4}|z_0^i|
 +\frac{2\epsilon\tilde{L}}{n}\frac{5}{4}\Big(t+\frac{h}{2}\Big)\max_{s\leq t}\sum_{j\neq i}|z_s^j|\Big) \nonumber
 \\ & \leq 2ht\Big(t+\frac{h}{2}\Big)\Big((L+2\epsilon\tilde{L})^2\frac{25}{16}|z_0^i|^2+\frac{25}{4}(L+2\epsilon\tilde{L})\epsilon\tilde{L}|z_0^i|\max_{s\leq t}\frac{1}{n}\sum_{j\neq i}|z_s^j| \nonumber
\\ & \indent +\frac{25}{4}\frac{(\epsilon\tilde{L})^2}{n^2}\max_{s\leq t}\Big(\sum_{j\neq i}|z_s^j|\Big)^2\Big) \nonumber
\\ & \leq 2ht\Big(t+\frac{h}{2}\Big)\Big(\frac{25}{8}(L+2\epsilon\tilde{L})^2|z_0^i|^2+\frac{25(\epsilon\tilde{L})^2}{2n^2}\max_{s\leq t}\Big(\sum_{j\neq i}|z_s^j|\Big)^2 \Big) \nonumber
\\ & \leq h\frac{25}{16}(L+2\epsilon\tilde{L})|z_0^i|^2+h\frac{5\epsilon\tilde{L}}{8n^2}\max_{s\leq t}\Big(\sum_{j\neq i}|z_s^j|\Big)^2  \label{eq:gamma}
\end{align}
where \eqref{eq:gamma_estimate} follows by \eqref{eq:position_estimate5} with $w^i=0$ and \eqref{eq:velocity_estimate5} with $w^i=0$ and since by \eqref{eq:KandL} $\kappa\leq (L+2\epsilon\tilde{L})$. Note that Young's product inequality is used to bound the cross term in the third step.

Therefore, by \eqref{eq:delta^i}, \eqref{eq:eps^i_3}, \eqref{eq:epsilon12} and \eqref{eq:gamma},
\begin{align}
\int_0^t \Big(&(t-r)|\varepsilon^i(r)|+|\delta^i(r)|\Big)\rmd r \nonumber
\\ &  \leq  \frac{t^2}{2}(\max_{s\leq t}(|\varepsilon^i_1(s)+\varepsilon_2^i(s)|)+\max_{s\leq t} |\varepsilon_3^i(s)|+\max_{s\leq t}|\varepsilon^i_4(s)|)+t\max_{s\leq t}|\delta^i(s)| \nonumber
\\ & \leq ht\Big(\frac{25}{32}(L+2\epsilon\tilde{L})+\frac{75}{64}(L+2\epsilon\tilde{L})+\frac{25}{16}(L+2\epsilon\tilde{L})+\frac{75}{16}L+15\epsilon\tilde{L}\Big)|z_0^i|^2 \nonumber
\\ & +ht\Big(\frac{5}{16}+\frac{15}{32}+\frac{5}{8}+6\Big)\frac{\epsilon\tilde{L}}{n^2}\max_{s\leq t}\Big(\sum_{j\neq i}|z_s^j|\Big)^2 \nonumber
\\ & = ht \Big(\frac{525}{64}L+\frac{235}{64}K\Big)|z_0^i|^2 +ht\frac{237}{32}\frac{\epsilon\tilde{L}}{n^2}\max_{s\leq t}\Big(\sum_{j\neq i}|z_s^j|\Big)^2 
 \label{eq:constant_C0_1}
\end{align}
where we used $\epsilon\tilde{L}<K/6$ in \eqref{eq:constant_C0_1}.
We note that by \eqref{eq:cond_h_t}, \eqref{eq:KandL} and since by assumption $\epsilon\tilde{L}<K/6$,
\begin{align} \label{eq:cond_h_t_2}
h\Big(\frac{525}{64}L+\frac{235}{64}K\Big)\leq \frac{Kt}{64}\leq \frac{\kappa t}{32}
\end{align}
and
\begin{align} \label{eq:cond_h_t_3}
h\frac{237}{32}\leq \frac{237}{32}\frac{Kt}{525L+235K} \leq \frac{1}{2}t.
\end{align}
Therefore, by \eqref{eq:constant_C0_1}, \eqref{eq:cond_h_t_2} and \eqref{eq:cond_h_t_3}
\begin{align}
\int_0^t \Big(&(t-r)|\varepsilon^i(r)|+|\delta^i(r)|\Big)\rmd r   
\leq t^2 \Big( \frac{\kappa}{32}|z_0|^2+\frac{\epsilon\tilde{L}}{2n^2}\max_{s\leq t}\Big(\sum_{j\neq i} |z_s^j|\Big)^2\Big). \label{eq:constant_C0}
\end{align}

Inserting \eqref{eq:integral_beta} and \eqref{eq:constant_C0} 
in \eqref{eq:uniquesolution2} and applying \eqref{eq:cond_t_num} yields,
\begin{align*}
a^i(t)&\leq \Big(1-\frac{5}{6} \kappa t^2\Big)|z_0^i|^2+\frac{25}{48} \kappa t^2 |z_0^i|^2+\Big(\frac{5}{96}+1\Big)\frac{\epsilon\tilde{L}t^2}{n^2}\Big(\max_{s\leq t}\sum_{j\neq i} |z_s^j|\Big)^2+\hat{C}t^2
\\& 
+t^2\Big(\frac{\kappa}{32}|z_0^i|^2+\frac{\epsilon\tilde{L}}{2n^2}\Big(\max_{s\leq t}\sum_{j\neq i} |z_s^j|\Big)^2\Big).
\end{align*}
By \eqref{eq:strong_conv_hatC}, 
we obtain for $x,y\in\mathbb{R}^{d n}$ with $|x^i-y^i|>\tilde{R}$,
\begin{align*}
|z_t^i|^2 & \leq \Big(1-\Big(\frac{5}{6}-\frac{25}{48}-\frac{2}{32}\Big)\kappa t^2\Big)|z_0^i|^2
 +\Big(1+\frac{5}{96}+\frac{1}{2}\Big)\frac{\epsilon\tilde{L}t^2}{n^2}\Big(\max_{s\leq t}\sum_{j\neq i} |z_s^j|\Big)^2
\\ &\leq \Big(1-\frac{1}{4} \kappa t^2\Big)|z_0^i|^2
+2\frac{\epsilon\tilde{L}t^2}{n^2}\Big(\max_{s\leq t}\sum_{j\neq i} |z_s^j|\Big)^2,
\end{align*}
as required.

\end{proof}

\section{Proofs of main results}\label{sec:proofs}
\subsection{Proof of main contraction result} \label{section_proof_oftheorems}

For the proof of \Cref{thm:generalcase}, we write ${R}^i$ and $r^i$ for $r^i(x,y)=|x^i-y^i|$ and ${R}^i(x,y)=|{\mathbf{X}}^i(x,y)-{\mathbf{Y}}^i(x,y)|$ for fixed $x,y\in\mathbb{R}^{dn}$.
Further, we write ${r}_s^i=|{q}_s^i(x,\xi)-{q}_s^i(y,\eta)|$ for the distance between the two positions at time $s$ satisfying \eqref{eq:hamdyn_num} 
where $\xi,\eta$ are the velocities coupled using the construction given in \Cref{subsec:coupling}. Further, we denote $z=x-y$ and $w=\xi-\eta$.

\begin{proof} [Proof of \Cref{thm:generalcase}] 
Note that \eqref{cond_T}, \eqref{cond_epsilon} and \eqref{eq:def_kappa} imply
\begin{align} \label{eq:conseq_from_cond_eps}
\kappa \geq (1/2)K \hspace{5mm}\text{and} \hspace{5mm} L+4\epsilon\tilde{L}\leq L+(2K/3)\leq (5/3)L.
\end{align}
Hence, we obtain by \eqref{cond_T}
\begin{align} \label{eq:conditionT_reformulated}
(L+4\epsilon\tilde{L})(T+h_1)^2\leq \min\Big(\frac{1}{4},\frac{\kappa}{L+4\epsilon\tilde{L}},\frac{1}{256(L+4\epsilon\tilde{L})\tilde{R}^2}\Big).
\end{align}
Moreover, the following inequalities are satisfied,
\begin{align}
\gamma T & \leq 1, \label{uneq_1}
\\ (L+4\epsilon\tilde{L})(T+h) &\leq \gamma/4, \label{uneq_2}
\\ \gamma \tilde{R} &\leq 1/4, \label{uneq_3}
\\ \exp(T^{-1}(R_1-\tilde{R})) & \geq 12. \label{uneq_6}
\end{align}
Inequalities \eqref{uneq_1} and \eqref{uneq_3} follow by \eqref{def_gamma}, \eqref{uneq_2} follows by \eqref{def_gamma} and \eqref{eq:conditionT_reformulated},  and the inequality \eqref{uneq_6} follows by \eqref{def_R_1}.

We first prove a bound on $\mathbb{E}[f(R^i)-f(r^i)]$ for each particle $i$ similarly to the strategy to bound $\mathbb{E}[f(R)-f(r)]$ in \cite[Proof of Theorem 2.4]{BoEbZi2020}.  We split the calculation of this expectation in two cases depending on the applied coupling. 

\textbf{Case 1:} $r^i=|x^i-y^i|\geq\tilde{R}$. In this case, the initial velocities of the $i$-th particles are synchronized, i.e., $w^i=0$.
By concavity of the function $f$, by \Cref{lemma:contraction_num} and since
\begin{align} \label{eq:hilfungleichung}
\sqrt{1-\mathsf{a}}\leq 1-\mathsf{a}/2 \hspace{5mm}\text{ for }\mathsf{a}\in[0,1),
\end{align}
we obtain
\begin{align}
\mathbb{E}[f(R^i)-f(r^i)]&\leq f'(r^i)\mathbb{E}[R^i-r^i] \nonumber
\\ & \leq f'(r^i)\Big(-\frac{1}{8}\kappa T^2\Big)r^i
 +f'(r^i)\sqrt{2\epsilon\tilde{L}}\frac{T}{n}\mathbb{E}\Big[\max_{s\leq T}\sum_{j\neq i} r_s^j\Big]. \label{eq:convcase_firstcase}
\end{align}

\textbf{Case 2:} $r^i=|x^i-y^i|< \tilde{R}$.
In this case, since the distance between the $i$-th particles is smaller than $\tilde R$, the initial velocities of the $i$-th particles satisfy $w^i = - \gamma z^i$ with maximal possible probability and otherwise a reflection is applied. These disjoint possibilities motivate splitting the expectation $E[f(R^i)-f(r^i)]$ as follows
\begin{align*}
\mathbb{E}[f(R^i)-f(r^i)]&=\mathbb{E}[f(R^i)-f(r^i),\{w^i=-\gamma z^i\}] 
\\ &+\mathbb{E}[f(R_1\wedge R^i)-f(r^i),\{w^i\neq -\gamma z^i\}] 
\\ &+\mathbb{E}[f(R^i)-f(R_1\wedge R^i),\{w^i\neq-\gamma z^i\}] =  \Romanbar{1}+  \Romanbar{2}+\Romanbar{3}.
\end{align*}

First, we bound the probability $\mathbb{P}[w^i\neq -\gamma z^i]$, which equals the total variation distance between a standard normal distribution with zero mean and a normal distribution with mean $\gamma z^i$ and unit variance, cf. Lemma 4.4 of \cite{BoEb2019}. Note using the coupling characterization of the TV distance, this representation shows that the coupling $\xi^i-\eta^i=-\gamma z^i$ holds with maximal probability.
By (\ref{uneq_3}),
\begin{align}
\mathbb{P}[w^i\neq-\gamma z^i]& 
=\int_{-\infty}^{\gamma|z^i|/2}\frac{1}{\sqrt{2\pi}}\Big(\exp\Big(-\frac{1}{2}\mathsf{x}^2\Big)-\exp\Big(-\frac{1}{2}\Big(\mathsf{x}-\gamma|z^i|\Big)^2\Big)\Big)^+\rmd \mathsf{x}  \nonumber
\\ & =\int_{-\infty}^{\gamma|z^i|/2}\frac{1}{\sqrt{2\pi}}\exp\Big(-\frac{1}{2}\mathsf{x}^2\Big)\rmd \mathsf{x}-\int_{-\infty}^{-\gamma|z^i|/2}\frac{1}{\sqrt{2\pi}}\exp\Big(-\frac{1}{2}\mathsf{x}^2\Big)\rmd \mathsf{x} \nonumber
\\ &\leq \frac{2}{\sqrt{2\pi}}\int_0^{\gamma|z^i|/2} \exp\Big(-\frac{1}{2}\mathsf{x}^2\Big)\rmd \mathsf{x}\leq \frac{2}{\sqrt{2\pi}}\frac{\gamma|z^i|}{2}< \frac{1}{10}. \label{eq:probability_reflected coupling}
\end{align}

Next, we bound $\Romanbar{1}$, $\Romanbar{2}$ and $\Romanbar{3}$.
For $\Romanbar{1}$, we note that on the set $\{w^i=-\gamma z^i\}$, by \eqref{eq:position_estimate4} and \eqref{uneq_2} 
\begin{align*}
R^i&\leq (1-\gamma T)|z^i|+(L+2\epsilon\tilde{L})(T^2+Th)|z^i|
 +\frac{2\epsilon\tilde{L}}{n}(T^2+Th)\max_{s\leq T}\sum_{j\neq i}|z_s^j|
\\ &\leq |z^i|\Big(1-\gamma T +\frac{\gamma T}{4}\Big)+\frac{2\epsilon\tilde{L}}{n}(T^2+Th)\max_{s\leq T} \sum_{j\neq i} r^j_s
\\ & =  \Big(1-\frac{3\gamma T}{4}\Big)r^i+\frac{2\epsilon\tilde{L}}{n}(T^2+Th)\max_{s\leq T} \sum_{j\neq i} r^j_s.
\end{align*}
Hence by concavity of $f$ and by \eqref{eq:probability_reflected coupling},
\begin{align}
\Romanbar{1} & \leq -f'(r^i)\frac{3}{4}\gamma T r^i \mathbb{P}[w^i=-\gamma z^i]  +f'(r^i)\mathbb{E}\Big[\frac{2\epsilon\tilde{L}}{n}(T^2+Th)\max_{s\leq T} \sum_{j\neq i} r^j_s\Big] \nonumber
\\ & \leq -f'(r^i)\frac{27}{40}\gamma T r^i +f'(r^i)\frac{2\epsilon\tilde{L}}{n}(T^2+Th)\mathbb{E}\Big[\max_{s\leq T} \sum_{j\neq i} r^j_s\Big].  \label{eq:exactcase_boundI}
\end{align}
To bound $\Romanbar{2}$, note that by \eqref{eq:definition_f} for $r,s\leq R_1$, 
\begin{align*}
f(s)-f(r)=\int_r^s e^{-t/T}\rmd t =T(e^{-r/T}-e^{-s/T})\leq Te^{-r/T}=Tf'(r). 
\end{align*}
Therefore, by \eqref{eq:probability_reflected coupling}
\begin{align}
\Romanbar{2} &\leq Tf'(r^i)\mathbb{P}[w^i\neq-\gamma z^i]   \leq Tf'(r^i)\frac{\gamma r^i}{\sqrt{2\pi}}< \frac{2}{5}\gamma T r^i f'(r^i). \label{eq:exactcase_boundII}
\end{align}
where we used the bound $1/\sqrt{2 \pi} < 2/5$.
For $\Romanbar{3}$, we get by concavity of $f$
\begin{align}
\Romanbar{3} \leq f'(R_1)\mathbb{E}[(R^i-R_1)^+,\{w^i\neq -\gamma z^i\}]. \label{eq:estimateIII_1}
\end{align}
If $w^i\neq-\gamma z^i$, then $w^i=2(e^i\cdot \xi^i)e^i$ with $e^i=z^i/|z^i|$ and hence $|z^i+Tw^i|=|r^i+2Te^i\cdot \xi^i |$. This computation and \eqref{eq:position_estimate5} yield
\begin{align*}
R^i\leq (1+(L+2\epsilon\tilde{L})(T^2+Th))\max(|r^i+2Te^i\cdot \xi^i |, r^i)+\frac{2\epsilon\tilde{L}(T^2+Th)}{n}\max_{s\leq T} \sum_{j\neq i} r^j_s.
\end{align*}
Hence by \eqref{eq:conditionT_reformulated}
and since $(5/4)r^i-R_1\leq (5/4)\tilde{R}-R_1\leq 0$,
\begin{align}
\mathbb{E}&[(R^i-R_1)^+, \{w^i\neq -z^i\gamma\}] \nonumber
\\ &\leq \mathbb{E}\Big[\Big(\frac{5}{4}\max(|r^i+2Te^i\cdot \xi^i |, r^i)+\frac{2\epsilon\tilde{L}(T^2+Th)}{n}\max_{s\leq T} \sum_{j\neq i} r^j_s-R_1\Big)^+, \{w^i\neq -\gamma z^i\}\Big] \nonumber
\\ &\leq \mathbb{E}\Big[\Big(\frac{5}{4}|r^i+2Te^i\cdot \xi^i|-R_1\Big)^+, \{w^i\neq -\gamma z^i \}\Big]+\mathbb{E}\Big[\frac{2\epsilon\tilde{L}(T^2+Th)}{n}\max_{s\leq T} \sum_{j\neq i} r^j_s\Big]. \label{eq:estimateIII_2}
\end{align}
For the first term, where only the $i$-th particle is involved, we follow the calculations in the proof of \cite[Theorem 2.4]{BoEbZi2020},
\begin{align}
\mathbb{E}\Big[&\Big(\frac{5}{4}|r^i+2Te^i\cdot \xi^i |-R_1\Big)^+, \{w^i\neq -\gamma z^i \}\Big] \nonumber
\\ &=\int_{-\infty}^\infty \Big(\frac{5}{4}|r^i+2 T\mathsf{u}|-R_1\Big)^+\frac{1}{\sqrt{2\pi}}\Big(\exp\Big(-\frac{\mathsf{u}^2}{2}\Big)-\exp\Big(-\frac{(\mathsf{u}+\gamma r^i)^2}{2}\Big)\Big)^+ \rmd \mathsf{u} \nonumber
\\ &= \int_{-\frac{\gamma r^i}{2}}^\infty \Big(\frac{5}{4}|r^i+2 T\mathsf{u}|-R_1\Big)^+\frac{1}{\sqrt{2\pi}}\Big(\exp\Big(-\frac{\mathsf{u}^2}{2}\Big)-\exp\Big(-\frac{(\mathsf{u}+\gamma r^i)^2}{2}\Big)\Big) \rmd \mathsf{u} \nonumber
\\ &=  \int_{-\frac{\gamma r^i}{2}}^{\frac{\gamma z^i}{2}} \Big(\frac{5}{4}|r^i+2 T\mathsf{u}|-R_1\Big)^+\frac{1}{\sqrt{2\pi}}\exp\Big(-\frac{\mathsf{u}^2}{2}\Big)\rmd \mathsf{u} \nonumber
\\ &+ \int_{\frac{\gamma r^i}{2}}^{\infty} \Big(\Big(\frac{5}{4}|r^i+2 T\mathsf{u}|-R_1\Big)^+-\Big(\frac{5}{4}|r^i+2(\mathsf{u}-\gamma r^i) T|-R_1\Big)^+ \Big) \frac{1}{\sqrt{2\pi}}\exp\Big(-\frac{\mathsf{u}^2}{2}\Big)\rmd \mathsf{u} \nonumber
\\ &\leq \int_{\frac{\gamma r^i}{2}}^{\infty} \Big(\frac{5}{4}2\gamma r^i T \Big)\frac{1}{\sqrt{2\pi}}\exp\Big(-\frac{\mathsf{u}^2}{2}\Big)\rmd \mathsf{u} 
\leq \frac{5}{4}\gamma T  r^i . \label{eq:estimateIII_3}
\end{align}
Hence by \eqref{uneq_6}, \eqref{eq:estimateIII_1}, \eqref{eq:estimateIII_2} and \eqref{eq:estimateIII_3},
\begin{align}
\Romanbar{3} & \leq f'(R_1)\Big(\frac{5}{4}\gamma T r^i+\mathbb{E}\Big[\frac{2\epsilon\tilde{L}(T^2+Th)}{n}\max_{s\leq T} \sum_{j\neq i} r^j_s\Big]\Big)   \nonumber
\\& \leq f'(r^i)\Big(\frac{5}{48}\gamma T r^i+\frac{\epsilon\tilde{L}(T^2+Th)}{6n}\mathbb{E}\Big[\max_{s\leq T} \sum_{j\neq i} r^j_s\Big]\Big). \label{eq:exactcase_boundIII}
\end{align}
We combine the bounds on \Romanbar{1}, \Romanbar{2} and \Romanbar{3} in \eqref{eq:exactcase_boundI}, \eqref{eq:exactcase_boundII} and \eqref{eq:exactcase_boundIII} respectively, to obtain for $r^i\leq\tilde{R}$,
\begin{align}
\mathbb{E}[f(R^i)-f(r^i)]&\leq  -f'(r^i)\frac{27}{40}\gamma T r^i +f'(r^i)\frac{2}{5}\gamma T r^i +f'(r^i)\frac{5}{48}\gamma T  r^i\nonumber
\\ &+f'(r^i)\Big(\frac{2\epsilon\tilde{L}(T^2+Th)}{n}+\frac{\epsilon\tilde{L}(T^2+Th)}{6n}\Big)\mathbb{E}\Big[\max_{s\leq T} \sum_{j\neq i} r^j_s\Big] \nonumber
\\ & \leq -f'(r^i)\frac{41}{240}\gamma T r^i + f'(r^i)\frac{13\epsilon\tilde{L}(T^2+Th)}{6n}\mathbb{E}\Big[\max_{s\leq T} \sum_{j\neq i} r^j_s\Big]. \label{eq:convcase_secondcase}
\end{align}

Next, we combine \eqref{eq:convcase_firstcase} and \eqref{eq:convcase_secondcase} and sum over $i$ to obtain
\begin{equation} \label{eq:estimate_expectedsum_1}
\begin{aligned}
\mathbb{E}\Big[\sum_i (f(R^i)-f(r^i))\Big] &\leq -\min\Big(\frac{41}{240}\gamma T ,\frac{1}{8}\kappa T^2\Big)\sum_ir^i f'(r^i)
\\ & + \max\Big(\frac{13\epsilon\tilde{L}(T^2+Th)}{6},\sqrt{2\epsilon\tilde{L}}T\Big)\frac{1}{n}\sum_if'(r^i)\mathbb{E}\Big[\max_{s\leq T} \sum_{j\neq i} r^j_s\Big].
\end{aligned}
\end{equation}
To bound the expectation in the last term of \eqref{eq:estimate_expectedsum_1} we note that
when $w^j \ne -\gamma z^j$, then $w^j = 2 (e^j \cdot \xi^j) e^j$ with $e^j = z^j/|z^j|$, and hence by \eqref{uneq_3},
\begin{align} 
\mathbb{E}\Big[|w^j|&1_{\{r^j<\tilde{R}\}\cap \{w^j\neq -\gamma z^j\}}\Big] \nonumber
\\ & = 1_{\{r^j\leq \tilde{R}\}}\int_{-\frac{\gamma r^j}{2}}^\infty\frac{2|\mathsf{x}|}{\sqrt{2\pi}}\Big(\exp\Big(-\frac{\mathsf{x}^2}{2}\Big)-\exp\Big(-\frac{(\mathsf{x}+\gamma r^j)^2}{2}\Big)\Big) \rmd \mathsf{x} \nonumber
\\ & =1_{\{r^j\leq \tilde{R}\}}\Big(\int_{-\frac{\gamma r^j}{2}}^\infty\frac{1}{\sqrt{2\pi}}2|\mathsf{x}|\exp\Big(-\frac{\mathsf{x}^2}{2}\Big) \rmd \mathsf{x}-\int_{\frac{\gamma r^j}{2}}^\infty\frac{1}{\sqrt{2\pi}}2|\mathsf{x}-\gamma r^j|\exp\Big(-\frac{\mathsf{x}^2}{2}\Big) \rmd \mathsf{x}\Big) \nonumber
\\ & \leq 1_{\{r^j\leq \tilde{R}\}}\Big(\int_{-\frac{\gamma r^j}{2}}^{\frac{\gamma r^j}{2}}\frac{1}{\sqrt{2\pi}}2|\mathsf{x}|\exp\Big(-\frac{\mathsf{x}^2}{2}\Big) \rmd \mathsf{x}+\int_{\frac{\gamma r^j}{2}}^\infty\frac{1}{\sqrt{2\pi}}2\gamma r^j\exp\Big(-\frac{\mathsf{x}^2}{2}\Big) \rmd \mathsf{x}\Big) \nonumber
\\ & \leq 1_{\{r^j\leq \tilde{R}\}}\Big((\gamma r^j)^2+ \gamma r^j\Big) \leq 1_{\{r^j\leq \tilde{R}\}}\Big(\frac{1}{4} \gamma r^j+ \gamma r^j\Big) 
\leq \frac{5}{4}\gamma  r^j. \label{eq:estimate_expectatedsum2}
\end{align}
Then we obtain by \eqref{eq:position_estimate6}, by \eqref{eq:conditionT_reformulated},
and since by \eqref{uneq_1} for $w^j=-\gamma z^j$, $|z^j+Tw^j|\leq |z^j|$,
\begin{align}
\mathbb{E}\Big[\max_{s\leq T} \sum_{j\neq i} r^j_s\Big]& \leq \frac{5}{4}\mathbb{E}\Big[\sum_{j}\max( |z^j+Tw^j|,|z^j|)\Big] \nonumber
\\ &\leq \frac{5}{4}\mathbb{E}\Big[\sum_{j} |z^j|+\sum_j T|w^j|1_{\{r^j<\tilde{R}\}\cap \{w^j\neq -\gamma z^j\}}\Big] \nonumber
\\ & = \frac{5}{4}\sum_{j} |z^j|+\frac{5}{4}T\sum_{j}\mathbb{E}\Big[|w^j|1_{\{r^j<\tilde{R}\}\cap \{w^j\neq -\gamma z^j\}}\Big]
\leq \frac{45}{16}\sum_{j} r^j, \label{eq:estimate_expectatedsum}
\end{align}
where last step holds by \eqref{eq:estimate_expectatedsum2} and \eqref{uneq_1}.
Hence inserting \eqref{eq:estimate_expectatedsum} in \eqref{eq:estimate_expectedsum_1}, 
\begin{equation} \label{eq_expectationsumR-r} \begin{aligned} 
\mathbb{E}\Big[\sum_i (f(R^i)-f(r^i))\Big] &\leq -\min\Big(\frac{41}{240}\gamma T ,\frac{1}{8}\kappa T^2\Big)\sum_i f'(r^i) r^i 
\\ & + \sum_i f'(r^i)\max\Big(\frac{13\epsilon\tilde{L}(T^2+Th)}{6},\sqrt{2\epsilon\tilde{L}}T\Big)\frac{1}{n}\frac{45}{16}\sum_j r^j. \end{aligned}
\end{equation}
Since by \eqref{uneq_2} $\kappa T^2\leq T\gamma/4$, the minimum in \eqref{eq_expectationsumR-r} is attained at $\frac{1}{8}\kappa T^2$. Since \eqref{eq:KandL}, \eqref{eq:conditionT_reformulated} and \eqref{cond_epsilon} imply \eqref{eq:bound_epstildeL} with $t=T$, it holds that $(13/6)\epsilon\tilde{L}(T^2+Th)\leq 13/(6\sqrt{40})\sqrt{\epsilon\tilde{L}(T^2+Th)}$.
Hence, the maximum in \eqref{eq_expectationsumR-r} is attained at $\sqrt{2\epsilon\tilde{L}}T$.
The minimum of $\frac{r^if'(r^i)}{f(r^i)}$ is attained at $R_1$ defined in \eqref{def_R_1},
\begin{align} \label{eq:generalcase_minimumf(r^i)}
\inf_{r^i}\frac{r^if'(r^i)}{f(r^i)}=\frac{R_1\exp(-R_1/T)}{T(1-\exp(-R_1/T))}\ge \frac{5}{4}\Big(\frac{\tilde{R}}{T}+2\Big)\exp\Big(-\frac{5\tilde{R}}{4T}\Big)\exp\Big(-\frac{5}{2}\Big),
\end{align}
and it holds by \eqref{def_R_1} that
\begin{align}
\sum_i f'(r^i) \frac{1}{n}r^j\leq \frac{f(r^j)}{ f'(r^j)}\leq\exp\Big(\frac{R_1}{T}\Big)f(r^j)= \exp\Big(\frac{5\tilde{R}}{4T}\Big)\exp\Big(\frac{5}{2}\Big)f(r^j) \label{eq:generalcase_maximumf(r^i)}
\end{align} 
where we used that $f(r^j) \ge r^j f^{\prime}( r^j)$ and $\exp(-R_1/T)\leq f^{\prime}(r^i) \le 1$. 
Hence,
\begin{align*}
\mathbb{E}\Big[\sum_i (f(R^i)-f(r^i))\Big] &\leq -\frac{1}{8}\kappa T^2\frac{5}{4}\Big(\frac{\tilde{R}}{T}+2\Big)\exp\Big(-\frac{5\tilde{R}}{4T}\Big)\exp\Big(-\frac{5}{2}\Big) \sum_i f(r^i) 
\\&+\sqrt{2\epsilon\tilde{L}}T \frac{45}{16}\exp\Big(\frac{5}{2}\Big)\exp\Big(\frac{5\tilde{R}}{4T}\Big)\sum_if(r^i)
\\ &\leq -\frac{1}{78} \kappa T^2\exp\Big(-\frac{5\tilde{R}}{4T}\Big)\sum_i f(r^i),
\end{align*}
where the last step holds by \eqref{cond_epsilon}.
\end{proof}

\subsection{Proofs of results from Section \ref{subsection:dimesionfree_distance_bounds} }\label{section:proofs_errorbounds}

\begin{proof} [Proof of \Cref{cor_Wassersteindist}]
This proof works analogously to the proof of \cite[Corollary 2.6]{BoEbZi2020} and uses essentially \cite[Lemma 6.1]{BoEbZi2020}. By \Cref{thm:generalcase}, the contractivity condition 
\begin{align} \label{eq:contractivity condition}
\mathbb{E}[\rho(\mathbf{X}(x,y),\mathbf{Y}(x,y))]\leq e^{-c}\rho(x,y)
\end{align}
is satisfied for the coupling $(\mathbf{X}(x,y),\mathbf{Y}(x,y))$. Let $\nu,\eta$ be probability measures on $\mathbb{R}^{d n}$ and let $\omega$ be an arbitrary coupling of $\nu$ and $\eta$. By \cite[Lemma 6.1]{BoEbZi2020}, there exists a Markov chain $(\mathbf{X}_m,\mathbf{Y}_m)_{m\geq 0}$ on a probability space $(\tilde{\Omega},\tilde{\mathcal{A}},\tilde{P})$ such that $(\mathbf{X}_0,\mathbf{Y}_0)\sim \omega$, $(\mathbf{X}_m)$, $(\mathbf{Y}_m)$ are Markov chains each having transition kernel $\pi_h$ and initial distributions $ \nu$ and $\eta$, respectively, and $M_m=e^{cm}\rho (\mathbf{X}_m,\mathbf{Y}_m)$ is a non-negative supermartingale. Then, for all $m\in\mathbb{N}$,
\begin{align*}
\mathcal{W}_\rho (\nu{\pi_h}^m,\eta{\pi_h}^m)\leq \mathbb{E}[\rho(\mathbf{X}_m,\mathbf{Y}_m)]\leq e^{-cm}\mathbb{E}[\rho(\mathbf{X}_0,\mathbf{Y}_0)]=e^{-cm}\int \rho \rmd \omega.
\end{align*}
Since $\omega$ is chosen arbitrary, we take the infimum over all couplings $\omega\in\Gamma(\nu,\eta)$ and obtain \eqref{eq:contraction_Wasserstein_1}. The bound \eqref{eq:contraction_Wasserstein_2} follows by
\eqref{eq:equiv_metric}.
The existence of a unique probability measure $\mu_h$ on $\mathbb{R}^{dn}$ holds by \eqref{eq:contraction_Wasserstein_2} and by Banach fixed-point theorem, cf. \cite[Theorem 3.9]{eberle2020}. Since $\mu_h{\pi_h}^m=\mu_h$ for all $m$, $\Delta(m)\leq e^{R_1/T} e^{-cm} \Delta(0)$. Hence, for a given $\tilde{\epsilon}>0$, $\Delta(m)\leq \tilde{\epsilon}$ holds for \eqref{eq_stepnumber} by \eqref{def_R_1}.
\end{proof}

\begin{proof}[Proof of \Cref{thm:flowdifference_exact_num}]
This proof uses essentially standard numerical analysis techniques and a priori estimates given in \Cref{lemma:positionestimate}.
Fix $x,\xi\in\mathbb{R}^{d n}$. Denote by $(x_s,v_s)=(q_s(x,\xi),p_s(x,\xi))$ the Hamiltonian dynamics driven by \eqref{hamiltoniandynamics}. Set $\mathbf{x}^i_{k}:=q^i_{kh}(x,\xi)$, $\tilde{\mathbf{x}}^i_{k}:=\tilde{q}^i_{kh}(x,\xi)$,  $\mathbf{v}^i_{k}:=p^i_{kh}(x,\xi)$ and $\tilde{\mathbf{v}}^i_{k}:=\tilde{p}^i_{kh}(x,\xi)$. By \eqref{hamiltoniandynamics} and \eqref{eq:hamdyn_num}, it holds
\begin{align*}
|\mathbf{x}^i_{k+1} -\tilde{\mathbf{x}}^i_{k+1}|&\leq |\mathbf{x}^i_{k} -\tilde{\mathbf{x}}^i_{k}|+h|\mathbf{v}_k^i-\tilde{\mathbf{v}}_k^i |+\Big|\int_{kh}^{(k+1)h}\int_{kh}^u \Big(\nabla_i U(x_r)-\nabla_i U (\tilde{\mathbf{x}}_k)\Big) \rmd r \rmd u \Big|,
\\ |\mathbf{v}^i_{k+1} -\tilde{\mathbf{v}}^i_{k+1}|&\leq |\mathbf{v}_k^i-\tilde{\mathbf{v}}_k^i |+\Big|\int_{kh}^{(k+1)h} \Big( \frac{1}{2}\nabla_i U (\tilde{\mathbf{x}}_k)-\nabla_i U(x_u)+\frac{1}{2}\nabla_i U (\tilde{\mathbf{x}}_{k+1}) \Big) \rmd u \Big|.
\end{align*}
By \eqref{eq:bound_deltaiU2} and \eqref{eq:hamdyn_num},
\begin{align*}
\sum_i\Big|\nabla_i U(x_r)-\nabla_i U (\tilde{\mathbf{x}}_k) \Big| &\leq \sum_i\Big|\nabla_i U(x_r)-\nabla_i U ({\mathbf{x}}_k) \Big|+\sum_i\Big|\nabla_i U(\mathbf{x}_k)-\nabla_i U (\tilde{\mathbf{x}}_k) \Big|
\\ & \leq \sum_i\Big|\int_{kh}^r v_s\cdot\nabla\nabla_i U(x_s)\rmd s\Big| +(L+4\epsilon\tilde{L})\sum_i |\mathbf{x}^i_{k} -\tilde{\mathbf{x}}^i_{k}|
\\ & \leq \sum_i \Big|\int_{kh}^r (L+4\epsilon\tilde{L})v_s^i \rmd s\Big|+(L+4\epsilon\tilde{L})\sum_i |\mathbf{x}^i_{k} -\tilde{\mathbf{x}}^i_{k}|
\\& \leq \sum_i (L+4\epsilon\tilde{L})\Big(h\Big(\frac{21}{16}|v_0^i|+\frac{5}{4}(L+4\epsilon\tilde{L})T|x_0^i|\Big) + |\mathbf{x}^i_{k} -\tilde{\mathbf{x}}^i_{k}|\Big),
\end{align*}
where \eqref{eq:velocity_estimate3} and $(L+4\epsilon\tilde{L})T^2\leq (1/4)$ is used in the last step. Analogously,
\begin{align}
\sum_i\Big(&-\nabla_i U(x_u)+\frac{1}{2}\nabla_i U (\tilde{\mathbf{x}}_k)+\frac{1}{2}\nabla_i U (\tilde{\mathbf{x}}_{k+1})\Big) \nonumber
\\ & \leq \sum_i (L+4\epsilon\tilde{L})\Big(h\Big(\frac{21}{16}|v_0^i|+\frac{5}{4}(L+4\epsilon\tilde{L})T|x_0^i|\Big)
 + \frac{1}{2}|\mathbf{x}^i_{k} -\tilde{\mathbf{x}}^i_{k}| + \frac{1}{2}|\mathbf{x}^i_{k+1} -\tilde{\mathbf{x}}^i_{k+1}|\Big). \label{eq:trapez_ineq_notapplied}
\end{align} 
Then for any initial position $x\in\mathbb{R}^{d n}$,
\begin{align*}
\mathbb{E}\Big[\sum_i |\mathbf{x}^i_{k+1} -\tilde{\mathbf{x}}^i_{k+1}|\Big]&\leq \Big(1+\frac{h^2(L+4\epsilon\tilde{L})}{2}\Big)\mathbb{E}\Big[\sum_i |\mathbf{x}^i_{k} -\tilde{\mathbf{x}}^i_{k}|\Big]  
\\&  +h\mathbb{E}\Big[\sum_i |\mathbf{v}_k^i-\tilde{\mathbf{v}}_k^i | \Big]
+\frac{h^3}{2}M_1,  
\end{align*}
and
\begin{equation}\label{eq:trapez_ineq_notapplied2}
\begin{aligned}
\mathbb{E}\Big[\sum_i |\mathbf{v}^i_{k+1} &-\tilde{\mathbf{v}}^i_{k+1}|\Big]\leq \mathbb{E}\Big[\sum_i |\mathbf{v}^i_{k} -\tilde{\mathbf{v}}^i_{k}|\Big]+h^2M_1 \\ & +\frac{(L+4\epsilon\tilde{L})h}{2}\Big(\mathbb{E}\Big[\sum_i |\mathbf{x}^i_{k+1} -\tilde{\mathbf{x}}^i_{k+1}|\Big] +
\mathbb{E}\Big[\sum_i |\mathbf{x}^i_{k} -\tilde{\mathbf{x}}^i_{k}|\Big]\Big)
\end{aligned}
\end{equation}
with $M_1:=\mathbb{E}_{\xi\sim \mathcal{N}(0,I_{d n})}[\sum_i (L+4\epsilon\tilde{L})(\frac{21}{16}|\xi^i|+\frac{5}{4}(L+4\epsilon\tilde{L})T|x^i|)]$.
Set $a_k:=\mathbb{E}[\sum_i |\mathbf{x}^i_{k} -\tilde{\mathbf{x}}^i_{k}|]$ and $b_k:=\mathbb{E}[\sum_i |\mathbf{v}^i_{k} -\tilde{\mathbf{v}}^i_{k}|]$. 
The goal is to bound $a_k$ from above using the discrete Gronwall lemma \cite[Proposition 3.2]{emmrich99}.
Note that this sequence $(a_k,b_k)$ with $a_0=b_0=0$ 
satisfies
\begin{align*}
{a}_{k+1}&\leq (1+(L+4\epsilon\tilde{L})h^2/2){a}_k+h{b}_k+(h^3 M_1/2) 
\\ {b}_{k+1}&\leq{b}_k+h^2M_1+((L+4\epsilon\tilde{L})h/2)({a}_{k+1}+{a}_{k}).
\end{align*}
We deduce for ${b}_{k+1}$
\begin{align*}
{b}_{k+1}&\leq (L+4\epsilon\tilde{L})h\sum_{l=1}^k {a}_l+ \frac{(L+4\epsilon\tilde{L})h}{2}{a}_{k+1} +(k+1)h^2M_1.
\end{align*}
Inserting this estimate in ${a}_{k+1}$ yields
\begin{align}
{a}_{k+1}&\leq (1+(L+4\epsilon\tilde{L})h^2){a}_k+(kh^3M_1+h^3M_1/2)+(L+4\epsilon\tilde{L})h^2\sum_{l=1}^{k-1}{a}_l. \label{eq:sequence_ak}
\end{align}
Note that the sequence $(\tilde{a}_k)$ satisfying 
\begin{align}
\tilde{a}_{k+1}&= (1+(L+4\epsilon\tilde{L})h^2)\tilde{a}_k+(k+(1/2))h^3M_1+(L+4\epsilon\tilde{L})h^2\sum_{l=1}^{k-1}\tilde{a}_l
\end{align}
is an upper bound of the sequence $(a_k)$, i.e. $a_k\leq \tilde{a}_k$. Moreover, it holds $\tilde{a}_{k}\leq \tilde{a}_{k+1}$. Hence,
\begin{align*}
\tilde{a}_{k+1}\leq (1+(L+4\epsilon\tilde{L})kh^2)\tilde{a}_k+(k+1/2)h^3M_1 \leq (1+(L+4\epsilon\tilde{L})Th)\tilde{a}_k+Th^2M_1.
\end{align*}
Applying the discrete Gr\"{o}nwall lemma to $\tilde{a}_k$ yields for all $k\leq (T/h)$,
\begin{align}
a_k \leq \tilde{a}_k&\leq \frac{1}{(L+4\epsilon\tilde{L})T}\Big((1+(L+4\epsilon\tilde{L})hT)^k-1\Big)ThM_1 \nonumber
\\& \leq h \frac{\exp((L+4\epsilon\tilde{L})T^2)-1}{(L+4\epsilon\tilde{L})}M_1\leq h \frac{\exp(1/4)-1}{(L+4\epsilon\tilde{L})}M_1,  \label{eq:sequence_a_k}
\end{align}
where we applied $(L+4\epsilon\tilde{L})T^2\geq 1/4$ in the last step.

Hence, there exists a constant $C_2$ depending on $L$, $\tilde{L}$, $\epsilon$ 
and $T$ such that for all $k\in \mathbb{N}$ with $kh\leq T$ and for any initial value $x\in\mathbb{R}^{d n}$,
\begin{align*}
\mathbb{E}\Big[\sum_i |\mathbf{x}_k^i-\tilde{\mathbf{x}}_k^i|\Big]\leq  h\cdot C_2\Big(d^{1/2}n+\sum_i|x^i|\Big)
\end{align*}
and so \eqref{eq:first_estimate} holds. Note that the term $d^{1/2}n$ comes from $\mathbb{E}[\sum|\xi^i|]$ since $\xi^i\sim\mathcal{N}(0,I_d)$.

If we assume additionally \Cref{ass_V_thirdfourthder} and \Cref{ass_W_thirdfourthder}, then we can instead of \eqref{eq:trapez_ineq_notapplied} bound using \eqref{eq:bound_deltaiU} and the trapezoidal rule,
\begin{align}
\Big|&\int_{kh}^{(k+1)h}\sum_i\Big(-\nabla_i U(x_u)+\frac{1}{2}\nabla_i U (\tilde{\mathbf{x}}_k)+\frac{1}{2}\nabla_i U (\tilde{\mathbf{x}}_{k+1})\Big)\rmd u\Big| \nonumber
\\ & \leq \Big|\int_{kh}^{(k+1)h}\sum_i\Big(-\nabla_i U(x_u)+\frac{1}{2}\nabla_i U ({\mathbf{x}}_k)+\frac{1}{2}\nabla_i U ({\mathbf{x}}_{k+1})\Big)\rmd u\Big| \nonumber
\\ & +\frac{h}{2}\sum_i(L+4\epsilon\tilde{L})( |{\mathbf{x}}_{k}-\tilde{\mathbf{x}}_{k}|+|{\mathbf{x}}_{k+1}-\tilde{\mathbf{x}}_{k+1}|) \nonumber
\\ & \leq \frac{h}{2}\sum_i(L+4\epsilon\tilde{L})(|\mathbf{x}_k^i-\tilde{\mathbf{x}}_k^i|+|\mathbf{x}_{k+1}^i-\tilde{\mathbf{x}}_{k+1}^i|)+\frac{h^3}{12}\sum_i\sup_{u\in[kh,(k+1)h]}\Big|\frac{\rmd^2}{\rmd u^2}\nabla_i U(x_u)\Big|.  \label{eq:trapeziodal_bound}
\end{align}
The last term is bounded using \eqref{eq:hamdyn_num}, \eqref{eq:bound_deltaiU}, \Cref{ass_V_thirdfourthder} and \Cref{ass_W_thirdfourthder} by
\begin{align*}
\sum_i\sup_{u\in[kh,(k+1)h]}\Big|\frac{\rmd^2}{\rmd u^2}\nabla_i U(x_u)\Big|\leq \sum_i (L_H+8\epsilon\tilde{L}_H)\max_{s\leq T}|v_s^i|^2+\sum_i (L+4\epsilon\tilde{L})^2\max_{s\leq T}|x_s^i|.
\end{align*}
Since we can bound $\sum_i\max_{s\leq T}|v_s^i|^2$ and $\sum_i\max_{s\leq T}|x_s^i|$ by \Cref{lemma:positionestimate} and Young's product inequality in terms of $\sum_i|\xi^i|$, $\sum_i|\xi^i|^2$, $\sum_i|x^i|$ and $\sum_i|x^i|^2$, we can bound the last term in \eqref{eq:trapeziodal_bound} after taking expectation over $\xi\sim\mathcal{N}(0,I_{dn})$ by a constant $h^3M_2$ where $M_2$ is a constant depending on $L$, $\tilde{L}$, $L_H$, $\tilde{L}_H$, $\epsilon$, $d$, $n$, $\sum_i|x^i|$ and  $\sum_i|x^i|^2$. More precisely, the dependence of $M_2$ is linear in $nd$, $\sum_i|x^i|$ and  $\sum_i|x^i|^2$.
Replacing $h^2M_1$ in \eqref{eq:trapez_ineq_notapplied2} by $h^3M_2$ leads to the fact that ${a}_k$ in \eqref{eq:sequence_a_k} is bounded from above by 
${a}_{k+1}\leq h^2(\exp(1/(4k))-1)/(L+4\epsilon\tilde{L})(M_2+M_1/(2T)).$
Hence, there exists a constant $\tilde{C}_2$ of order $\mathcal{O}(T^{-1})$ depending on $L$, $\tilde{L}$, $\epsilon$, $L_H$ and $\tilde{L}_H$ such that for all $k\in\mathbb{N}$ with $kh\leq T$ and for any initial value $x\in\mathbb{R}^{dn}$ \eqref{eq:second_estimate} holds, which concludes the proof.
\end{proof}

\begin{proof} [Proof of \Cref{thm:error_invmeas_stepnumber}]
Let $\nu$ be an arbitrary probability measure on $\mathbb{R}^{dn}$.
Recall that by \Cref{cor_Wassersteindist}, it holds $\mathcal{W}_{\ell^1}(\mu_h{\pi_h}^m,\nu{\pi_h}^m)\leq  \exp((5/4)(2+(\tilde{R}/T)))\exp(-cm)\mathcal{W}_{\ell^1}(\mu_h,\nu)$. By \eqref{eq:equiv_metric} and \Cref{cor_Wassersteindist},
\begin{align*}
\Delta(m)&:=\mathcal{W}_{\ell^1}(\mu,{\nu}{\pi_h}^m)\leq \mathcal{W}_{\ell^1}(\mu,\mu_h)+\mathcal{W}_{\ell^1}(\mu_h,{\nu}{\pi_h}^m)\leq \Romanbar{1}+\Romanbar{2}, \hspace{1cm} \text{where}
\\ \Romanbar{1} & = \exp\Big(\frac{5}{4}\Big(2+\frac{\tilde{R}}{T}\Big)\Big)\mathcal{W}_\rho(\mu,\mu_h)
\\ \Romanbar{2} & = \exp\Big(\frac{5}{4}\Big(2+\frac{\tilde{R}}{T}\Big)-cm\Big)\mathcal{W}_{\ell^1}(\mu_h,\nu).
\end{align*}
For $m$ chosen as in \eqref{eq:estimateofstepnumbers_num},  $\Romanbar{2}\leq \tilde{\epsilon}/2$.
To obtain $\Romanbar{1}\leq \tilde{\epsilon}/2$, we use the results of \Cref{cor:bound_mu}.
Then there exists $h_2$ such that for $h\leq \min(h_1,h_2)$, $\Romanbar{1}\leq \tilde{\epsilon}/2$ holds. In particular, we choose $h_2^{-1}=2C_2(d^{1/2}n+\int\sum_i|x^i|\mu(\rmd x))/(c\tilde{\epsilon})$. Hence, for fixed $L$, $\tilde{L}$, $\epsilon$, $K$, $R$, $T$, $h_2^{-1}$ is of order $\mathcal{O}(\tilde{\epsilon}^{-1}( d^{1/2}n+\int\sum_i|x^i|\mu(\rmd x)))$. 
If additionally \Cref{ass_V_thirdfourthder} and \Cref{ass_W_thirdfourthder} are assumed, then for $h\leq \min(h_1,\tilde{h}_2)$ where $\tilde{h}_2^{-1}=(2\tilde{C}_2(dn+\int\sum_i|x^i|\mu(\rmd x)+\int\sum_i|x^i|^2\mu(\rmd x))/(c\tilde{\epsilon}))^{1/2}$, $\Romanbar{1}\leq \tilde{\epsilon}/2$ holds. Note that $\tilde{h}_2^{-1}$ is for fixed $L$, $\tilde{L}$, $L_H$, $\tilde{L}_H$, $\epsilon$, $K$, $R$, $T$ of order $\mathcal{O}(\tilde{\epsilon}^{-1/2}((nd)^{1/2}+(\int\sum_i|x^i|\mu(\rmd x))^{1/2}+(\int\sum_i|x^i|^2\mu(\rmd x))^{1/2}))$.

Let us finally remark that $\sum|x^i|\mu(\rmd x)=\int |x^1|\mu(\rmd x)$ and $\int\sum_i|x^i|^2\mu(\rmd x)$ are finite. This holds, since by \Cref{ass_V_strongconv} and \Cref{ass_W_lipschitz} $\exp(-U(x))$ can be bounded from above by a density function of a Gaussian product measure which has finite first and second moments.
\end{proof}

\subsection{Proofs of results from Section \ref{sec:dimension_free_bounds} }\label{section:proofs_dim_free_bounds}

\begin{proof}[Proof of \Cref{thm:ergodicaverages}] The proof follows \cite[Proof of Theorem 3.17]{eberle2020}. 
It holds for $m,b\in\mathbb{N}$ by \eqref{eq:estimator},
\begin{align*}
\mathbb{E}_\nu[A_{m,b} g]=\frac{1}{m}\sum_{k=b}^{b+m-1} (\nu{\pi_h}^k)(g).
\end{align*}
For all $g\in\mathcal{C}_b^1(\mathbb{R}^{nd})$ with $\max_{l\in \{1,...,n\}}\|\nabla_l g \|\leq \infty$,
\begin{align*}
|g(x)-g(y)|&=\sum_{i=1}^n |g(x^1,...,x^i,y^{i+1},...,y^n)-g(x^1,...,x^{i-1},y^{i},...,y^n)| 
\\ & \leq \max_l \|\nabla_l g\| \sum_{i=1}^n |(x^1,...,x^i,y^{i+1},...,y^n)-(x^1,...,x^{i-1},y^{i},...,y^n)|\\ & =\max_l \|\nabla_l g\| \sum_{i=1}^n|x^i-y^i|.
\end{align*}
Then for all $k\in\mathbb{N}$ and for all couplings $\omega\in\Gamma(\nu\pi_h^k,\mu)$,
\begin{align*}
|(\nu\pi_h^k)(g)-\mu(g)|\leq \max_l \|\nabla_l g\| \int\sum_{i=1}^n|x^i-y^i|\omega(\rmd x \rmd y).
\end{align*}
Hence by the triangle inequality, by \eqref{eq:contraction_Wasserstein_3} and by \eqref{eq:equiv_metric},
\begin{align*}
|\mathbb{E}_\nu [A_{m,b}g]&-\mu(g)|
\\&\leq \frac{1}{m}\sum_{k=b}^{b+m-1}|(\nu{\pi_h}^k)(g)-\mu(g)|
 \leq \frac{1}{m}\sum_{k=b}^{b+m-1}  \max_i\|\nabla_i g\|_\infty \mathcal{W}_{\ell^1}(\nu{\pi_h}^k,\mu)
\\ & \leq \frac{1}{m}\sum_{k=b}^{b+m-1}  \max_i\|\nabla_i g\|_\infty \mathcal{W}_{\ell^1}(\nu{\pi_h}^k,\mu_h)+ \max_i\|\nabla_i g\|_\infty \mathcal{W}_{\ell^1}(\mu_h,\mu)
\\ & \leq \frac{1}{m}\sum_{k=b}^{b+m-1}  \max_i\|\nabla_i g\|_\infty  M e^{-ck} \mathcal{W}_{\ell^1}(\nu,\mu_h)+ \max_i\|\nabla_i g\|_\infty \mathcal{W}_{\ell^1}(\mu_h,\mu)
\\ &\leq  \frac{1}{m} \max_i\|\nabla_i g\|_\infty  M \frac{e^{-cb}}{1-e^{-c}} \mathcal{W}_{\ell^1}(\nu,\mu_h)+ \max_i\|\nabla_i g\|_\infty \mathcal{W}_{\ell^1}(\mu_h,\mu)
\end{align*}
with $M=\exp(\frac{5}{4}(2+\frac{\tilde{R}}{T}))$.
Applying \Cref{cor:bound_mu} yields the result.
\end{proof}

\begin{appendix}

\section{Contractivity of uHMC for $K$-Strongly Convex and $L$-gradient Lipschitz $V$} \label{appendixA}

Here, we consider the special case of a single particle with potential $V$ that is $K$-strongly convex and $L$-gradient Lipschitz.
In this case, we prove that the uHMC transition kernel is contractive with respect to the $L^p$-Wasserstein distance for $p\in[1,\infty)$, which is given by 
\begin{align*}
\mathcal{W}^p(\nu,\eta)=\inf_{\omega\in\Gamma(\nu,\eta)}\Big(\int|x-y|^p\omega(\rmd x\rmd y)\Big)^{1/p}
\end{align*}
for two probability measures $\nu,\eta$ on $\mathbb{R}^d$ with finite $p$-th moment, where $\Gamma(\nu,\eta)$ denotes the set of all couplings of $\nu$ and $\eta$.

\begin{theorem}[Contractivity of uHMC under global strong convexity] \label{thm:contr_uHMC}
Suppose that \Cref{ass_V_locmin}-\Cref{ass_V_strongconv} with $R=0$ hold. Let $T>0$ and $h \ge 0$ be such that \eqref{eq:CT} holds and $T/h \in \mathbb{Z}$ if $h>0$.  Then for any $p\in[1,\infty)$, any probability measures $\nu, \eta$ on $\mathbb{R}^d$ with finite $p$-th moment and $m \in \mathbb{N}_0$, \begin{align}
\mathcal{W}^p(\nu \tilde{\pi}^m, \eta \tilde{\pi}^m) \, &\le \, (1 - c)^m \, \mathcal{W}^p(\nu, \eta) \qquad \text{where} \\
c \, & = \, K \, T^2 \, / \, 10 \;.
\end{align}
\end{theorem}

For fixed duration hyperparameter, note that the $\mathcal{W}^p$ contraction rate $c$ is uniform in the timestep hyperparameter.  

\medskip

To prove this theorem, we introduce the following piecewise quadratic interpolation of the Verlet flow $(\tilde{q}_t(x,v), \tilde{v}_t(x,v))$  \begin{equation}
\label{verlet_flow}
	\frac{d}{dt} \tilde{q}_t   \, = \,   \tilde{v}_{\lb{t}}- (t-\lb{t})  \ \nabla
	V(\tilde{q}_{\lb{t}}) \;,  \quad  \frac{d}{dt} \tilde{v}_t \, = \, -\frac{1}{2}\left(\nabla
	V(\tilde{q}_{\lb{t}})\ +\ \nabla V(\tilde{q}_{\ub{t}})\right)
\end{equation}
with initial condition $(\tilde{q}_0(x,v),\tilde{v}_0(x,v)) = (x,v) \in \mathbb{R}^{2d}$.  The following lemma states that $|\tilde{q}_T(x,v) - \tilde{q}_T(y,v)|^2$ is itself contractive provided that the duration $T$ is sufficiently small as indicated, and $h \le T$ (which follows from $T/h \in \mathbb{Z}$).  This result extends the contractivity of the exact Hamiltonian flow from Lemma 2.1 of \cite{chenvempala19} to the velocity Verlet integrator.

\begin{lemma}[Contractivity of velocity Verlet under global strong convexity] \label{lem:contr_verlet_flow}
Suppose that \Cref{ass_V_locmin}-\Cref{ass_V_strongconv} with $R=0$ hold. Let $T>0$ and $h \ge 0$ satisfy $T/h \in \mathbb{Z}$ if $h>0$ and \begin{align} 
 L T^2 & \, \le \,  20^{-1} \label{eq:CT} \;.
\end{align}
Then for all $x, y, v \in \mathbb{R}^d$, \begin{equation}
|\tilde{q}_T(x,v) - \tilde{q}_T(y,v)|^2 \ \le \ \left( 1 - K \, T^2 \, / \, 5 \right) |x - y |^2 \;.
\end{equation}
\end{lemma}

\begin{proof}[Proof of Theorem~\ref{thm:contr_uHMC}]
By synchronously coupling the random initial velocities in two copies of uHMC and applying  Lemma~\ref{lem:contr_verlet_flow}, it immediately follows that the transition kernel of uHMC is contractive in the $L^p$-Wasserstein distance with respect to the Euclidean distance on $\mathbb{R}^d$ with the given contraction rate.
\end{proof}

\begin{remark} \label{rmk:cocoercivity}
If $V$ is continuously differentiable, convex, and $L$-gradient Lipschitz, then $\nabla V$ satisfies the following `co-coercivity' property \begin{equation} \label{eq:cocoercive}
|\nabla V(x) - \nabla V(y)|^2 \, \le \, L \, (\nabla V(x) - \nabla V(y)) \cdot (x- y) \;, \quad \text{for all $x, y \in \mathbb{R}^d$} \;. 
\end{equation}
This  property  plays a crucial role in proving Lemma~\ref{lem:contr_verlet_flow}.
\end{remark}

\begin{proof}[Proof of Lemma~\ref{lem:contr_verlet_flow}]
The proof parallels the proof of \Cref{lemma:contraction_num}, but employs the sharper argument from Lemma 2.1 of \cite{chenvempala19}.  Fix $t>0$ and $h \ge 0$ such that $t/h \in \mathbb{Z}$ for $h>0$. Introduce the shorthand $x_t = \tilde{q}_t(x,v)$ and $y_t = \tilde{q}_t(y,v)$.  Let $z_t := x_t - y_t$ and $w_t := v_t(x,v) - v_t(y,v)$.  Let $a_t \, :=\, \norm{z_t}^2$ and $b_t := 2 z_t \cdot w_t$. Our goal is to obtain an upper bound for $a_t$.  To this end, define \[
\rho_t \, := \,  \Phi_t \cdot (x_t - y_t)  \;, \qquad \Phi_t \, := \,  \nabla V(x_t) - \nabla V(y_t)  \;,  
\] 
and note by \Cref{ass_V_lipschitz}, \Cref{ass_V_strongconv} and \eqref{eq:cocoercive}, \begin{equation} \label{ieq:rho}
    K a_t \, \le 
    \, \rho_t  \, \le 
    \, L a_t \, \;, \qquad  |\Phi_t|^2 \, \le \, L \rho_t  \;, \qquad \text{for all $t \ge 0$} \;.
\end{equation} 
Moreover, by \eqref{verlet_flow}, note that \begin{align} 
		\frac{d}{dt} z_t \, &= \, w_{\lb{t}} - (t-\lb{t})  \Phi_{\lb{t}} ,\label{eq:dotz} \\
		~ \frac{d}{dt} w_t \, &= \, - \frac{1}{2} \left( \Phi_{\lb{t}} +\Phi_{\ub{t}} \right). \label{eq:dotw}
\end{align}
Let $\alpha>0$ be a parameter, which we specify shortly.
 A straightforward computation shows that 
	\begin{align} \label{ivp_a}
		\frac{d}{dt} a_t & \ = \ b_t + \delta_t \;, \\  \label{ivp_b} \frac{d}{dt} b_t & \ = \  \ - \alpha K a_t  + \beta_t \;,
		\end{align}
		where we have introduced \begin{align*}
		   \delta_t \, &:= \, 2 z_t \cdot ( w_{\lb{t}} - w_t - (t-\lb{t}) \Phi_{\lb{t}} ) = (t-\lb{t}) z_t \cdot (\Phi_{\ub{t}} - \Phi_{\lb{t}}) \;, \\ 
		   \beta_t  \, &:= \, \alpha K a_t + 2  | w_{\lb{t}} - (t-\lb{t}) \Phi_{\lb{t}} |^2 - z_t \cdot (\Phi_{\lb{t}} + \Phi_{\ub{t}}) \\
		   & \qquad -  (t-\lb{t}) \left( w_{\lb{t}} - (t-\lb{t}) \Phi_{\lb{t}} \right) \cdot (\Phi_{\ub{t}} - \Phi_{\lb{t}}) \;.
		\end{align*}  Note that $\delta_t$ is piecewise smooth satisfying \[
		\delta_t' \, = \, z_t \cdot (\Phi_{\ub{t}} - \Phi_{\lb{t}}) + (t-\lb{t}) \left( w_{\lb{t}} - (t-\lb{t}) \Phi_{\lb{t}} \right) \cdot (\Phi_{\ub{t}} - \Phi_{\lb{t}})
		\]
		between consecutive grid points, and having jump discontinuities at the grid points where $\delta_{t_k+} = 0$ and \begin{equation} \label{eq:delta_tk_minus}
	 \delta_{t_k-} = h ( \rho_{t_k} - \rho_{t_{k-1}} ) - h \left( h w_{t_{k-1}} - (h^2/2)  \Phi_{t_{k-1}} \right) \cdot \Phi_{t_{k-1}}  \;.
		\end{equation} Set $s_{t-r} \, := \, \sin( \sqrt{\alpha K} (t-r) ) / \sqrt{\alpha K}$ and $c_{t-r} \, := \,  \cos(\sqrt{\alpha K} (t-r) )$ such that $c_{t-r} = - \frac{\rmd}{\rmd r} s_{t-r}$ and \begin{equation} 
0 \le s_{t-r} \le s_{t-s} \quad  \text{for $s \le r \le t \le 1/(2 \sqrt{K})$} \;. \label{eq:s_t}
\end{equation}  By, first, variation of parameters, and second, integration by parts for piecewise smooth functions with jump discontinuities on the evenly spaced time grid $\{t_k\}$, 
\begin{align} 
    a_t & \, = \, c_t a_0 + \int_0^t c_{t-r} \delta_r dr + \int_0^t s_{t-r} \beta_r dr \nonumber \\
    & \, = \, c_t a_0 
     + \sum_{k:~t_k \le t} \left[ s_{t-r} \delta_r \right]_{r = t_{k}-}^{r = t_{k}+} 
    + \left[ -s_{t-r} \delta_r \right]_{r=0+}^{r=t-}
    + \int_0^t s_{t-r} (\delta'_r + \beta_r )  dr \, \nonumber \\
    & \, = \, c_t a_0 - \sum_{k:~t_k \le t} s_{t-t_k} \delta_{t_k-} + \int_0^t s_{t-r} (-2 \rho_{\lb{r}} + \alpha K a_{\lb{r}} + \epsilon_r  )  dr \label{vop_at} 
    \end{align}
where $\epsilon_t \, := \, \epsilon^1_t + \epsilon^2_t + \epsilon^3_t$ and
		\begin{align*}
		\epsilon^1_t \, &:= \,  2 |w_{\lb{t}} - (t-\lb{t}) \Phi_{\lb{t}}|^2 , \\
		\epsilon^2_t  \, &:= \,- 2 (z_t - z_{\lb{t}}) \cdot
		\Phi_{\lb{t}} , \\ \epsilon^3_t  \, &:= \, \alpha K (z_t - z_{\lb{t}}) \cdot (z_t + z_{\lb{t}}) .
	\end{align*} 
To upper bound $\epsilon_t^1$,  apply the Peter-Paul inequality with parameter $\sqrt{2}$, \begin{align}
    \epsilon^1_t & \, \le \, 6 |w_{\lb{t}}|^2 + 3 (t- \lb{t})^2 | \Phi_{\lb{t}}|^2 \, \overset{\eqref{ieq:rho}}{\le} \, 6 |w_{\lb{t}}|^2 + 3 L h^2 \rho_{\lb{t}} \;.  \label{ieq:epsi1} 
    \end{align}
Similarly, for $\epsilon_t^2$, apply  \eqref{eq:dotz} and Peter-Paul inequality with parameter $\sqrt{2}$,
\begin{align}
    \epsilon^2_t & \, = \, -2 (t-\lb{t}) w_{\lb{t}} \cdot \Phi_{\lb{t}} +  (t-\lb{t})^2 |\Phi_{\lb{t}}|^2 \overset{\eqref{ieq:rho}}{\le} 2 |w_{\lb{t}}|^2 + \frac{3}{2} L h^2  \rho_{\lb{t}}  \;.  \label{ieq:epsi2} 
    \end{align}
Finally, for $\epsilon^3_t$, apply \eqref{eq:dotz} and Young's product inequality 
    \begin{align}
    \epsilon^3_t & \, = \, \alpha K \left( |z_t - z_{\lb{t}}|^2 + 2 z_{\lb{t}} \cdot (z_t - z_{\lb{t}}) \right) \nonumber \\
   &\, = \, \alpha K  | (t-\lb{t}) w_{\lb{t}}  -  (1/2) (t-\lb{t})^2 \Phi_{\lb{t}}  |^2  \nonumber \\
 & \qquad  + 2 \alpha K z_{\lb{t}} \cdot ( (t-\lb{t}) w_{\lb{t}} - (1/2)  (t-\lb{t})^2 \Phi_{\lb{t}} )  \nonumber \\
&\, \overset{\eqref{ieq:rho}}{\le} \, (1 + 2 \alpha K h^2) |w_{\lb{t}}|^2 + \alpha K ( \alpha h^2 + (1/2) L h^4  ) \rho_{\lb{t}} \;. \label{ieq:epsi3}
\end{align}
To upper bound the sum in \eqref{vop_at} coming from integration by parts, expand the sum using \eqref{eq:delta_tk_minus}, apply summation by parts, and Young's product inequality \begin{align}
  & - \sum_{k:~t_k < t} s_{t-t_k} \delta_{t_k-}  \ = \ 
   - \sum_{k:~t_1 \le t_k < t} h s_{t-t_k} (\rho_{t_k} - \rho_{t_{k-1}}) \nonumber \\
   & \qquad + \sum_{k :~t_1 \le t_k < t}  h s_{t-t_k} (h w_{t_{k-1}} - \frac{h^2}{2} \Phi_{t_{k-1}} ) \cdot \Phi_{t_{k-1}}  \nonumber \\
   & \ = \     \sum_{k :~ t_2 \le t_k < t} h ( s_{t-t_k} - s_{t-t_{k-1}} ) \rho_{t_{k-1}} - h s_{t-\lb{t}} \rho_{\lb{t}} + h s_{t- t_1} \rho_{0} \nonumber \\
   & \qquad + \sum_{k:~t_1 \le t_k < t}  h s_{t-t_k} (h w_{t_{k-1}} - \frac{h^2}{2} \Phi_{t_{k-1}} ) \cdot \Phi_{t_{k-1}} \nonumber \\
   & \ \overset{\eqref{ieq:rho}}{\le} \  \int_0^{t_1} s_{t- r} \rho_{\lb{r}} dr + (1/2) \int_0^t s_{t-r} ( |w_{\lb{r}}|^2 +  L h^2 \rho_{\lb{r}} ) dr  \label{ieq:bdry}
\end{align}
where in the last step we used Young's product inequality and \eqref{eq:s_t}.  
To estimate the terms in \eqref{ieq:epsi1}-\eqref{ieq:bdry} involving $|w_{\lb{t}}|^2$, by \eqref{eq:dotw} and since $w_0=0$,  \begin{align}
&  2 |w_{\lb{t}}|^2 \, = \,  2 \left| \frac{1}{2} \int_0^{\lb{t}}  \left( \Phi_{\lb{s}} + \Phi_{\ub{s}} \right) ds \right|^2 \nonumber  ~ \, \le \, 	\left| \int_0^{\lb{t}}  \Phi_{\lb{s}}   ds \right|^2 + 	\left| \int_0^{\lb{t}} \Phi_{\ub{s}}  ds \right|^2 \nonumber  \\
& \quad \, \le \, t \left( \int_0^{\lb{t}} | \Phi_{\lb{s}} |^2 ds  +  \int_0^{\lb{t}} | \Phi_{\ub{s}} |^2  ds  \right)
		\, \overset{\eqref{ieq:rho}}{\le} \, 2 L t \int_0^{t}  \rho_{\lb{s}} ds 	
	\nonumber 
\end{align}
where in the second to last step we used Cauchy-Schwarz inequality.  By \eqref{eq:s_t} and Fubini's Theorem,  \begin{align}
\int_0^t s_{t-r} |w_{\lb{r}}|^2 dr & \,  \overset{\eqref{eq:s_t}}{\le} \,
 L \int_0^t \int_0^r  r s_{t-s}   \rho_{\lb{s}} ds dr
=  L \int_0^t \int_s^t  r s_{t-s}  \rho_{\lb{s}}  dr ds  \nonumber \\
& \ \le \,  \frac{L t^2}{2} \int_0^t   s_{t-s} \rho_{\lb{s}}  ds \;.  \label{ieq:wt2}
\end{align} Combining \eqref{ieq:epsi1}-\eqref{ieq:bdry},  \begin{align}
&  \sum_{k:~t_k < t} [-s_{t-t_k} \delta_{t_k-}]  + \int_0^t s_{t-r} (-2 \rho_{\lb{r}} + \alpha K a_{\lb{r}} + \epsilon_r  )  dr \nonumber \\
&   \, \le \,  \int_0^{t_1} s_{t-r}  \rho_{\lb{r}} dr + 
\int_0^t s_{t-r} ( - 2 \rho_{\lb{r}} +  \alpha K a_{\lb{r}}  ) dr \nonumber \\ 
& \qquad +\int_0^t s_{t-r} ( (5 L h^2 + \alpha^2 K h^2 + \frac{\alpha}{2} K L h^4)  \rho_{\lb{r}} + ( \frac{19}{2} + 2 \alpha K h^2) |w_{\lb{r}}|^2 ) dr \nonumber  \\
& ~~  \, \overset{\eqref{ieq:wt2}}{\le} \,  \int_0^{t_1} s_{t-r}  \rho_{\lb{r}} dr + 
\int_0^t s_{t-r} ( - 2 \rho_{\lb{r}} +  \alpha K a_{\lb{r}}  ) dr \nonumber \\ 
&  \qquad +\int_0^t s_{t-r} ( (5 L h^2 + \alpha^2 K h^2 + \frac{\alpha}{2} K L h^4)   + ( \frac{19}{4} + \alpha K h^2) L t^2 ) \rho_{\lb{r}} dr \nonumber  \\
& ~~  \, \overset{\eqref{ieq:rho}}{\le} \,   \int_{0}^t s_{t-r}  \left[ \alpha + ( \frac{39}{4} + \alpha^2 + \frac{3}{2} \alpha L t^2) L t^2  - 1 \right] \rho_{\lb{r}}  dr \le 0 \quad \text{with $\alpha=4/9$} \label{eq:nonposrem} 
\end{align}
where in the last step we used $t/h \in \mathbb{Z}$, $K \le L$ and condition~\eqref{eq:CT} (i.e.,  $L t^2 \le 20^{-1}$). The required estimate is then obtained by inserting \eqref{eq:nonposrem} into \eqref{vop_at} and then  using the elementary inequality \begin{align*}
c_t \, &\le \, 1 - (\alpha/2) K t^2 + (1/6) \alpha^2 K^2 t^4 \le 1 - (\alpha/2) K t^2 + (1/120) \alpha^2 K t^2 (K/L)\\  
 \, &\le \, 1 - (\alpha/2 - \alpha^2/120) K t^2 \le 1 -  K t^2 / 5  \quad \text{with $\alpha=4/9$} 
\end{align*} which follows from  condition~\eqref{eq:CT} (i.e., $L t^2 \le 20^{-1}$) and $K \le L$. \end{proof}

 \section{Perturbation of the product model}\label{appendix}
If the confinement potential is a quadratic potential, i.e., $V(\mathsf{x})=K/2|\mathsf{x}|^2$ for all $\mathsf{x}\in\mathbb{R}^d$, the mean-field model can be treated as a perturbation of the product model.
Given $x,y\in\mathbb{R}^{dn}$ we consider the synchronous coupling of four transition kernels $\pi_h(x,\cdot)$, $\pi_h(y,\cdot)$, $\pi_h^{prod}(x,\cdot)$ and $\pi_h^{prod}(y,\cdot)$, where $\pi_h(x,\cdot)$ and $\pi_h(y,\cdot)$ denote the two transition kernels with a mean-field interaction, i.e., $\epsilon>0$, and $\pi_h^{prod}(x,\cdot)$ and $\pi_h^{prod}(y,\cdot)$ are transition kernels of the product model, i.e., $\epsilon=0$.
Then the coupling HMC step is given by
\begin{align*}
&\mathbf{X}(x,y)=q_T(x,\xi), && \mathbf{Y}(x,y)=q_T(y,\xi), \\& \mathbf{X}^{prod}(x,y)=\hat{q}_T(x,\xi), && \mathbf{Y}^{prod}(x,y)=\hat{q}_T(y,\xi),
\end{align*}
where $\xi\sim\mathcal{N}(0,I_{dn})$ and $\hat{q}_T$ denotes the position component of the Hamiltonian dynamics given by \eqref{eq:hamdyn_num} for the product model. 
\begin{theorem} \label{thm:app}
Suppose that $V(\mathsf{x})=(K/2)|\mathsf{x}|^2$ for all $\mathsf{x}\in\mathbb{R}^d$ and \Cref{ass_W_lipschitz} hold. Let $T\in(0,\infty)$, $h_1\in[0,\infty)$ and $\epsilon\in(0,\infty)$ satisfy
\begin{align} \label{eq:app_condT}
K(T^2+Th_1) \le 1.
\end{align}
Then for any $h\in[0,h_1]$ such that $h=0$ or $T/h\in\mathbb{N}$ and any $x,y\in\mathbb{R}^{dn}$,
\begin{align*}
\sum_{i=1}^n|\mathbf{X}^i(x,y)-\mathbf{Y}^i(x,y)&-(\mathbf{X}^{i,prod}(x,y)-\mathbf{Y}^{i,prod}(x,y))|
 \le 8\epsilon \tilde{L}(T^2+Th)\sum_{i=1}^n |x^i-y^i|.
\end{align*}
\end{theorem}

\begin{proof}
Fix $x,y,v\in\mathbb{R}^d$. For $t\in[0,T]$, we write $x_t^i=q_t^i(x,v)$ and $y_t^i=q_t^i(y,v)$ for the $i$-th position component of the solution to \eqref{eq:hamdyn_num} with initial values $(x,v)$ and $(y,v)$, respectively, and with potential $U(x)=\sum_{i=1}^n((K/2)|x^i|+\epsilon n^{-1} \sum_{j=1, j\neq i}^n W(x^i-x^j))$. Analogously, we write $\hat{x}_t^i=\hat{q}_t^i(x,v)$ and $\hat{y}_t^i=\hat{q}_t^i(y,v)$ for the $i$-th position component of the solution to \eqref{eq:hamdyn_num} with initial values $(x,v)$ and $(y,v)$, respectively, and with potential $\hat{U}(x)=\sum_{i=1}^n(K/2)|x^i|$. We set $z_t^i=x_t^i-y_t^i$ and $\hat{z}_t^i=\hat{x}_t^i-\hat{y}_t^i$ for all $i=1,\ldots, n$ and $t\in[0,T]$.
By \eqref{eq:hamdyn_num} and \Cref{ass_W_lipschitz} it holds for $t\in[0,T]$,
\begin{align}
\max_{s \le t} \sum_{i=1}^n&|z_s^i-\hat{z}_s^i|  \nonumber
\\ &= \max_{s\le t} \sum_{i=1}^n\Big|\int_0^s \int_0^r \Big(-\frac{1}{2}(\nabla_i U(x_{ \lfloor u \rfloor}^i)-\nabla_i U(y_{ \lfloor u \rfloor}^i)+\nabla_i U(x_{ \lceil u \rceil}^i)-\nabla_i U(y_{ \lceil u \rceil}^i)) \nonumber
\\ & +\frac{1}{2}(\nabla_i \hat{U}(\hat{x}_{ \lfloor u \rfloor}^i)-\nabla_i \hat{U}(\hat{y}_{ \lfloor u \rfloor}^i) +\nabla_i \hat{U}(\hat{x}_{ \lceil u \rceil}^i)-\nabla_i \hat{U}(\hat{y}_{ \lceil u \rceil}^i))\Big)\rmd u \rmd r  \nonumber
\\ & -\frac{h}{2}\int_0^s  \Big(\nabla_i U(x_{ \lfloor u \rfloor}^i)-\nabla_i U(y_{ \lfloor u \rfloor}^i)-(\nabla_i \hat{U}(\hat{x}_{ \lfloor u \rfloor}^i)-\nabla_i \hat{U}(\hat{y}_{ \lfloor u \rfloor}^i))\Big) \rmd u\Big| \nonumber
\\ & \le \frac{K}{2}(t^2+th)\max_{s\le t} \sum_{i=1}^n|z_s^i-\hat{z}_s^i|+ 2\epsilon\tilde{L}(t^2+th) \max_{s\le t } \sum_{i=1}^n|z_s^i|. \label{eq:app_estimate}
\end{align}
By \eqref{eq:app_condT} and \eqref{eq:position_estimate6},
\begin{align*}
\max_{s\le t} \sum_{i=1}^n |z_s^i-\hat{z}_s^i|\le 4\epsilon \tilde{L}(t^2+th)\max_{s\le t} \sum_{i=1}^n |z_s^i|\le 8 \epsilon \tilde{L}(t^2+th) \sum_{i=1}^n |x^i-y^i|.
\end{align*}
Thus, the result holds for $t=T$.
\end{proof}
We note that the step \eqref{eq:app_estimate} uses crucially that the third derivative of $V$ vanishes. 

As some calculations simplify in the product case with quadratic confinement potential, \eqref{eq:contraction_num_bound} in \Cref{lemma:contraction_num} holds for all $i=1,\ldots,n$ provided $K(t^2+th)\le 1/4$ and $h\le (4/165)t$ is satisfied.
Hence by \eqref{eq:hilfungleichung},
\begin{align*}
\sum_{i=1}^n|\mathbf{X}^{i,prod}(x,y)-\mathbf{Y}^{i,prod}(x,y)|\le (1-(1/8)KT^2)\sum_{i=1}^n|x^i-y^i|
\end{align*}
for $K(T^2+Th)\le 1/4$ and $h\le (4/165)T$.
Combining the contraction result for the product model with the  perturbation result yields the following consequence.

\begin{korollar} \label{kor:app}
Suppose that $V(\mathsf{x})=(K/2)|\mathsf{x}|^2$ for all $\mathsf{x}\in\mathbb{R}^d$ and \Cref{ass_W_lipschitz} hold. Let $T\in(0,\infty)$, $h_1\in(0,\infty)$ and $\epsilon\in(0,\infty)$ satisfy
\begin{align}
&K(T^2+Th_1)\le 1/4, \qquad h\le (4/165)T, \qquad \text{and} \nonumber
\\ & \epsilon\tilde{L}\le K/256. \label{eq:app_epsilon}
\end{align}
Then, for any $h\in[0,h_1]$ such that $h=0$ or $T/h\in\mathbb{N}$ and for any $x,y\in\mathbb{R}^{dn}$
\begin{align*}
\sum_{i=1}^n|\mathbf{X}^i(x,y)-\mathbf{Y}^i(x,y)|\le (1-KT^2/16)\sum_{i=1}^n |x^i-y^i|,
\end{align*}
and for any two probability measures $\nu$ and $\eta$ on $\mathbb{R}^{dn}$ and any $m\in\mathbb{N}$, 
\begin{align*}
\mathcal{W}_{l^1} (\nu \pi_h^m, \eta \pi_h^m)\le e^{-KT^2m/16} \mathcal{W}_{l^1} (\nu, \eta).
\end{align*}
\end{korollar}

\begin{proof}
The result is a direct consequence of the contraction result and \Cref{thm:app}, i.e.,
\begin{align*}
\sum_{i=1}^n|\mathbf{X}^i(x,y)-\mathbf{Y}^i(x,y)|&\le \sum_{i=1}^n|\mathbf{X}^{i,prod}(x,y)-\mathbf{Y}^{i,prod}(x,y)|
\\ & +\sum_{i=1}^n|\mathbf{X}^i(x,y)-\mathbf{Y}^i(x,y)-(\mathbf{X}^{i,prod}(x,y)-\mathbf{Y}^{i,prod}(x,y))|
\\ & \le (1-KT^2/8)\sum_{i=1}^n |x^i-y^i|+(8\epsilon\tilde{L}(T^2+Th)\sum_{i=1}^n |x^i-y^i|
\\ & \le  (1-KT^2/16)\sum_{i=1}^n |x^i-y^i|,
\end{align*}
where the last step follows by \eqref{eq:app_epsilon}. The second bound in \Cref{kor:app} holds in the same line as the proof of \Cref{cor_Wassersteindist}.
\end{proof}

\end{appendix}

 \section*{Acknowledgments}
 The authors would like to thank Andreas Eberle for his insights and advice during the development of this work.
 
N.~B.-R. was supported by the National Science Foundation under Grant No.~DMS-1816378 and the Alexander von Humboldt Foundation.

K.~S. was supported by \textit{Bonn International Graduate School of Mathematics}.
 Gef\"ordert durch die Deutsche Forschungsgemeinschaft (DFG) im Rahmen der Exzellenzstrategie des Bundes und der L\"ander - GZ 2047/1, Projekt-ID 390685813.

%
%
 


\bibliographystyle{imsart-number} 
\bibliography{bibliography}       

%

\end{document}